\documentclass[reqno]{amsart}

\usepackage{amsmath}
\usepackage{amscd}
\usepackage{amssymb}
\usepackage{enumerate}
\usepackage{amsfonts}
\usepackage{graphicx}

\DeclareFontEncoding{OT2}{}{}

\numberwithin{equation}{section}

\newcommand{\rar}[1]{\stackrel{#1}{\longrightarrow}}

\newcommand{\al}{\alpha} \newcommand{\be}{\beta} \newcommand{\ga}{\gamma}
 \newcommand{\de}{\delta} \newcommand{\De}{\Delta}
\newcommand{\la}{\lambda} \newcommand{\La}{\Lambda}
 \newcommand{\ka}{\kappa}
\newcommand{\eps}{\epsilon} \newcommand{\sg}{\sigma}
  
\newcommand{\om}{\omega} \newcommand{\Om}{\Omega}

\newcommand{\bC}{{\mathbb C}}

\newcommand{\bN}{{\mathbb N}}

\newcommand{\bQ}{{\mathbb Q}}
\newcommand{\bR}{{\mathbb R}}

\newcommand{\bZ}{{\mathbb Z}}

\newcommand{\cH}{{\mathcal H}}

\newcommand{\cK}{{\mathcal K}}
\newcommand{\cL}{{\mathcal L}}
\newcommand{\cM}{{\mathcal M}}

\newcommand{\cU}{{\mathcal U}}

\newcommand{\fC}{{\mathfrak C}}

\newcommand{\fa}{{\mathfrak a}}

\newcommand{\fg}{{\mathfrak g}}
\newcommand{\fh}{{\mathfrak h}}

\newcommand{\fo}{{\mathfrak o}}

\newcommand{\Bb}{\overline{B}}

\newcommand{\Gb}{\overline{G}}

\newcommand{\Lb}{\overline{L}}

\newcommand{\Xb}{\overline{X}}
\newcommand{\Yb}{\overline{Y}}

\newcommand{\fgb}{\overline{\fg}}

\newcommand{\dd}{\partial}

\newcommand{\eq}[1]{(\ref{#1})}

\newcommand{\abs}[1]{\vert #1\vert}
\newcommand{\norm}[1]{\vert\vert #1\vert\vert}

\newcommand{\Ker}{\operatorname{Ker}}

\newcommand{\pr}{\mathrm{pr}}
\newcommand{\ad}{\operatorname{ad}}
\newcommand{\Ad}{\operatorname{Ad}}

\newcommand{\Span}{\operatorname{span}}
\newcommand{\Supp}{\operatorname{supp}}
\newcommand{\supp}{\operatorname{supp}}
\newcommand{\meas}{\operatorname{meas}}
\newcommand{\rank}{\operatorname{rank}}
\newcommand{\Ind}{\operatorname{Ind}}

\newcommand{\st}{\,\big\vert\,}

\newcommand{\sbr}{\smallbreak}
\newcommand{\mbr}{\medbreak}

\newcommand{\pare}[1]{\bigl( #1\bigr)}
\newcommand{\brak}[1]{\bigl[ #1\bigr]}
\newcommand{\cur}[1]{\bigl\{ #1\bigr\}}

\newcommand{\bra}{\langle}
\newcommand{\ket}{\rangle}

\newcommand{\pair}[2]{\bra #1, #2\ket}

\newcommand{\matr}[4]{\left(\begin{array}{cc} \! #1 & \! #2 \! \\ \! #3 & \! #4 \! \end{array}\right)}

\newcommand{\sch}{Schr\"odinger\ }

\newcommand{\Uj}{U_j}
\newcommand{\Uoj}{U_{1j}}
\newcommand{\xj}{x^j}

\newcommand{\Pso}{\Psi_0}
\newcommand{\Pst}{\Psi^*}

\newcommand{\Psti}{\tilde{\Psi}}
\newcommand{\dti}{\tilde{\de}}

\newcommand{\lan}{\langle}
\newcommand{\ran}{\rangle}

\newcommand{\vone}{\textbf{1}}

\newtheorem{thm}{Theorem}[section]
\newtheorem{cor}[thm]{Corollary}
\newtheorem{lem}[thm]{Lemma}
\newtheorem{prop}[thm]{Proposition}

\theoremstyle{remark}
\newtheorem{rem}[thm]{Remark}

\newtheorem{defin}[thm]{Definition}

\newcommand{\lemm}[1]{Lemma \ref{#1}}
\newcommand{\theor}[1]{Theorem \ref{#1}}

\newcommand{\sect}[1]{Section \ref{#1}}
\newcommand{\subsect}[1]{\S\ref{#1}}

\newcommand{\rn}{\bR^n}

\newcommand{\unirrep}{unitary irreducible representation\ }
\newcommand{\unirreps}{unitary irreducible representations\ }

\title[Beyond Weyl and Colin de Verdi\`ere's formulas]{Beyond the classical Weyl
and Colin de Verdi\`ere's formulas \\ for Schr\"odinger operators
with polynomial \\ magnetic and electric fields}

\author[M.~Boyarchenko and S.Z.~Levendorskii]{Mitya Boyarchenko\address{Mitya Boyarchenko: Department of
Mathematics, University of Chicago, Chicago, IL 60637, \\ e-mail:
mitya@math.uchicago.edu} \and Sergei Levendorskii\address{Sergei
Levendorskii: Department of Economics, University of Texas,
Austin, TX, \\ e-mail:leven@eco.utexas.edu}}

\date{May 21, 2004}

\begin{document}

\begin{abstract}
We present a pair of conjectural formulas that compute the leading
term of the spectral asymptotics of a Schr\"odinger operator on
$L^2(\bR^n)$ with quasi-homogeneous polynomial magnetic and
electric fields. The construction is based on the orbit method due
to Kirillov. It makes sense for any nilpotent Lie algebra and is
related to the geometry of coadjoint orbits, as well as to the
growth properties of certain ``algebraic integrals,'' studied by
Nilsson. By using the direct variational method, we prove that the
formulas give the correct answer not only in the ``regular'' cases
where the classical formulas of Weyl or Colin de Verdi\`ere are
applicable but in many ``irregular'' cases, with different types
of degeneration of potentials.
\end{abstract}

\maketitle

\setcounter{tocdepth}{1}

\tableofcontents

\section{Introduction}\label{s:intro}

\subsection{}  Let
\begin{equation}\label{e:schr}
H =H(a) + V = -\sum_{j=1}^n \left(\frac{\dd}{\dd x_j} +
\sqrt{-1}\cdot a_j(x)\right)^2 + V(x)
\end{equation}
be a Schr\"odinger operator in $\rn\ $with a real semi-bounded
electric potential $V $ and  magnetic potential $
a(x)=(a_{1}(x),\ldots , a_{n}(x)) \in C^{2}(\rn , \rn ).$ If the
electric potential $V$ regularly grows at infinity, it is
well-known that $H(0)+V\ $is a self-adjoint operator with discrete
spectrum, and the counting function of the spectrum obeys the
classical Weyl formula
\begin{equation}\label{e:weyl1}
 N(\la , H(0)+ V)\sim (2\pi )^{-n}\int_{a(x, \xi)<\la}dxd\xi.
 \end{equation}
Here $a(x, \xi)=||\xi||^2+V(x)$ is the symbol of $H(0)+V$, and
$f(\la)\sim g(\la)$ means that $f(\la)/g(\la)\to 1$ as $\la\to
+\infty$. One
 easily rewrites \eq{e:weyl1} in the form
 \begin{equation}\label{e:weyl2}
 N(\la , H(0)+ V)\sim (2\pi )^{-n}|v_{n}|\int_{\rn }
(\la - V(x))_{+}^{n/2} dx,
\end{equation}
where $|v_{n}| $ is the volume of the unit ball of $\rn $ and
$a_{+}=\max\{ 0, a\}\  $(see e.g. \cite{Roz, RSS}).

In \cite{Sim, R, S, L1, G}, it is shown that the spectrum of
$H(0)+V$ can be discrete even if $V$ does not grow in some
directions, and for wide classes of degenerate potentials, the
leading term of the asymptotics of $N(\la,H(a)+V)$ is computed.
The results of these papers agree with the general "uncertainty
principle" stated in \cite{Fe}; it seems that this principle
provides upper and lower bounds but it is difficult to use it to
study spectral asymptotics. Note that in many cases, asymptotic
formulas are non-classical in the sense that they do not agree
with the "classical" formula \eq{e:weyl1}. Three cases are
possible: the classical Weyl formula holds (the so-called {\em
weak degeneration case}); an analog of the classical Weyl formula
with the operator-valued symbol parameterized by points of a set
with a measure inherited from $T^*\bR^{n}$ is valid ({\em strong
degeneration case}); the classical Weyl formula fails but the
leading term of the asymptotics is expressed in terms of an
auxiliary scalar function, and no operator-valued symbol is
involved ({\em intermediate degeneration case}). In simple strong
degeneration  cases, an operator valued symbol is parameterized by
the cotangent bundle over a manifold of degeneration of $V$, call
it $M$, and the operator-valued analog of \eq{e:weyl1} is of the
form
\begin{equation}\label{e:weyl3}
 N(\la , H(0)+ V)\sim (2\pi )^{-n+r}\int_{T^*M}N(\la, a(x',
 \xi'))dx'd\xi',
 \end{equation}
 where $r={\rm codim}M$, and for each $(x', \xi')\in T^*M$, $a(x', \xi')$ is an operator
 in $L_2(\bR^r)$. Similar types of asymptotic formulas holds
 for many other classes of differential operators,
 pseudodifferential operators and boundary value problems (see
 \cite{L1,L2,L3,L4} and the bibliography therein).

\subsection{}
 If $V=0\ $and the magnetic tensor $
B=[b_{jk}],\quad
b_{jk}(x)=\partial_{k}a_{j}(x)-\partial_{j}a_{k}(x)$, grows
regularly at infinity, the leading term of the asymptotics was
obtained in \cite{CdV} (see also \cite{T, I}):
\begin{equation}\label{e:CdV}
N(\lambda , H(a))\sim\int_{\bR^n }v_{B(x)}(\lambda )dx,
\end{equation}
where $v_B(\la)$ is defined as follows. Let $\rank B=2r$, and let
$b_1\geq b_2\geq \dotsb\geq b_r>0$ be the positive eigenvalues of
$iB$. Then
\[
v_B(\la)=(2\pi)^{-n+r} \abs{v_{n-2r}} b_1\dotsm b_r
\sum_{n_1,\dotsc,n_r\geq 0} ( \la-\sum_{j=1}^r (2n_j+1)\cdot b_j
)_+^{n/2-r}.
\]
Note that $B$, $r$ and the $b_j$'s depend on $x$. However, in the
case of a \sch operator with polynomial potentials, there is a
dense open subset of $\bR^n$ of full measure on which $B(x)$ has
maximal rank, so one can replace the integral in \eqref{e:CdV} by
the integral over this subset. Then $r$ will remain constant
throughout the integration.

\subsection{}
In the general case, only upper and lower bounds for $N(\la,
H(a)+V)\ $ are known \cite{MN}. They are given in terms of a
function $\Pst= \Pst_{a,V}\ $constructed in \cite{HM}; for
polynomial $V(\geq 0)\ $ and $b_{jk} $,
\begin{equation}\label{e:psistar}
\Pst(x)=\sum_{\al}|\partial^{\al}V(x)|^{1/(|\al |+2)} +\sum_{\al
,j,k}|\partial^{\al}b_{jk}(x)|^{1/(|\al |+2)}.
\end{equation}
It was shown that under certain conditions - in particular, for
polynomial $V\ $and $b_{jk}\ $- the following statements hold:

(\cite{HM}) {\em the spectrum of $H(a)+V\ $is discrete if and only
if}
\begin{equation}\label{e:necsuf}
\Pst(x)\to +\infty\ as\ |x|\to +\infty,
\end{equation}
and (\cite{MN}) {\em if \eq{e:necsuf} holds then there exists
$C>0\ $ such that
\begin{equation}\label{e:est}
C^{-1}G(C^{-1}\lambda , \Pst, \rn )\leq N(\lambda , H(a)+V)\leq
CG(C\lambda , \Pst, \rn ),
\end{equation}
where}
\begin{equation}\label{e:defG}
G(\lambda , \Pst, U)=\int_{U}(\lambda -\Pst(x)^{2})^{n/2}_{+}dx.
\end{equation}
Thus in the case $B\neq 0\ $not growing in some directions, the
leading term of the asymptotics is unknown apart from a special
case of Schr\"odinger operator (and Dirac operator) in 2D with
homogeneous potentials \cite{L4}.

\subsection{}\label{ss:pushforward} Notice the difference between
the formulas \eq{e:weyl1}, \eq{e:weyl3} and \eq{e:CdV}: the first
two are written in an invariant form, whereas the last one is
similar to \eq{e:weyl2}, which is a realization of the invariant
formula \eq{e:weyl1}. This observation suggests that there should
be an invariant formula of which \eq{e:CdV} is a realization.
Moreover, one should expect that there is a general formula, with
\eq{e:weyl1}, \eq{e:weyl3} and \eq{e:CdV} as special cases. The
following observations indicate the direction where one should
look for such a formula.

\sbr

Let $H$ be a \sch operator \eqref{e:schr} with polynomial
potentials. We consider the real Lie algebra $\fg$ of polynomial
differential operators on $\bR^n$ generated by the operators
$L_j=\dd/\dd x_j+\sqrt{-1}\cdot a_j(x)$ and the operator $L_0$ of
multiplication by $\sqrt{-1}\cdot V(x)$. It is easy to see that
$\fg$ is a finite dimensional nilpotent Lie algebra. We prove in
Theorem \ref{t:reps} that, after possibly replacing $H$ with a
gauge equivalent operator (a process which changes neither the
spectrum of $H$ nor the Lie algebra $\fg$, up to isomorphism), the
``tautological'' representation of $\fg$ on $L^2(\bR^n)$ lifts to
a unitary representation $\rho$ of the corresponding connected and
simply connected nilpotent Lie group $G=\exp\fg$. Moreover, $H$
has discrete spectrum if and only if $\rho$ is irreducible. Note
that $H$ can be naturally viewed as the image of the element
$H^\circ:=-\sum_{j=1}^n L_j^2-\sqrt{-1}\cdot L_0$ of
$\cU(\fg)_\bC=\cU(\fg)\otimes_\bR \bC$, the complexification of
the universal enveloping algebra of $\fg$, under the
representation of $\fg$ induced by $\rho$.

\sbr

Assume that $\sg(H)=\sg_d(H)$, so $\rho$ is irreducible. The orbit
method, due to Kirillov \cite{Ki}, provides a natural one-to-one
correspondence between (unitary equivalence classes of) unitary
irreducible representations of $G$ and orbits of the coadjoint
action of $G$ on $\fg^*$. In particular, we let
$\Om_\rho\subset\fg^*$ denote the coadjoint orbit corresponding to
$\rho$. Suppose now that the magnetic potential $a=0$, and that
$V(x)$ grows regularly at infinity. The values of the symbol
$a(x,\xi)$ appearing in the classical Weyl formula \eqref{e:weyl1}
can be interpreted as the images of $H^\circ$ in a family of
representations of $G$ on the one-dimensional space $L^2(\bR^0)$.
The family is parameterized by points of the orbit $\Om_\rho$, and
the measure $(2\pi)^{-n} dx d\xi$ coincides with the canonical
(Kostant) measure on $\Om_\rho$.

\sbr

On the other hand, assume that $V=0$ and the magnetic tensor
$B(x)$ grows regularly at infinity. It can (and will) be shown
that the formula of Colin de Verdi\`ere \eq{e:CdV} can be written
in the form
\begin{equation}\label{e:sd0}
N(\la, H)\sim \int_Q N(\la,H_\Theta) \,d\nu(\Theta) \quad
\text{as}\ \la\to +\infty,
\end{equation}
where $H_\Theta$ is the image of $H^\circ$ in a certain unitary
irreducible representation of $G$ on $L^2(\bR^r)$, $Q$ is a
manifold parameterizing a family of such representations, and the
measure $d\nu(\Theta)$ can be obtained in the following way. Let
$\tilde{Q}\subset\fg^*$ be the union of the orbits corresponding
to the representations parameterized by the points of $Q$. There
is a natural ``projection map'' $p:\Om_\rho\to\tilde{Q}$, such
that the pushforward $\tilde{\nu}$ of the canonical measure on
$\Om_\rho$ is a $G$-invariant measure on $\tilde{Q}$. One can
decompose $\tilde{\nu}$ as an integral of the canonical measures
on the orbits contained in $\tilde{Q}$, with respect to a certain
``quotient'' measure on $Q=\tilde{Q}/G$. Then we take $\nu$ to be
this quotient measure.

\sbr

The classical Weyl formula also can be written in the form
\eqref{e:sd0}, with $Q$ parameterizing a family of one-dimensional
representations (in this case, $Q=\tilde{Q}$, so one does not need
to decompose the pushforward measure). It is tempting to
conjecture that for any magnetic Schr\"odinger operator with
discrete spectrum one can find a family of irreducible
representations of $G$ and the pushforward measure $d\nu(\Theta)$
on $Q$ such that \eq{e:sd0} holds. As it turns out, this
construction can be realized in many cases albeit not all, and the
first goal of the paper is to suggest a general way of
construction of the family $Q$ and the pushforward measure
$d\nu(\Theta)$. Naturally enough (cf. generalizations of the
classical Weyl formula in \cite{L1, L2, L3}), we have two similar
(but a bit different) algorithms: one for the strong degeneration
case, and one for the weak and intermediate degeneration case; and
in the intermediate degeneration case, one has to introduce
additional logarithmic factors into \eq{e:sd0}. To verify our
conjecture for several classes of magnetic Schr\"odinger
operators, we use a modification of the variational technique from
\cite{L1, L2, L3, L4}.

\subsection{}\label{ss:scaling} Let us keep the same notation as above, and write
$\mu_{\Om_\rho}$ for the canonical (Kostant) measure on the orbit
$\Om_\rho$. In trying to turn the vague ideas of the previous
paragraph into a precise formula that applies to \sch operators
with degenerate potentials, one meets two considerable
difficulties. The first one lies in the fact that there seems to
be no natural general way of defining a projection map
$p:\Om_\rho\to\tilde{Q}\subset\fg^*$ such that the pushforward
$p_*(\mu_{\Om_\rho})$ will always be a $G$-invariant measure. The
second difficulty, which is more serious, is that in certain cases
(such as the intermediate degeneration example studied in
\S\ref{ss:intermed}) there exists an asymptotic formula of the
form \eqref{e:sd0} (with additional logarithmic factors), but the
measure $\nu$ cannot be obtained from a {\em pushforward} measure
arising from a process described above.

\sbr

Thus, one has to look for a different construction of the subset
$\tilde{Q}\subset\fg^*$ and the $G$-invariant measure
$\tilde{\nu}$ on $\tilde{Q}$. In our paper we suggest a
construction which has the advantage of being canonical (i.e.,
independent of any choices). Moreover, the measure $\tilde{\nu}$
it provides is automatically $G$-invariant. Thus, both problems
mentioned in the last paragraph are solved at once. To the best of
our knowledge, no similar construction has been used in this or
any related context before.

\sbr

Let us give a brief description of our idea. For each $\la>0$, we
let $\mu_\la=\mu_{\la,\Om_\rho}$ denote the positive Borel measure
on $\fg^*$ defined by $\mu_\la(A)=\mu_{\Om_\rho}(\Om_\rho\cap
\la\cdot A)$ for every Borel subset $A\subset\fg^*$. Note that
$\mu_\la$ is supported on $\la^{-1}\cdot\Om_\rho$, which is
another coadjoint orbit in $\fg^*$. Now $\Om_\rho$ is closed in
$\fg^*$, and there is a coordinate system on $\Om_\rho$ which
identifies $\Om_\rho$ with $\bR^{2n}$, such that $\mu_{\Om_\rho}$
corresponds to the usual Lebesgue measure under this
identification (both of these statements hold for arbitrary
nilpotent Lie algebras). In particular, we see that each $\mu_\la$
can be identified with a positive linear functional on the space
$C_c(\fg^*)$ of compactly supported continuous functions on
$\fg^*$. Note also that, if $A$ is a neighborhood of $0$ in
$\fg^*$, then, as $\la\to +\infty$, the sets $\Om_\rho\cap\la\cdot
A$ exhaust all of $\Om_\rho$; thus, $\mu_\la(A)\to +\infty$. Let
us now suppose that there exists a function $f(\la)$ such that the
functionals $f(\la)\cdot\mu_\la\in C_c(\fg^*)^*$ have a nonzero
weak-* limit $f_0\in C_c(\fg^*)^*$. By the Riesz representation
theorem, $f_0$ corresponds to a positive Borel measure $\mu_0$ on
$\fg^*$. We define $\tilde{Q}=\Supp(\mu_0)$, and
$\tilde{\nu}=\mu_0\big\vert_{\tilde{Q}}$. Then $\tilde{Q}$ is a
conical $G$-invariant subset of $\fg^*$, and the $G$-invariance of
$\tilde{\nu}$ is automatic, since each of the measures $\mu_\la$
is $G$-invariant.

\subsection{}\label{ss:quasihomog} For simplicity, we will refer to
the construction described in the previous paragraph as the {\em
scaling construction}. Due to its ``homogeneous'' nature, it is
not surprising that in applying the construction to the
computation of spectral asymptotics of \sch operators, one has to
require a certain homogeneity condition on the potentials. We will
say, somewhat imprecisely, that \eqref{e:schr} is a {\em \sch
operator with quasi-homogeneous potentials} if $V(x)$ and $B(x)$
are {\em quasi-homogeneous} polynomials {\em of the same weight};
this means that there exists an $n$-tuple of positive rational
numbers $\ga=(\ga_1,\dotsc,\ga_n)$ such that for all $t\in\bR$,
$t>0$, and all $x\in\bR^n$, we have
\[
V\pare{t^{\ga_1}x_1,\dotsc,t^{\ga_n}x_n}=t\cdot V(x)
\quad\text{and}\quad
B\pare{t^{\ga_1}x_1,\dotsc,t^{\ga_n}x_n}=t\cdot B(x).
\]
Similarly, if $V(x)$ and all the components $b_{jk}(x)$ of the
magnetic tensor are {\em homogeneous} polynomials {\em of the same
degree}, we will refer to $H=H(a)+V$ as a {\em \sch operator with
homogeneous potentials}.

\sbr

We will prove that in the quasi-homogeneous situation where the
classical formulas of Weyl and Colin de Verdi\`ere are applicable,
our construction gives the same result as the ``pushforward''
construction described in \S\ref{ss:pushforward}. On the other
hand, in the intermediate degeneration examples that we have
studied, it also produces the ``correct'' measure space $(Q,\nu)$,
even though the pushforward construction no longer applies.

\subsection{} We remark that our scaling construction makes sense
for any nilpotent Lie algebra. Indeed, let $\fg$ be a finite
dimensional nilpotent Lie algebra over $\bR$ and $\Om\subset\fg^*$
a coadjoint orbit. It is known (cf., e.g., \cite{reps}, Ch.~I)
that $\Om$ is a closed (in fact, Zariski closed) submanifold of
$\fg^*$. Moreover, it follows from the explicit parameterization
obtained in \cite{Bon} that there exists a {\em polynomial} map
$\varphi:\bR^{2n}\to\fg^*$ which is a diffeomorphism onto $\Om$,
and such that under this diffeomorphism $\mu_\Om$ corresponds to
the standard Lebesgue measure on $\bR^{2n}$.

\sbr

As before, for every $\la>0$, we define a positive Borel measure
$\mu_\la$ on $\fg^*$ by \[\mu_\la(A)=\mu_\Om(\Om\cap\la\cdot
A)=\meas\pare{\varphi^{-1}(\la\cdot A)},\] where $\meas$ is the
Lebesgue measure. Since $\varphi$ is proper, we see that
$C_c(\fg^*)\subset L^1(d\mu_\la)$ for each $\la>0$. In particular,
once again, we can identify $\mu_\la$ with a positive linear
functional on $C_c(\fg^*)$, and the rest of our construction goes
through without any changes. It will be apparent from the
computations of explicit examples in Sections \ref{s:ex1} and
\ref{s:examples} below that the scaling construction is closely
related to the geometry of the embedding
$\Om\hookrightarrow\fg^*$.

\sbr

Furthermore, let us choose an arbitrary Euclidean structure on
$\fg^*$ and a corresponding orthonormal basis, and let
$\pare{P_1(x),\dotsc,P_t(x)}$ be the coordinates of $\varphi(x)$
with respect to this basis. Fix $A\subset\fg^*$ to be the unit
ball around the origin. We have seen in \S\ref{ss:scaling} that
the growth of $\mu_\la$ as $\la\to +\infty$ is related to the
behavior of the function
\[
G(\la)=\mu_\la(A)=\meas\cur{x\in\bR^{2n}\st \sum_{j=1}^t
P_j(x)^2\leq\la^2}.
\]
Now $P(x):=\sum P_j(x)^2$ is a polynomial function on $\bR^{2n}$
with $P(x)\to +\infty$ as $\norm{x}\to\infty$, and it follows from
the results of Nilsson on the growth of ``algebraic integrals''
\cite{Ni1,Ni2} that there exist positive reals $c,C,\al$ and a
nonnegative integer $\be$ such that
\[
C^{-1}\cdot\la^\al\cdot(\log\la)^\be\leq G(\la)\leq
C\cdot\la^\al\cdot(\log\la)^\be \quad\text{for all }\la>c.
\]
This result is important for the formulation of our conjectures.
(The work of Nilsson was used, in a similar situation, by Manchon
\cite{Ma1}.)

\subsection{}
The idea of applying representation-theoretic methods to the study
of partial differential operators is not new (see, for example,
\cite{HN} and the references therein). Several authors have
studied extensions of the known results about \sch operators to
the differential operators arising from unitary representations of
general nilpotent Lie groups. In \cite{LMN}, an analogue of
\eq{e:est} for the image under an irreducible representation of
the ``sublaplacian'' on a stratified nilpotent Lie algebra was
obtained. Manchon in \cite{Ma1} has generalized the approximate
spectral projection method of Tulovski\v{i} and Shubin \cite{TS}
to prove a Weyl-type asymptotic formula for elliptic operators
associated to representations of arbitrary nilpotent Lie groups.
In \cite{Ma2, Ma3}, this result was generalized to arbitrary Lie
groups, more precisely, to the representations corresponding to
closed tempered coadjoint orbits for which Kirillov's character
formula is valid. Notice, however, that \cite{Ma1, Ma2, Ma3} uses
the initial form of the approximate spectral projection method,
which requires high regularity of the symbol. In particular, if a
degeneration of any kind is present, this form of the approximate
spectral projection method does not work at all. For a general
version of the approximate spectral projection method, and
applications to various classes of degenerate and hypoelliptic
operators, see \cite{L1, L2, L3}.

\subsection{} The plan of the paper is as follows. In Section
\ref{s:conj}, we recall the main necessary definitions, and
formulate our conjectures. We also state several useful theorems
and propositions which are needed to explain our construction and
apply it to computing the leading term of the asymptotics of
Schr\"odinger operators (cf. Sections \ref{s:ex1} and
\ref{s:examples}). The proofs of these results are given in
Appendix A. We do not claim any originality here: some of these
theorems and propositions can be deduced from either the general
results on nilpotent Lie groups \cite{Bon,reps} or the previous
work on applications of representation theory to differential
operators \cite{HN,LMN}, while others are straightforward
extensions of known facts. However, there are several reasons for
including a complete account of these results.

\sbr

First, we explicitly isolate the class of Lie algebras that are
``responsible'' for \sch operators with polynomial magnetic and
electric fields: these are precisely the nilpotent Lie algebras
$\fg$ such that the commutator $[\fg,\fg]$ is abelian. Indeed, if
$H$ is a \sch operator \eqref{e:schr} with polynomial potentials
$a,V$, and $\fg$ is the associated Lie algebra as defined in
\S\ref{ss:pushforward}, then it is clear that $\fg$ enjoys this
property. A ``converse'' to this statement is given in
\S\ref{ss:reps}. It seems that this observation has not explicitly
appeared in the literature before.

\sbr

Moreover, most of the works relating differential operators to
representation theory of nilpotent Lie groups deal only with {\em
stratified} Lie algebras (cf. \cite{HN,LMN}), i.e., Lie algebras
$\fg$ admitting a decomposition
$\fg=\fg_1\oplus\fg_2\oplus\dotsb\oplus\fg_s$ as a direct sum of
vector subspaces, such that $[\fg_j,\fg_k]\subseteq\fg_{j+k}$
($\fg_j=(0)$ for $j>s$), and $\fg$ is generated by $\fg_1$ as a
Lie algebra. However, there are situations where the Lie algebra
arising from a \sch operator with polynomial potentials admits no
natural grading, as the simple example studied in Section
\ref{s:ex1} already shows. On the other hand, the theory we
develop in Section \ref{s:conj} makes no use of a grading on
$\fg$.

\sbr

Lastly, we state the results of Section \ref{s:conj} in exactly
the form one needs to be able to use our conjectural formula for
practical computations. The fact that we are only dealing with a
special class of nilpotent Lie algebras allows us to give rather
simple proofs and make our text essentially self-contained. In
some sense, this is easier than trying to deduce all the results
in the form we need from the works of other authors.

\sbr

In Section \ref{s:ex1}, we use an example of the Schr\"odinger
operator in 2D with zero electric potential and magnetic tensor
$b(x)=x_1^2-x_2$ (this is an example of strong degeneration), both
to illustrate in detail the use of our conjectural formula, and to
explain the direct variational method of the calculation of the
asymptotics of the spectrum, which can be applied to many other
classes of \sch operators with degenerate potentials. In Section
\ref{s:examples}, we study the weak degeneration case for
operators without either  magnetic or electrical potential, and
deduce from our conjecture the classical Weyl formula and Colin de
Verdi\`ere's formula, respectively. At the end of the section,
using the direct variational method, we prove a general theorem
for \sch operators with polynomial electric and magnetic fields,
which gives the leading term of the asymptotics under fairly weak
conditions on the function $\Pst$. These conditions are satisfied
in many cases of weak and intermediate degeneration but not in the
strong degeneration case. In particular, we prove that in the case
of a quasi-homogeneous electric potential, the classical Weyl
formula holds if and only the integral in this formula converges;
and our general conjectural formula also gives the classical Weyl
formula if and only if this condition is satisfied. In
\sect{s:example2}, we consider the \sch operator in 2D with
magnetic tensor $b(x)=x_1^k x_2^l$ and zero electric potential. In
the case $k\neq l$  we have the strong degeneration, and in the
case $k=l$ - the intermediate one. Finally, we briefly mention the
example of the \sch operator in 3D, $H=-\De+x_1^{2k} x_2^{2l}
x_3^{2p}$, with $p<k\leq l$. The most technical parts of the
proofs are delegated to Appendix B.

\subsection{Notation and terminology}
Our interest lies mainly in the degenerate cases where the
classical formulas of Weyl or Colin de Verdi\`ere are not
applicable. However, as we show in \S\ref{ss:Weyl} and
\S\ref{ss:CdV}, our formula works just as well in the regular
situations. In terms of the distinction between ``strong
degeneration'' and ``weak/intermediate degeneration'', the cases
with no degeneration at all belong to the latter type.

\sbr

We use the following notation: $\bN=\{1,2,3,\dotsc\}$ (the set of
natural numbers), $\bZ_+=\{n\in\bZ\big\vert n\geq 0\}$,
$\bR_+=\{x\in\bR\big\vert x\geq 0\}$,
$\bR^\times=\bR\setminus\{0\}$. If $S$ is any finite set, we write
$\# S$ for the number of elements of $S$. If $d\in\bZ_+$, we
denote by $\abs{v_d}$ the volume of the unit ball in $\bR^d$ (our
convention is that $\abs{v_0}=1$). Given $a\in\bR$, we define
$(a)_+$ by $(a)_+=a$ if $a>0$, $(a)_+=0$ if $a\leq 0$.

\section{Main results and conjectures}\label{s:conj}

\subsection{Schr\"odinger operators and unitary
representations}\label{ss:reps} As explained in the introduction,
our goal is to write down a conjectural formula for the leading
term of the spectral asymptotics of a \sch operator with
quasi-homogeneous polynomial potentials. However, in this
subsection the quasi-homogeneity condition plays no role. Thus, we
fix a \sch operator $H=H(a)+V$ with polynomial potentials $a,V$.
We define $\fg_H$, as in the introduction, to be the real Lie
algebra generated by the polynomial differential operators
$L_j=\dd/\dd x_j + \sqrt{-1}\cdot a_j(x)$ ($1\leq j\leq n$) and
$L_0=\sqrt{-1}\cdot V(x)$. It is clear that $\fg_H$ is a finite
dimensional nilpotent Lie algebra. Moreover, the commutator
$[\fg_H,\fg_H]$ consists only of multiplication operators, and
thus $[\fg_H,\fg_H]$ is an abelian ideal of $\fg_H$.

\sbr

Recall that two \sch operators, $H=H(a)+V$ and $H'=H(a')+V'$, are
said to be {\em gauge equivalent} if $V=V'$ and the corresponding
magnetic tensors are the same: $B=B'$. By Poincar\'e's lemma, the
last condition is equivalent to the existence of a differentiable
function $\phi:\bR^n\to\bR$ such that
$a_j'(x)=a_j(x)+\dd_j\phi(x)$ for all $1\leq j\leq n$. If such a
$\phi$ exists, it is easy to check that the unitary operator
$\exp(i\phi)$ conjugates $H$ into $H'$; in particular, $H$ and
$H'$ have the same spectrum. On the other hand, if $H$, $H'$ are
\sch operators with polynomial potentials that are gauge
equivalent, then the corresponding Lie algebras $\fg_H$,
$\fg_{H'}$ are isomorphic, because the commutation relations in
$\fg_H$ depend only on $V(x)$, $b_{jk}(x)$ and their derivatives.

\sbr

By the ``tautological representation'' of $\fg_H$ we will mean the
representation of $\fg_H$ on $L^2(\bR^n)$ by (unbounded)
skew-adjoint operators that takes every element of $\fg_H$ to the
polynomial differential operator it represents. We note that,
unlike the case of finite dimensional representations, the problem
of lifting the tautological representation to a unitary
representation of the connected and simply connected nilpotent Lie
group $G=\exp\fg_H$ is not trivial. We will address this issue in
the theorem below.

\sbr

From a more abstract point of view, let $\fg$ be an arbitrary
finite dimensional nilpotent Lie algebra over $\bR$ such that
$[\fg,\fg]$ is abelian. A {\em sublaplacian} for $\fg$ is an
element $S\in\cU(\fg)_\bC$ which has the form
$S=-(L_1^2+\dotsb+L_N^2)-\sqrt{-1}\cdot L_0$, where
$L_0,L_1,\dotsc,L_N\in\fg$ are linearly independent elements that
generate $\fg$ as a Lie algebra, and $L_0$ commutes with
$[\fg,\fg]$. Note that we have extended the standard definition of
a sublaplacian (which does not contain the $L_0$ term) to include
the case of a \sch operator with nonzero electric potential. Then
we have the following result.

\begin{thm}\label{t:reps}
\begin{enumerate}[(a)]
\item Every unitary irreducible
representation of $G=\exp\fg$ has a natural realization in a space
$L^2(\bR^n)$, $n\leq N$, such that each element of $\fg$ maps to a
polynomial differential operator of order $\leq 1$; $L_0$ and all
elements of $[\fg,\fg]$ map to multiplication operators; and $S$
maps to a \sch operator with polynomial potentials which has
discrete spectrum if the image of $-\sqrt{-1}\cdot L_0$ is a
polynomial that is bounded below.
\item Conversely, if $H$ is a \sch operator \eqref{e:schr} with
polynomial potentials, there exists a \sch operator $H_0$ with
polynomial potentials which is gauge equivalent to $H$, such that
if $\fg=\fg_{H_0}$ and $S\in\cU(\fg)_\bC$ is the element
corresponding to $H_0$, then the tautological representation of
$\fg$ on $L^2(\bR^n)$ can be lifted to a unitary irreducible
representation of $G=\exp\fg$, which is irreducible if and only if
$H$ has discrete spectrum.
\end{enumerate}
\end{thm}

\sbr

It is important to have a concrete realization of each of the
representations of $\fg$ that arises from a unitary irreducible
representation of $G$. These will be discussed in detail in
\S\ref{ss:realizations}.
\begin{rem}
We would like to explain the reason for our change of notation: we
use $S$ to denote a sublaplacian in $\fg$, whereas we use
$H^\circ$ everywhere else to denote the sublaplacian that
naturally corresponds to a \sch operator $H$. The reason is that
in our proof of Theorem \ref{t:reps}, as well as some other proofs
presented in Appendix A, we often use the letter $H$ to denote a
closed subgroup of the group $G=\exp\fg$. This notation is used
throughout the literature on representation theory, so it does not
make sense to change it. On the other hand, we do not wish to
confuse a subgroup $H$ with a \sch operator $H$.
\end{rem}

\subsection{Preliminary version of the conjecture}\label{ss:vague}
Let us now formulate a preliminary version of our conjecture. Let
$H$ be a \sch operator \eqref{e:schr} with discrete spectrum and
quasi-homogeneous polynomial potentials, and let $\fg=\fg_H$ be
the associated Lie algebra. Since we are interested in $\sg(H)$,
we may assume, by Theorem \ref{t:reps}, that the tautological
representation of $\fg$ lifts to a unitary representation of $G$
on $L^2(\bR^n)$; moreover, this representation is then
irreducible, whence corresponds to a coadjoint orbit
$\Om\subset\fg^*$. Let $\mu_\Om$ be the Kostant measure on $\Om$;
for the precise normalization, see Definition \ref{d:canonical}.
Then we have the ``dilates'' $\mu_\la$ of the measure $\mu_\Om$,
as defined in the introduction:
$\mu_\la(A)=\mu_\Om(\Om\cap\la\cdot A)$, for every Borel subset
$A\subseteq\fg^*$. Furthermore, $H$ naturally defines an element
$H^\circ\in\cU(\fg)_\bC$, and the definition of $\fg$ implies that
$H^\circ$ is a sublaplacian for $\fg$. For any coadjoint orbit
$\Theta\subset\fg^*$, we denote by $H_\Theta$ the image of
$H^\circ$ in the \unirrep of $G$ that corresponds to $\Theta$ via
Kirillov's theory. By Theorem \ref{t:reps}, each $H_\Theta$ can be
naturally realized as a \sch operator with polynomial potentials.

\sbr

\noindent
{\sc Conjecture 1.} There exist a positive real number $\al$ and a
nonnegative integer $\be$ such that the weak limit
$\mu_0=\lim_{\la\to +\infty} \la^{-\al}\cdot(\log\la)^{-\be}\cdot
\mu_\la$ exists and is nonzero. Then $\mu_0$ is automatically
$G$-invariant; let $Q=\pare{\Supp\mu_0}/G$, and let
$\rho:\Supp\mu_0\to Q$ be the natural projection. Let $\nu$ be the
measure on $Q$ such that for every function $F\in C_c(\fg^*)$, we
have
\[
\int_{\fg^*} F\,d\mu_0 = \int_Q d\nu(q)\cdot \int_{\rho^{-1}(q)}
F(x)\, d\mu_q(x),
\]
where $d\mu_q$ denotes the Kostant measure corresponding to the
orbit $\rho^{-1}(q)$ (the existence of $\nu$ is proved in
Proposition \ref{p:quotient}). Then there exists a constant
$\kappa\geq 1$ such that
\begin{equation}\label{e:conj1}
N(\la,H)\sim \kappa\cdot (\log\la)^\be\cdot \int_Q N(\la,H_\Theta)
\,d\nu(\Theta) \quad \text{as}\ \la\to +\infty.
\end{equation}

\subsection{Precise version of the conjecture}\label{ss:precise}
We now formulate a more precise form of our conjecture---one that
essentially provides a formula for the constant $\kappa$ that
appears in \eqref{e:conj1}. To that end, we introduce the function
\begin{equation}\label{e:phistar}
\Phi^*(x)=\sum_{\al} \abs{\dd^\al V(x)}^{1/2} + \sum_{\al,j,k}
\abs{\dd^\al b_{jk}(x)}^{1/2};
\end{equation}
this is to be compared with the function $\Pst$ used in
\cite{HM,MN} (see \eq{e:psistar}). If, for example, $V\equiv 0$
and $B(x)$ grows regularly at infinity, then the terms
corresponding to $\al=0$ dominate both $\Psi^*$ and $\Phi^*$, so
we see that these two functions have the same asymptotic behavior
as $\norm{x}\to\infty$. However, in general, it may happen that
the function $\Psi^*(x)$ grows slower than the function
$\Phi^*(x)$.

\sbr

We keep the same notation and assumptions as in Conjecture 1. In
particular, since $H$ has discrete spectrum, both $\Phi^*$ and
$\Psi^*$ tend to $+\infty$ as $\norm{x}\to\infty$, so it makes
sense to define the functions
\[
G_1(\la)=\meas\cur{ x\in\bR^n \st \Phi^*(x)\leq\la}
\]
and
\[
G_2(\la)=\meas\cur{ x\in\bR^n \st \Psi^*(x)\leq\la},
\]
where $\meas$ stands for the usual Lebesgue measure.

\sbr

\noindent
{\sc Conjecture 2.} Assume that $H$ is a \sch operator on
$L^2(\bR^n)$ with discrete spectrum and quasi-homogeneous
potentials. Let $(Q,\nu)$ be defined as in Conjecture 1. Then one
of the following situations occurs.
\begin{enumerate}[(a)]
\item We have $G_2(\la)/G_1(\la)\to\infty$ as $\la\to +\infty$.
This is the {\em strong degeneration case}. Then Conjecture 1 is
valid with the normalization constant $\kappa=1$.
\item We have $G_2(\la) = O(G_1(\la))$ as $\la\to +\infty$. This
is the {\em weak/intermediate degeneration case.} Then there
exists a limit $\lim\limits_{\la\to +\infty} G_2(\la)/G_1(\la)$,
and Conjecture 1 is valid with $\kappa$ equal to the value of this
limit.
\end{enumerate}

\subsection{Concrete realization of
representations}\label{ss:realizations} Until the end of the
section, the quasi-homogeneity assumption will play no role. Let
$\fg$ be a real finite dimensional nilpotent Lie algebra such that
$[\fg,\fg]$ is abelian, and let
$S=-(L_1^2+\dotsb+L_N^2)-\sqrt{-1}\cdot L_0\in\cU(\fg)_\bC$ be a
sublaplacian. We wish to obtain concrete realizations of the
representations of $\fg$ induced by \unirreps of $G=\exp\fg$. Let
$\fh\subseteq\fg$ be a Lie subalgebra, and $H=\exp\fh$ the
corresponding connected and simply connected subgroup of $G$.
(\sch operators do not appear until the end of the section, so the
notation should not cause any confusion.) Fix $f\in\fg^*$. We say
that $\fh$ is {\em subordinate to} $f$ if
$f\big\vert_{[\fh,\fh]}\equiv 0$. Under this condition, $f$
defines a unitary character $\chi_f$ of $H$ via $\chi_f(\exp
h)=\exp(i\cdot f(h))$. Thus we may form the induced representation
$\rho_{f,\fh}=\Ind_H^G(\chi_f)$. This construction is reviewed in
\S\ref{ss:induced}. Kirillov's classification \cite{Ki} of
\unirreps of $G$ can be summarized as follows.

\sbr

Let us say that $\fh$ is a {\em polarization} of $\fg$ at $f$ if
$\fh$ is of maximal dimension among the subalgebras of $\fg$ that
are subordinate to $f$. Then $\rho_{f,\fh}$ is irreducible if and
only if $\fh$ is a polarization at $f$. Moreover, in this case,
$\rho_{f,\fh}$ does not depend on the choice of $\fh$, up to
unitary equivalence. Also, at every $f\in\fg^*$ there exists at
least one polarization. Thus, we write $\rho_f=\rho_{f,\fh}$ for
any choice of a polarization $\fh$ at $f$. Finally, every \unirrep
of $G$ is unitarily equivalent to $\rho_f$ for some $f\in\fg^*$,
and $\rho_{f_1}$, $\rho_{f_2}$ are unitarily equivalent if and
only if $f_1$, $f_2$ lie in the same coadjoint orbit of $G$.

\sbr

Let us define the alternating bilinear form
\begin{equation}\label{e:form}
B_f:\fg\times\fg\to\bR,\quad B_f(X,Y)=\pair{f}{[X,Y]}.
\end{equation}
Thus, a subalgebra $\fh\subseteq\fg$ is subordinate to $f$ if and
only if $\fh$ is isotropic with respect to $B_f$. One can prove
that $\fh$ is a polarization at $f$ if and only if $\fh$ is
maximally isotropic with respect to $B_f$ {\em as a linear
subspace}. In particular, all polarizations at $f$ have the same
dimension,
\[
\dim\fh=\frac{1}{2}\cdot\pare{\dim\fg+\dim\fg(f)},
\]
where
\[
\fg(f)=\Ker B_f = \cur{ X\in\fg\st f\pare{[X,Y]}=0\ \forall\,
Y\in\fg}.
\]
In our situation, we can give an elementary proof of the existence
of polarizations of a special form:
\begin{lem}\label{l:polarization}
Let $\fg$, $S$ be as above, and $f\in\fg^*$. Then there exists a
polarization $\fh$ of $\fg$ at $f$ such that
$\fh\supseteq\fg(f)+\bR L_0+[\fg,\fg]$, and hence $\fh$ is an
ideal of $\fg$. Moreover, $[\fh,\fh]\subseteq\fg(f)$, so $\fg(f)$
is an ideal of $\fh$.
\end{lem}

\sbr

Let us now fix a subalgebra $\fh\subset\fg$ subordinate to $f$,
but not necessarily a polarization at $f$, which satisfies the
requirement of the lemma: $\fg(f)+\bR L_0+[\fg,\fg]\subseteq\fh$.
Since $L_0,L_1,\dotsc,L_N$ generate $\fg$ as a Lie algebra, we
have $\fg=[\fg,\fg]+\Span_\bR\{L_0,L_1,\dotsc,L_N\}$, and hence, a
fortiori, $\fg=\fh+\Span_\bR\{L_1,\dotsc,L_N\}$. After reindexing,
we may assume that for some $0\leq n\leq N$, the elements
$L_1,\dotsc,L_n$ form a complementary basis to $\fh$ in $\fg$. (We
 allow $n=0$, which means that $\fh=\fg$.) For every element
$h\in\fh$, let us define a real polynomial $p_h(x)$ in $n$
variables $x=(x_1,\dotsc,x_n)$ by
\begin{equation}\label{e:action-of-h}
p_h(x)=\sum_{\al_1,\dotsc,\al_n\geq 0}
\frac{1}{\al_1!\dotsm\al_n!}\cdot f\pare{(\ad
L_1)^{\al_1}\dotsm(\ad L_n)^{\al_n}(h)}\cdot x_1^{\al_1}\dotsm
x_n^{\al_n}.
\end{equation}
\begin{prop}\label{p:realization}
There exists a realization of the representation
$\rho_{f,\fh}=\Ind_H^G(\chi_f)$ of the Lie group $G$ in the space
$L^2(\bR^n,dm)$ (where $dm$ is the Lebesgue measure) such that the
induced representation of $\fg$ takes every $h\in\fh$ to the
operator of multiplication by $\sqrt{-1}\cdot p_h(x)$, and takes
$L_j$, for $1\leq j\leq n$, to the operator $\dd/\dd
x_j+\sqrt{-1}\cdot a_j(x)$, where $a_j(x)\in\bR[x_1,\dotsc,x_n]$
is a certain polynomial.
\end{prop}
The practical applications of this proposition are based on the
obvious analogy between \eqref{e:action-of-h} and the usual
Taylor's formula.

\subsection{Coadjoint orbits and Kostant
measures}\label{ss:orbits} Let $G$ be any connected Lie group, and
$\fg$ its Lie algebra. If $f\in\fg^*$, we denote by $G(f)$ the
stabilizer of $f$ in $G$ (with respect to the coadjoint action),
and by $\fg(f)$ the Lie algebra of $G(f)$. If $\Om\subset\fg^*$ is
a coadjoint orbit, then for any point $f\in\Om$, the orbit map
$G\to\Om$, $g\mapsto (\Ad^* g)(f)$, identifies $\Om$ with the
homogeneous space $G/G(f)$, and hence identifies the tangent space
$T_f\Om$ with the quotient $\fg/\fg(f)$. The notation is
consistent with the one used in \S\ref{ss:realizations}: if $B_f$
is the alternating bilinear form on $\fg$ given by
$B_f(X,Y)=\pair{f}{[X,Y]}$, then it is easy to see that $\fg(f)$
is precisely the kernel of $B_f$. Moreover, $B_f$ induces an
alternating bilinear {\em nondegenerate} form $\om_f$ on
$\fg/\fg(f)\cong T_f\Om$. One then proves the following facts
(see, e.g., \cite{reps}, Ch.~II):
\begin{enumerate}[1)]
\item the forms $\om_f$ vary smoothly with $f$, thus defining a
nondegenerate differential $2$-form $\om_\Om$ on $\Om$;
\item the form $\om_\Om$ is closed, and thus a symplectic form on
$\Om$;
\item the form $\om_\Om$ is $G$-invariant.
\end{enumerate}
\begin{defin}\label{d:canonical}
The form $\om_\Om$ is called the {\em canonical symplectic form}
on the orbit $\Om$. The {\em Kostant measure} (or the {\em
canonical measure}) on the orbit $\Om$ is the positive Borel
measure $\mu_\Om$ associated with the volume form
\[
(2\pi)^{-n}\cdot \frac{1}{n!}\cdot \om_\Om^n \qquad \bigl(
n=\frac{1}{2}\dim\Om \bigr).
\]
(Note that $\dim\Om$ is even because $\Om$ admits a symplectic
form.)
\end{defin}

\sbr

It is clear that the Kostant measure is $G$-invariant. In the
remainder of this subsection we obtain an explicit
parameterization of the coadjoint orbits for the Lie algebras of
the type considered in \S\ref{ss:reps}, and we derive formulas for
the corresponding canonical symplectic forms and Kostant measures.
We note that explicit parameterizations of the dual space of a
(not necessarily nilpotent) Lie algebra have been studied by
various authors: see, e.g., \cite{long}. More recently, a very
fine stratification of $\fg^*$ for nilpotent $\fg$ has been
obtained in \cite{Bon}. We will use a result from loc. cit. in the
next subsection.

\sbr

In our subsequent computations (especially the ones that appear in
the concrete examples of Sections \ref{s:ex1}, \ref{s:examples}
and \ref{s:example2}) we will implicitly use the following result.
The proof is completely straightforward and is therefore omitted.
Let $\fg$ be a Lie algebra and $\fa\subset\fg$ an ideal. Write
$\fa^\perp$ for the annihilator of $\fa$ in $\fg^*$. The quotient
map $\fg\to\fg/\fa$ induces an isomorphism of vector spaces
$(\fg/\fa)^*\rar{\cong}\fa^\perp\hookrightarrow\fg^*$. Let $G$ be
a connected Lie group with Lie algebra $\fg$, and $A\subset G$ the
closed connected normal subgroup corresponding to $\fa$. The
adjoint action of $G$ on $\fg$ leaves $\fa$ stable, whence $G$
also acts on $\fg/\fa$ and on $(\fg/\fa)^*$. Then we have the
following
\begin{prop}\label{p:trivial}
\begin{enumerate}[(a)]
\item The isomorphism $(\fg/\fa)^*\to\fa^\perp$ above is
$G$-equivariant, and the action of $G$ on $(\fg/\fa)^*$ factors
through the quotient group $G/A$; thus the $G$-orbits in
$(\fg/\fa)^*$ are the same as the coadjoint orbits of $G/A$ in
$(\fg/\fa)^*$.
\item If $\Om\subset\fg^*$ is any coadjoint orbit, then either
$\Om\cap\fa^\perp=\emptyset$, or $\Om\subset\fa^\perp$. In the
latter case, $\Om$ is the image of a coadjoint orbit in
$(\fg/\fa)^*$. Moreover, the canonical symplectic form and the
Kostant measure on $\Om$ are the same whether we regard $\Om$ as a
coadjoint orbit for $G$ or as a coadjoint orbit for $G/A$.
\item If $G$ is simply connected and nilpotent, then the bijection
between the coadjoint orbits in $\fg^*$ that meet $\fa^\perp$ and
the coadjoint orbits in $(\fg/\fa)^*$, defined above, corresponds,
via Kirillov's theory, to the natural bijection between the
unitary irreducible representations of $G$ that are trivial on
$A$, and all unitary irreducible representations of $G/A$.
\end{enumerate}
\end{prop}

We return to the situation considered in \S\ref{ss:realizations}.
Thus $G$ is a connected and simply connected nilpotent Lie group
with Lie algebra $\fg$ such that $[\fg,\fg]$ is abelian. Fix a
point $f_0\in\fg^*$. We wish to parameterize the $G$-orbit $G\cdot
f_0\subset\fg^*$. As before, we assume we are given a sublaplacian
$S=-(L_1^2+\dotsb+L_N^2)-\sqrt{-1}\cdot L_0$ for $\fg$, and we let
$\fh$ be a real polarization of $\fg$ at $f_0$ provided by Lemma
\ref{l:polarization}: $\fg(f_0)+\bR L_0+[\fg,\fg]\subset\fh$.
Furthermore, we suppose that for some $1\leq n\leq N$,
$L_1,\dotsc,L_n$ is a complementary basis for $\fh$ in $\fg$.

\sbr

From now on, we also assume that $\fh$ is an {\em abelian} ideal
of $\fg$. To justify this assumption, we note that since $\fh$ is
an ideal of $\fg$, so is $\fa:=[\fh,\fh]$; on the other hand, by
the definition of a polarization, $f_0$ annihilates $\fa$. Thus,
$f_0$ induces a linear functional $\bar{f}_0$ on $\fg/\fa$. By
Proposition \ref{p:trivial}, the canonical inclusion
$(\fg/\fa)^*\hookrightarrow\fg^*$ gives an isomorphism of the
coadjoint orbit of $\bar{f}_0$ in $(\fg/\fa)^*$ onto the coadjoint
orbit of $f_0$ in $\fg^*$; moreover, this isomorphism preserves
the canonical symplectic form and the Kostant measure. Lastly,
note that since $\fa\subset\fg(f)$ by Lemma \ref{l:polarization},
it is clear that $\fh/\fa$ is a maximal isotropic subspace of
$\fg/\fa$ with respect to the form $B_{\bar{f}_0}$. Thus, from the
point of view of either the coadjoint orbit of $f_0$, or of the
corresponding unitary irreducible representation, nothing is lost
by passing from $\fg$ to $\fg/\fa$.
\begin{prop}\label{p:orbit}
With the notation above, assume that $\fh$ is abelian. The map
\[
\varphi : \bR^n\times\bR^n \rar{} \fg^*
\]
defined by
\[
\pair{\varphi(\xi,x)}{L_j}=\xi_j\qquad \text{for } 1\leq j\leq n,
\]
\[
\pair{\varphi(\xi,x)}{Y}= \sum_{\al_1,\dotsc,\al_n\geq 0}
\frac{1}{\al_1!\dotsm\al_n!}\cdot f_0\pare{(\ad
L_1)^{\al_1}\dotsm(\ad L_n)^{\al_n}(Y)}\cdot x_1^{\al_1}\dotsm
x_n^{\al_n}\quad \text{for all } Y\in\fh,
\]
is a diffeomorphism of $\bR^{2n}$ onto the coadjoint orbit of
$f_0$ in $\fg^*$.
\end{prop}
By a slight abuse of notation, we identify $\Om$ with $\bR^{2n}$
using the diffeomorphism $\varphi$, and in particular, we view
$(\xi,x)$ as coordinates on the orbit $\Om$. Let us define
polynomials $b_{jk}(x)$ by
\[
b_{jk}(x) = \pair{\varphi(0,x)}{[L_j,L_k]};
\]
note that if $\fg$ arises from a \sch operator with polynomial
potentials, and if $f_0$ restricts to the linear functional
$\sqrt{-1}\cdot P(x)\mapsto P(0)$ on the subspace of $\fg$
consisting of multiplication operators, then the $b_{jk}(x)$ are
precisely the components of the magnetic tensor of the operator.
We will prove the following
\begin{prop}\label{p:measure}
The canonical symplectic form and the Kostant measure on the orbit
$\Om$ are given by
\begin{equation}\label{e:symplform}
\om_\Om = \sum_{j=1}^n d\xi_j\wedge dx_j + \sum_{1\leq j<k\leq n}
b_{kj}(x)\, dx_j\wedge dx_k
\end{equation}
and
\[
\mu_\Om = (2\pi)^{-n}\cdot d\xi_1\dotsm d\xi_n dx_1\dotsm dx_n.
\]
In other words, if we identify $\mu_\Om$ with its extension by
zero to $\fg^*$, then we can write
\begin{equation}\label{e:measure}
\mu_\Om = (2\pi)^{-n}\cdot \varphi_*(d\xi_1\dotsm d\xi_n
dx_1\dotsm dx_n),
\end{equation}
where $\varphi_*$ denotes the pushforward by the map
$\varphi:\bR^{2n}\to\fg^*$.
\end{prop}

\subsection{Polynomial measures}\label{ss:polymeas} In this
subsection we collect a few results that will be used in the
computations of Sections \ref{s:ex1}, \ref{s:examples} and
\ref{s:example2}, and a few others that help motivate our
conjectures in \S\ref{ss:vague} and \S\ref{ss:precise}.

\sbr

Let $\bR^N$ and $\bR^n$ be two Euclidean spaces. A {\em polynomial
measure} on $\bR^N$ is a Borel measure of the form
$\mu=c_0\cdot\phi_*(dm)$, where $c_0$ is a positive constant,
$\phi:\bR^n\to\bR^N$ is a polynomial map, and $dm$ is the Lebesgue
measure on $\bR^n$. Recall that the {\em pushforward measure}
$\phi_*(dm)$ is defined by the formula
\[
\bigl(\phi_*(dm)\bigr)(A):=m(\mu^{-1}(A)),
\]
for every Borel subset $A\subseteq\bR^N$. We say that $\mu$ is a
{\em regular polynomial measure} if it is finite on compact sets,
or, equivalently, if $\phi^{-1}(K)$ has finite Lebesgue measure
for every compact subset $K\subset\bR^N$ (a sufficient condition
for this is that $\phi$ is a proper map, but this condition is not
necessary). In practice, we will often think of regular polynomial
measures as positive linear functionals on the space $C_c(\bR^N)$
of all compactly supported continuous functions on $\bR^N$.
Explicitly, $\mu$ corresponds to the linear functional $\La$
defined by
\[
\La(F)=\int_{\bR^N} F\,d\mu = c_0\int_{\bR^n} F(\phi(x))\, dx,
\]
where the second integral is the usual Lebesgue integral.
Positivity here refers to the statement that if $F\in C_c(\bR^N)$
is nonnegative, then $\La(F)\geq 0$. It follows from Propositions
\ref{p:orbit} and \ref{p:measure} that if $\fg$ is a nilpotent Lie
algebra such that $[\fg,\fg]$ is abelian, then the Kostant measure
$\mu_\Om$ corresponding to each coadjoint orbit $\Om\subset\fg^*$
is a regular polynomial measure on $\fg^*$.

\sbr

Suppose that $\{\mu_j\}_{j=1}^\infty$ is a sequence of regular
polynomial measures on $\bR^N$. We will say that this sequence has
a {\em weak limit} if the limit
\[
\La(F):=\lim\limits_{j\to\infty} \int_{\bR^N} F\,d\mu_j
\]
exists for every $F\in C_c(\bR^N)$. Note that, if this is the
case, then the limit automatically defines a positive linear
functional on the space $C_c(\bR^N)$. It then follows from the
Riesz Representation Theorem that $\La$ itself corresponds to a
regular Borel measure $\mu$ on $\bR^N$. We then say that the
sequence $\{\mu_j\}$ {\em converges weakly} to $\mu$. The
following two results on weak limits will be used in Sections
\ref{s:ex1}, \ref{s:examples} and \ref{s:example2}.
\begin{prop}\label{p:measlim}
Suppose that $N=N_1+N_2$, that we are given a polynomial map
$\phi:\bR^n\to\bR^{N_1}$, and a sequence of polynomial maps
$\psi_j:\bR^n\to\bR^{N_2}$. Assume moreover that $\phi^{-1}(K)$
has finite Lebesgue measure for every compact subset
$K\subset\bR^{N_1}$, and that the sequence $\psi_j$ converges
pointwise to a polynomial map $\psi:\bR^n\to\bR^{N_2}$. Form the
polynomial maps $\tilde{\phi}_j:\bR^n\to\bR^N$ and
$\tilde{\phi}:\bR^n\to\bR^N$ by defining
\[
\tilde{\phi}_j(x)=\bigl(\phi(x),\psi_j(x)\bigr)\quad\text{and}\quad
\tilde{\phi}(x)=\bigl(\phi(x),\psi(x)\bigr)\quad\forall\,x\in\bR^n.
\]
Further, set $\mu_j=(\tilde{\phi}_j)_*(dm)$ and
$\mu=\tilde{\phi}_*(dm)$, where $dm$ is the Lebesgue measure on
$\bR^n$. Then $\mu_j$, $\mu$ are regular polynomial measures on
$\bR^N$, and $\mu_j\to\mu$ weakly.
\end{prop}

\mbr

As we will see, Proposition \ref{p:measlim} is useful for carrying
out the computation of the limiting measure $\mu_0$ that appears
in our conjectured formula in the cases of strong or weak
degeneration. In the intermediate degeneration case, the
computation turns out to be slightly more complicated, and the
following result (Proposition \ref{p:limborel}) is needed. We keep
the same notation as in \S\ref{ss:polymeas}. We say that a
positive Borel measure on $\bR^N$ is {\em c-finite} if it is
finite on compact sets; such a measure is then automatically
regular. We define weak limits for the c-finite measures in the
same way as for polynomial measures. Further, we introduce the
following notation. Given $a,b\in\bR^N$, we write
\[
[a,b)=\cur{ x\in\bR^N \st a_j\leq x_j<b_j \text{  for  } 1\leq
j\leq n}.
\]
Also, for $x\in\bR^N$, we set $\norm{x}_\infty=\max\limits_{1\leq
j\leq n}\abs{x_j}$. Finally, put $\vone=(1,1,\dotsc,1)\in\bZ^N$.

\begin{prop}\label{p:limborel}
Consider a family of $c$-finite positive Borel measures
$\{\mu_\la\}_{\la>0}$ on $\bR^N$. Assume that there exists a
countable set $E\subset\bR$ such that for each $a,b\in\bR^N$ with
$a_j\not\in E$, $b_j\not\in E$ for all $1\leq j\leq N$, there
exists a limit
\[
\nu\pare{[a,b)}=\lim\limits_{\la\to +\infty} \mu_\la\pare{[a,b)}.
\]
Then there exists a unique c-finite positive Borel measure $\mu_0$
on $\bR^N$ such that
\begin{equation}\label{e:star}
\int_{\bR^N} F(x)\, d\mu_0(x) = \lim\limits_{\la\to +\infty}
\int_{\bR^N} F(x)\, d\mu_\la(x) \quad \forall\, F\in C_c(\bR^N).
\end{equation}
Moreover, if $\nu_0$ is any c-finite positive Borel measure on
$\bR^N$ such that $\nu_0\pare{[a,b)}=\nu\pare{[a,b)}$ for all
$[a,b)$ as above, then \eqref{e:star} is satisfied with $\nu_0$ in
place of $\mu_0$; in particular, $\nu_0=\mu_0$.
\end{prop}

\mbr

Let us now recall a result of Nilsson \cite{Ni1} and explain its
relation to our conjectural formula. With the same notation as
above, let $\mu=\phi_*(dm)$ be a regular polynomial measure on
$\bR^N$, where $\phi:\bR^n\to\bR^N$ is a polynomial map given by
$\phi(x)=\pare{p_1(x),\dotsc,p_N(x)}$, and $dm$ is the Lebesgue
measure on $\bR^n$. Given $\la\in\bR$, $\la>0$, we define a
polynomial measure $\mu_\la$ on $\bR^N$ by
$\mu_\la(A)=\mu(\la^{-1}\cdot A)$, for each Borel subset
$A\subseteq\bR^N$. In particular, let $P(x)=\sum_{j=1}^N
p_j(x)^2$. If $A=B_1(0)$ is the unit ball around the origin in
$\bR^N$, then $\mu_\la(B_1(0))=\meas\{x\in\bR^n\big\vert
P(x)\leq\la^2\}$. We now state a result which is a special case of
Theorem 1 in \cite{Ni1}; the latter, in turn, is based on the
results of \cite{Ni2}.
\begin{thm}[Nilsson] Let $P(x)$ be a real polynomial on $\bR^n$
such that $P(x)\to +\infty$ as $\norm{x}\to\infty$, and set
\[
G(\la)=\meas\cur{x\in\bR^n \st P(x)\leq\la}.
\]
Then there exist positive reals $c,C,\al$ and a nonnegative
integer $\be$ such that
\[
C^{-1}\cdot\la^\al\cdot(\log\la)^\be \leq G(\la)\leq
C\cdot\la^\al\cdot(\log\la)^\be \quad\text{for all } \la>c.
\]
\end{thm}
It is clear that $\al$ and $\be$ are uniquely determined by the
polynomial $P(x)$. If $P=\sum_j p_j^2$ is as above, and $\al,\be$
are as in the theorem, then we see that
$\mu_\la(B_1(0))=O\pare{\la^{2\al}\cdot(\log\la)^\be}$ as $\la\to
+\infty$. But then, since the Lebesgue measure is a doubling
measure, we obtain that for every fixed $R>0$, we have
$\mu_\la(B_R(0))=O\pare{\la^{2\al}\cdot(\log\la)^\be}$ as $\la\to
+\infty$. So, the existence of a nonzero weak limit of
$\la^{-2\al}\cdot(\log\la)^{-\be}\cdot\mu_\la$ as $\la\to +\infty$
is plausible, which explains part of the statements of our
Conjectures 1 and 2.

\mbr

The final topic of the section is the decomposition of the
invariant measure on the dual space of a Lie algebra into an
integral of the Kostant measures. Let $\fg$ be an arbitrary finite
dimensional nilpotent Lie algebra over $\bR$, and $\fg^*$ its dual
space, with the coadjoint action of the Lie group $G=\exp\fg$. Let
$\fg^*/G$ denote the space of coadjoint orbits, and
$\rho:\fg^*\to\fg^*/G$ the quotient map. Write $\cM$ for the
$\sg$-algebra on $\fg^*/G$ defined by $A\in\cM$ $\iff$
$\rho^{-1}(A)\subseteq\fg^*$ is a Borel subset.
\begin{prop}\label{p:quotient}
Let $\mu$ be a positive Borel measure on $\fg^*$ which is
invariant under the coadjoint action of $G$. Then there exists a
unique positive measure $\nu$ on $\cM$ such that for every $F\in
C_c(\fg^*)$, we have
\[
\int_{\fg^*} F(x)\,d\mu(x) = \int_{\fg^*/G} d\nu(y)\cdot
\int_{\rho^{-1}(y)} F(x)\,d\mu_y(x),
\]
where $\mu_y$ is the Kostant measure corresponding to the
coadjoint orbit $\rho^{-1}(y)\subset\fg^*$.
\end{prop}

\section{An inhomogeneous example}\label{s:ex1}

\subsection{Application of the conjectural
formula}\label{ss:inhomog} In this section we consider in detail
the two-dimensional \sch operator
\[
H=-\frac{\dd^2}{\dd x_1^2} - \left(\frac{\dd}{\dd
x_2}+\sqrt{-1}\cdot
\pare{x_1^3/3-x_1 x_2}\right)^2
\]
with zero electric potential and magnetic tensor $b(x)=x_1^2-x_2$.
Note that $b(x)$ is quasi-homogeneous, but not homogeneous.
Furthermore, $b(x)$ vanishes along the parabola $x_2=x_1^2$ and,
in particular, does not grow at infinity in all directions (so the
formula of Colin de Verdi\`ere is not applicable). However, the
condition of the criterion of \cite{HM} is verified, so $H$ has
discrete spectrum. We wish to apply our conjecture to compute the
leading term of the asymptotics of $N(\la,H)$ as $\la\to +\infty$.

\sbr

First we must decide if we are in the strong degeneration case or
the intermediate degeneration case. To that end, we consider the
functions $\Phi^*$, $\Psi^*$ given by \eqref{e:phistar},
\eqref{e:psistar}. In our situation they become
\[
\Phi^*(x)=\abs{x_1^2-x_2}^{1/2}+\abs{2x_1}^{1/2}+c_1
\]
and
\[
\Psi^*(x)=\abs{x_1^2-x_2}^{1/2}+\abs{2x_1}^{1/3}+c_2,
\]
where $c_1$, $c_2$ are constants. We then have
\begin{eqnarray*}
G_1(\la)&=&\meas\cur{ x\in\bR^2 \st \Phi^*(x)\leq\la} \\
&=&\meas\cur{ x\in\bR^2 \st
\abs{x_1^2-x_2}^{1/2}+\abs{2x_1}^{1/2}\leq\la-c_1} \\
&=& \frac{1}{3}\cdot (\la-c_1)^4 \sim \la^4/3 \quad \text{as }
\la\to +\infty
\end{eqnarray*}
(the straightforward computation using Fubini's theorem is
omitted). Similarly,
\begin{eqnarray*}
G_2(\la)&=&\meas\cur{ x\in\bR^2 \st \Psi^*(x)\leq\la} \\
&=&\meas\cur{ x\in\bR^2 \st
\abs{x_1^2-x_2}^{1/2}+\abs{2x_1}^{1/3}\leq\la-c_2} \\
&=& \frac{1}{5}\cdot (\la-c_2)^5 \sim \la^5/5 \quad \text{as }
\la\to +\infty.
\end{eqnarray*}
Thus, we are in the strong degeneration case.
\begin{rem}
In general, it will be difficult or impossible to obtain precise
formulas for the functions $G_1(\la)$ and $G_2(\la)$, such as
above. However, it is usually not very hard to either prove an
estimate showing that $G_2(\la)/G_1(\la)\to\infty$ as $\la\to
+\infty$, or find $\lim_{\la\to +\infty} G_2(\la)/G_1(\la)$ in
case the limit exists and is finite.
\end{rem}

\mbr

Next we let $\fg=\fg_H$ be the Lie algebra associated to $H$, as
in Section \ref{s:conj}. Let us fix, once and for all, the
following basis of $\fg$:
\[
L_1=\frac{\dd}{\dd x_1},\ \ L_2=\frac{\dd}{\dd
x_2}+\sqrt{-1}\cdot\pare{x_1^3/3-x_1 x_2},\ \ X=\sqrt{-1}\cdot
(x_1^2-x_2),\ \ Y=\sqrt{-1}\cdot x_1,\ \ Z=\sqrt{-1}.
\]
We use this basis to identify $\fg^*$ with $\bR^5$; explicitly,
the identification is given by
\[
f\mapsto\pare{f(L_1),f(L_2),f(X),f(Y),f(Z)}.
\]
Write $\fh=\Span_\bR\{X,Y,Z\}$, and define a linear functional
$f_0:\fg\to\bR$ by $f_0(L_j)=0$ for $j=1,2$, and
$f_0(\sqrt{-1}\cdot P(x))=P(0)$ for all $\sqrt{-1}\cdot
P(x)\in\fh$. It follows from Proposition \ref{p:realization} that
the induced representation $\rho_{f_0,\fh}$ coincides with the
tautological representation of $\fg$ (at least up to changing the
action of $L_1$ and $L_2$, so that $H$ is replaced by a
gauge-equivalent operator). In particular, $\rho_{f_0,\fh}$ is
irreducible (by Theorem \ref{t:reps}), so $\fh$ is a polarization
of $\fg$ at $f_0$. Let $\Om=G\cdot f_0$ be the coadjoint orbit of
$f_0$ in $\fg^*$, and let $\mu_\Om$ be the corresponding Kostant
measure. From now on, to simplify notation, we will implicitly
identify $\mu_\Om$ with its extension by zero to all of $\fg^*$.

\sbr

By Proposition \ref{p:measure}, the orbit $\Om$ is parameterized
by the map $\varphi:\bR^4\to\bR^5\cong\fg^*$ given by
$\varphi(\xi,x)=(\xi_1,\xi_2,x_1^2-x_2,x_1,1)$, and, moreover, we
have
\[
\mu_\Om = (2\pi)^{-2}\cdot \varphi_*(d\xi dx).
\]
(As usual, we write $\xi=(\xi_1,\xi_2)$, $d\xi=d\xi_1 d\xi_2$,
etc.) To apply our scaling construction, we fix $F\in C_c(\fg^*)$
and study the integrals
\[
\int_{\bR^4}
F\pare{\la^{-1}\xi_1,\la^{-1}\xi_2,\la^{-1}(x_1^2-x_2),\la^{-1}x_1,\la^{-1}}\,
d\xi dx
\]
for $\la\gg 0$. Making the change of variables $\xi_j=\la\xi_j'$,
$\la^{-1}(x_1^2-x_2)=x_1'$, $\la^{-1}x_1=x_2'$, the integral above
becomes
\[
\la^4\cdot \int_{\bR^4} F\pare{\xi_1',\xi_2',x_1',x_2',\la^{-1}}\,
d\xi' dx'.
\]
Now $\la^{-1}\to 0$ as $\la\to +\infty$, and, using Proposition
\ref{p:measlim}, we find that there exists a weak limit
\[
\mu_0=\lim\limits_{\la\to +\infty} \la^{-4}\cdot\mu_\la =
(2\pi)^{-2}\cdot \psi_*(d\eta dy),
\]
where $\psi:\bR^4\to\bR^5\cong\fg^*$ is given by $(\eta,y)\mapsto
(\eta_1,\eta_2,y_1,y_2,0)$.

\sbr

In particular, $\mu_0$ is supported on the annihilator
$\fa^\perp\subset\fg^*$, where $\fa=\bR\cdot Z$ is an ideal of
$\fg$. Using Proposition \ref{p:trivial}, we pass from $\fg$ to
the quotient algebra $\fgb=\fg/\fa$. Naturally, we use the
following basis of $\fgb$: $\Lb_1$, $\Lb_2$, $\Xb$, $\Yb$, where
the bar denote the image of an element of $\fg$ under the quotient
map. The commutation relations in $\fgb$ are determined by
$[\Lb_1,\Lb_2]=\Xb$, $[\Lb_1,\Xb]=\Yb$, and all the other brackets
are zero. Furthermore, we see that $\mu_0$ is, up to the multiple
$(2\pi)^{-2}$, just the Lebesgue measure on $\fgb^*$. In
particular, the set of elements $f\in\fgb^*$ that vanish on $\Yb$
has measure zero with respect to $\mu_0$, and is also invariant
under the coadjoint action because $\Yb$ is central in $\fgb$.
Therefore this set can be ignored in the subsequent computations.

\sbr

Fix an element $f\in\fgb^*$ such that $f(\Yb)\neq 0$. It is clear
that $\Span_\bR\{\Lb_2,\Xb,\Yb\}$ is a polarization of $\fgb$ at
$f$. Using Proposition \ref{p:orbit}, we see that the coadjoint
orbit of $f$ in $\fgb^*$ is parameterized by the map
\[
\varphi_f:\bR^2\to\bR^4\cong\fgb^*,\quad (\zeta,z)\mapsto
\pare{\zeta,\, f(\Lb_2)+z f(\Xb)+ z^2 f(\Yb),\, f(\Xb)+2z
f(\Yb),\, f(\Yb)}.
\]
Moreover, by Proposition \ref{p:measure}, the corresponding
Kostant measure $\mu_f=\mu_{\Om_f}$ is given by
$\mu_f=(2\pi)^{-1}\cdot\pare{\varphi_f}_*(d\zeta dz)$. Now we must
find a convenient parameterization of the space $Q$ that appears
in our conjectures, which, in this case, is $\fgb^*/\Gb$.

\sbr

Given a ``generic'' coadjoint orbit $\Theta\subset\fgb^*$ (i.e.,
one that does not meet $\Yb^\perp$), it follows from the previous
paragraph that there is a unique $f\in\Theta$ with $f(\Xb)=0$.
Moreover, $\Theta$ is then determined uniquely by the value
$f(\Lb_2)$. In other words, if we define a map
\[
\fgb^*\setminus\Yb^\perp\to\bR^\times \times\bR
\]
by
\[
f\mapsto\left( f(\Yb),\,
f(\Lb_2)-\frac{3}{4}\cdot\frac{f(\Xb)^2}{f(\Yb)} \right),
\]
then the fibers of this map are precisely the generic coadjoint
orbits in $\fgb^*$. In coordinate form, the map is given by
$(z_1,z_2,z_3,z_4)\mapsto \pare{z_4,z_2-(3z_3^2)/(4z_4)}$. This
map clearly admits a smooth section given by $(a,b)\mapsto
(0,b,0,a)$. Thus, we have identified $Q$ with
$\bR^\times\times\bR$ (up to a set of measure zero, which we have
ignored in the beginning).

\sbr

Now the Kostant measure for the coadjoint orbit corresponding to
the point $(a,b)\in\bR^\times\times\bR$ is given by
\[
\mu_{(a,b)}=(2\pi)^{-1}\cdot \pare{\varphi_{(a,b)}}_*(d\zeta dz),
\]
where $\varphi_{(a,b)}:(\zeta,z)\mapsto (\zeta,az^2+b,2az,a)$. To
complete our construction, we must decompose the measure $\mu_0$
as an integral of the measures $\mu_{(a,b)}$ with respect to a
certain measure $\nu$ on $\bR^\times\times\bR$. More explicitly,
this means that we must find a measure $\nu$ such that, for every
$F\in C_c(\bR^3\times\bR^\times)$, we have
\[
(2\pi)^{-1}\int_{\bR^3\times\bR^\times}
F\pare{\zeta,az^2+b,2az,a}\,d\zeta\, dz\,
d\nu(a,b)=(2\pi)^{-2}\int_{\bR^3\times\bR^\times}
F(\eta_1,\eta_2,y_1,y_2)\,d\eta\, dy.
\]
Using the change of variables $\xi_1=\zeta$, $\xi_2=az^2+b$,
$x_1=2az$, $x_2=a$, we see that the measure $\nu$ making the above
equation work is given by $d\nu(a,b)=\pi^{-1}\cdot\abs{a}\,
da\,db$.

\sbr

Finally, we plug the result into \eqref{e:conj1}. Using
Proposition \ref{p:realization} to find the image of
$-\Lb_1^2-\Lb_2^2$ in the \unirreps of $\Gb=\exp\fgb$, we obtain
the following result:
\begin{equation}\label{e:formula-inhomog}
N(\la,H)\sim \frac{1}{\pi}\int_{\bR^\times\times\bR}
\abs{a}\,da\,db\cdot N\left(\la,-\frac{\dd^2}{\dd
z^2}+\pare{az^2+b}^2\right),
\end{equation}
or equivalently, in the form which will appear naturally in the
result of the direct variational argument at the end of the
section,
\begin{equation}\label{e:finhom2}
N(\la, H)\sim \frac{1}{\pi}\int_0^\infty d a
\int_{-\infty}^{+\infty} db \cdot N\left(\la,-\frac{d^2}{d
z^2}+\pare{\sqrt{a}z^2+b}^2\right).
\end{equation}
Change the variables $z\mapsto \la^{-1/2}z, \xi_2\mapsto
\la^{1/2}\xi_2$, $a\mapsto \la^{3/2} a$:
\begin{equation}\label{e:finhom3}
N(\la, H)\sim \ka(H)\la^{7/2}, \end{equation} where
\begin{equation}\label{e:defka}
\ka(H)=\frac{1}{\pi}\int_0^\infty da \int_{-\infty}^{+\infty} db
\cdot N\left(1,-\frac{d^2}{d z^2}+\pare{\sqrt{a}z^2+b}^2\right).
\end{equation}
The convergence of the integral \eq{e:defka} will be proved in
Appendix B.

\subsection{Verification of \eq{e:finhom3}: the general outline}\label{ss:Vgenoutline}
We modify the variational method of calculation of the leading
term of the spectral asymptotics applied to  Scr\"odinger
operators without magnetic fields in \cite{L1} and \cite{L3}, and
 Scr\"odinger operator in 2D, with homogeneous magnetic
potentials, in \cite{L4}. The proofs are based on several standard
variational lemmas; see, e.g., Appendix in \cite{L2} for the list.
We will construct open sets $\Uj\subset \Uoj$ and non-negative
functions $\psi_j\in C^1$, $j=0,\pm 1, \pm 2,\ldots$, such that
\begin{equation}\label{e:partunm}
\cup_j\overline{U_{j}}=\rn;
\end{equation}
\begin{equation}\label{e:partun}
U_{j}\cap U_{k}=\emptyset;
\end{equation}
\begin{equation}\label{e:partun0}
\forall\ x\in\rn,\quad {\rm card} \{j\ |\ x\in \Uoj\}\le C;
\end{equation}
\begin{equation}\label{e:partun1}
\supp~ \psi_j\subset \Uoj;
\end{equation}
\begin{equation}\label{e:partun2}
|\nabla \psi_j(x)|^2\le  \eps\Pst(x)^2,
\end{equation}
where $C, \eps$ are independent of $j$, and  $\eps$ is small, and
\begin{equation}\label{e:partun3}
\sum_j \psi^2_j(x)=1,\quad x\in\rn.
\end{equation}
Let $N_D(\lambda, A, U)\ $stand for the counting function of the
Dirichlet problem for an operator $A\ $on an open set $U$. Since
 the open sets $\Uj$ do not intersect
(pairwise), the standard variational argument gives an estimate
from below
\begin{equation}\label{e:bbelow}
N(\la; H(a)+V)\ge \sum_{j}N_D(\lambda, H(a)+V, \Uj).
\end{equation}
To derive a similar upper bound, we use
 the IMS-localization formula (see e.g. \cite{CFKS})
\begin{equation}\label{e:IMS}
H(a)+ V= \sum_{j}\psi_{j}(H(a) + V)\psi_{j} - \sum_{j\geq 0}
|\nabla\psi_{j}|^{2}.
\end{equation}
Applying \eq{e:partun2} and \eq{e:partun0} to last term on the
RHS, we obtain
\begin{equation}\label{e:Hbound}
H(a)+ V\ge  \sum_{j}\psi_{j}(H(a) + V- C\eps(\Pst)^2)\psi_{j}.
\end{equation}
Employing the standard variational considerations, we derive from
\eq{e:Hbound} and \eq{e:partun1} an estimate from above
\begin{equation}\label{e:babove}
N(\la; H(a)+V)\le \sum_{j}N_D(\lambda, H(a)+V-C\eps(\Pst)^2,
\Uoj).
\end{equation}
We will construct the partition of unity, depending on $\la$, so
that $\eps=\eps(\la)\to +0$ as $\la\to+\infty$. Then on the
strength of the lower bound obtained in \cite{HM}, \cite{MN}:
\begin{equation}\label{e:lbound}
H(a)+V\ge c(\Pst)^2,
\end{equation}
 one should expect that the leading term of the
asymptotics of the RHS in \eq{e:babove} will not change if the
error term $C\eps(\Pst)^2$ is omitted. Hence, the lower bound
\eq{e:bbelow} and upper bound \eq{e:babove} may have the same
leading term of the asymptotics (which is necessary if we want to
derive the leading term of the asymptotics of $N(\la, H(a)+V)$
from \eq{e:bbelow} and \eq{e:babove}) provided the difference
$\Uoj\setminus \Uj$ is relatively small w.r.t. $\Uj$. However,
\eq{e:partun2} implies that this difference cannot be arbitrary
small. Moreover, \eq{e:partun2} prevents the sets $\Uj$ themselves
from being too small, which is desirable on the next step: the
``freezing" of the electric and magnetic tensor at a point
$\xj\in\Uj$. Indeed, the larger the size of $\Uoj$, the larger is
the error due to the freezing.

\sbr

These general arguments explain that one should construct the
partition of unity carefully, by taking into account the competing
requirements described above. The next general observation
explains the difference between non-degenerate cases (as in
\cite{CdV}), and degenerate ones, from the point of view of this
scheme. If $\Uj$ is close to the degeneration set, where
$\Pso:=V^{1/2}+\sum_{jk}|b_{jk}|^{1/2}$ is zero (call this set
$\Sigma$), then one cannot freeze $V$ and $B$ completely, and
still have a small error; but one can freeze $V$ and $B$ in the
directions tangent to $\Sigma$. When convenient, one may freeze
$V$ and $B$ in directions close to tangent ones. In  the special
case which we consider in this section, $\Sigma=\{(x_1, x_2)\ |\
x_1^2=x_2\}$, and in a vicinity of $\xj\in\Sigma$, it is
convenient to freeze the dependence on $x_2$. Notice also that in
order to justify the freezing procedure, one has to work with the
magnetic potential $a$, and not with the magnetic tensor only.
This requires a careful choice of a gauge transformation.

\sbr

The final general remark is: if it is known (or anticipated) that
one or the other term in \eq{e:bbelow} and/or \eq{e:babove} does
not contribute to the leading term of the asymptotics, then one
may omit the corresponding term in \eq{e:bbelow}, and obtain a
crude upper bound for the term in \eq{e:babove}. In particular,
from the classical results for the Laplacian in a bounded domain,
it follows that a ball $U_0=B(0, M_0)$ does not contribute to the
leading term of the asymptotics, if $M_0=M_0(\la)$ grows slower
than any power of $\la$. (Here, $B(y, r)$ denotes an open ball of
radius $r$ centered at $y$.) We will use this observation when we
construct the partition of unity. Further, in all examples, which
we considered, in the weak and intermediate degeneration cases,
only the $\Uj$ which are separated from $\Sigma$ contribute to the
leading term of the asymptotics; on the contrary, in the strong
degeneration case, only the $\Uj$ adjacent to $\Sigma$ do. The
same holds for many other classes of degenerate operators (see
\cite{L1}-\cite{L4}). Hence, in the former case, it suffices to
obtain a crude estimate from above for terms with $\Uj$ adjacent
to $\Sigma$, and in the latter, for terms with $\Uj$ separated
from $\Sigma$. In the special case which we consider in this
section, the degeneration is strong.

\subsection{Construction of a partition of unity}\label{ss:AVpartunity}
We adopt the standard policy on constants: unless otherwise
stated, $C, c, C_1, c_1$, etc., denote  positive constants
(different at different stages of the proof), which can be chosen
the same for all $\la$ and $j$, and the other parameters they can
depend on. The dependence on $n$ is possible but  $n$ is fixed.
(In fact, in our special case, $n=2$, but we will use $n$ since
the scheme below works in the general case as well.)

We will construct open sets $\Uj\subset \Uoj$ and functions
$\chi_j$, which enjoy the properties
\eq{e:partunm}-\eq{e:partun2}. In addition,
\begin{equation}\label{e:partun4}
0\le \chi_j(x)\le 1,
\end{equation}
and
\begin{equation}\label{e:partun5}
\chi_j(x)=1,\quad x\in\Uj.
\end{equation} Then the functions
 \begin{equation}\label{e:defpsij}
 \psi_j=\frac{\chi_j}{(\sum_k \chi_k^2)^{1/2}}
 \end{equation}
 satisfy \eq{e:partun1}-\eq{e:partun3}.

 For large $M$ and $M_0$, denote by
$\Sigma_M$ the intersection of the union of balls $B(y,
M\Pst(y)^{-1})$, over all $y\in\Sigma$, with the set $\{x\ |\
||x||>M_0\}$, and set
$\Sigma_{-M}=\rn\setminus\overline{\Sigma_M\cup U_0}$. For the
general scheme to work, it is convenient that $M$ tends to
infinity with $\la$ but slower than any power of $\la$ or $|x|$
provided the latter in not in $U_0$; and $M_0$ should grow slower
than any power of $\la$. For instance, $M_0=M_0(\la)=\log\la$ and
$M=M(\la)=\log\log\la$ will do. Notice that we will use the
notation $\Sigma_{\pm M}$ not with the subscript $\pm M$ only but
with $\pm CM$ or $\pm cM$ as well. We will also need a function
which grows slower than any power of $M$; so we introduce
$M_1=\log M$.

 The function $\chi_0$ is a
``smoothed" indicator function of $U_0$: we take $\mu\in
C^\infty(\bR)$ such that $\mu(x)=1, x\le 1$, $\mu(x)=0, x\ge 1$,
and set $\chi_0(x)=\mu(M_1^{-1/2}(|x|-M_0))$, $U_{1,0}=B(0,
M_0+2M_1^{1/2})$. Next, we construct the sets and functions such
that $\Uj\subset \Sigma_{-M}$ (we use negative subscripts for
these functions), and finally, the sets and functions such that
$\Uj\subset \Sigma_{M}$ (we use positive subscripts for these
functions).

Set $\Psi=\Pst-\Pso$. The first useful observation is
\begin{lem}\label{l:estpsi0} For any $c>0$, there exist $C, p>0$
such that for all $x\in\Sigma_{-cM}$,
\begin{equation}\label{e:psidif}
\Psi(x)\le CM^{-p}\Pso(x).
\end{equation}
\end{lem}
\begin{proof}  In the
case $V(x)=0, B(x)=x_1^2-x_2$, we have
\begin{eqnarray*}
\Pst(x)&=&|x_1^2-x_2|^{1/2}+|2x_1|^{1/3}+1,\\
\Pso(x)&=&|x_1^2-x_2|^{1/2},\\
\Psi(x)&=&|2x_1|^{1/3}+1,
\end{eqnarray*}
and for $x\in\Sigma$,
\begin{equation}\label{e:psisg}
\Pst(x)=|2x_1|^{1/3}+1=|2^2x_2|^{1/6}+1.
\end{equation}
Let $x\in \Sigma_{-M}$. If $|x_1|\le 2x_2^{1/2}$, then
 \[\sqrt{|x_1^2-x_2|}\ge
\sqrt{|x_1|-x_2^{1/2}}\sqrt{|x_1|+x_2^{1/2}}\ge
M^{1/2}x_2^{-1/4+1/12}=M^{1/2}x^{1/2}\ge
M^{1/2}(|x_1|^{1/3}+1)/2,\] and if $|x_1|\ge 2x_2^{1/2}$, then
$|x_1^2-x_2|^{1/2}\ge |x_1|\ge M_0^{2/3}(|x_1|^{1/3}+1)/2$.
\end{proof}
Notice that if the degeneration is  complicated, the verification
of \eq{e:psidif} may be more involved.

The estimate \eq{e:psidif} being established, the proof of the
next crucial lemma is the same in any dimension, and for any
(polynomial) $V$ and $B$.
\begin{lem}\label{l:m-bounds}
 For any $c, C>0$, there exist $C_1,  c_1>0$ such that
if $x\in\Sigma_{-cM}$, and $|x-y|\le CM_1\Pso(x)^{-1}$, then
$y\in\Sigma_{-c_1M}$,
\begin{equation}\label{e:slsgmm}
\Pso(x)(1-C_1M^{-p})\le \Pso(y)\le \Pso(x)(1+C_1M^{-p}),
\end{equation}
and each of the functions $b_{jk}, V$ satisfies the estimate
\begin{equation}\label{e:difW}
|W(x)-W(y)|\le C_1M^{-p}\Pso(x)^2.
\end{equation}
\end{lem}
\begin{proof}
The Taylor formula gives
\[
|W(x)-W(y)|\le C_3\sum_{|\al|>
0}|W^{(\al)}(x)|(M_1\Pso(x)^{-1})^{|\al|},\] and by definition of
$\Psi$, for $|\al|>0$,
\[
|W^{(\al)}(x)|\le \Psi(x)^{|\al|+2}.
\]
Applying \eq{e:psidif}, and taking into account that $M_1=\log M$,
we obtain \eq{e:difW}. The estimate \eq{e:slsgmm} and inclusion
$y\in\Sigma_{-c_1M}$ follow  from \eq{e:difW}.
\end{proof}
The estimate \eq{e:slsgmm} implies that
$g_x(z)=M_1^{-2}\Pso(x)^2|z|^2$ is a slowly varying metric on
$\Sigma_{-cM}$ in the sense of Section 18.4 of \cite{H}, uniformly
in $\la\ge \la_0$: if $|x-y|\le M_1\Pso(x)^{-1}$, then
$g_x(\cdot)/C_1\le g_y(\cdot)\le C_1g_x(\cdot)$.  The proof of
Lemma 18.4.4 in \cite{H} provides constants $0< c_1<c_2\le 1$, and
for each $\la$, points $\xj\in \Sigma_{-M}$, $j=-1,-2,\ldots$,
such that the balls $B_{j, c_2}:=B(\xj, c_2M_1\Pso(\xj)^{-1})$
cover $\Sigma_{-M}$, and the balls $B_{j, c_1}$ do not intersect
(pairwise). Set $U_{-1}=B_{-1,c_2}$, and for $j\le -2$, set
$U_{j}=B_{j, c_2}\setminus \cup_{j<k<0} \bar B_{k, c_2`}$. Then
intersect with $\{x\ |\ |x|>M_0\}$ (which is the interior of the
complement to $U_0$), and keep the same notation for the
intersections. By construction, for non-positive $j\neq k, j,k\le
0$, \eq{e:partun} holds.

For $j\le -1$, and a subset $V$ of $\Uj$, denote by $V_c$ the
$1-$neighborhood of $V$ in the metric
$c^2\Pso(\xj)^{-2}|\cdot|^2$, and introduce
$\Uoj=(U_{j})_{M_1^{1/2}}$. From the properties of $B_{j, c_k}$,
$k=1,2,$ and \eq{e:slsgmm}, we derive that there exist $C, c>0$
such that for all $\la\ge \la_0$ and $j\le -1$,
\begin{equation}\label{e:u2}
\Sigma_{-M}\subset \cup_{j\le- 1}\overline{\Uj}\subset \cup_{j\le-
1}\overline{\Uoj}\subset \Sigma_{-cM};
\end{equation}
\begin{equation}\label{e:u3}
{\rm meas}(\Uoj\setminus\Uj)_{M_1^{1/2}}\leq C M_1^{n -
1/2}\Pso(\xj)^{-n};
\end{equation}
\begin{equation}\label{e:u4}
\forall\ x\in \rn,\quad {\rm card}\{ j\ |\ x\in \Uoj\}\leq C;
\end{equation}
\begin{equation}\label{e:u5}
B(\xj, cM_1\Pso(\xj)^{-1})\subset \Uj\subset \Uoj \subset B(\xj,
CM_1\Pso(\xj)^{-1}).
\end{equation}
Fix $\chi\in C^\infty_0(B(0, 2))$ such that $0\le \chi\le 1$,
$\chi(x)=1$ for $|x|\le 1$, and for $\al\in \bZ^n$ and $j\le-1$,
define
\begin{eqnarray*}
\chi_{0j,\al}(x)&=&\chi(2M_1^{-1/2}\Pso(\xj)x+\al),\\
 \chi_{j,\al}&=&\frac{\chi_{0j,\al}}{\sum_{\be\in
\bZ^n}\chi_{0j,\be}},
\end{eqnarray*}
and
\begin{equation}\label{e:defchij}
\chi_j=\sum_{\al: {\rm supp}\, \chi_{j,\al}\cap \Uj\ne \emptyset
}\chi_{j,\al}.
\end{equation}
It is immediate from \eq{e:u2}--\eq{e:u5} that $\Uj, \Uoj$ and
$\chi_j$ satisfy \eq{e:partun}--\eq{e:partun1}, \eq{e:partun4},
and \eq{e:partun5}. In addition,
\begin{equation}\label{e:psi4}
|\dd^\al \chi_j|\le C_\al (M_1^{1/2}\Pso(\xj))^{-|\al|}, \quad
\forall\ \al\in\bZ_+^n,
\end{equation}
where the constants $C,  C_\al$ are independent of $\la$ and $j$.

Now we construct $\chi_j, j\ge 1$. In the case under
consideration, $\Sigma_M$ is a disjoint union of two sets; we
denote the one contained in the half-plane $x_1>0$ by
$\Sigma_M^+$, and the other one by $\Sigma_M^-$. The functions
related to the former will have odd subscripts, and the ones
related to the latter---even subscripts. Clearly, the construction
is the same for both parts; we consider $\Sigma_M^+$. Change the
variable $z=x_1-\sqrt{x_2};$ then $\Sigma^+$ becomes a coordinate
half-axis. Construct points $\xj\in \Sigma^+\setminus U_0$ and
open intervals $U_{0j}\subset \Sigma^+\setminus U_0,
j=1,3,\ldots,$ such that
\begin{equation}\label{e:v1}
U_{0j}\cap U_{0k}=\emptyset, \ \forall\ j\ne k;
\end{equation}
\begin{equation}\label{e:v2}
\cup_j\overline{U_{0j}}\supset \Sigma^+\setminus U_0;
\end{equation}
\begin{equation}\label{e:v3}
B(\xj, cM_1\Pst(\xj)^{-1})\subset U_{0j}\subset B(\xj,
CM_1\Pst(\xj)^{-1}).
\end{equation}
Here $B(a, r)$ denotes a ball in the manifold of degeneration
$\Sigma^+\subset\{z=0, x_2>0\}$. For $j\ge 1$, and a subset $V$ of
$\Sigma^+\setminus U_0$, denote by $V_c$ the $1-$neighborhood of
$V$ in the metric $c^2\Pst(\xj)^{-2}|\cdot|^2$. Evidently,
\begin{equation}\label{e:v4}
\forall\ x_2\in \Sigma^+\setminus U_0, \quad {\rm card} \{j\ |\
x_2\in (U_{0j})_{M_1^{1/2}}\}\le C.
\end{equation}
Further, construct functions $\chi_{0j}, j=1,3,\ldots,$ which
satisfy \eq{e:partun4} and the following conditions
\begin{equation}\label{e:chi01}
\chi_{0j}(x_2)=1,\quad x_2\in U_{0j};
\end{equation}
\begin{equation}\label{e:chi02}
{\rm supp}~ \chi_{0j}\subset (U_{0j})_{M_1^{1/2}};
\end{equation}
\begin{equation}\label{e:chi03}
\forall\ x_2, \quad {\rm card}~\{j\ | x_2\in
(U_{0j})_{M_1^{1/2}}\}\le C;
\end{equation}
\begin{equation}\label{e:chi04}
| \chi_{0j}^{(s)}|\le C_\al (M_1^{1/2}\Pst(\xj))^{-s}, \quad
\forall\ s\in\bZ_+,
\end{equation}
where the constants $C, C_\al$ are independent of $\la$ and $j$.
 In our 2D-example, ${\rm dim}\, \Sigma=1$, and the construction of
open sets $U_{0j}$ and functions $\chi_{0j}$ is straightforward;
if ${\rm dim}\, \Sigma>1$,
 one checks that the metric $M_1^{-2}\Pst(x)^2|\cdot|^2$ is slowly
 varying on $\Sigma\setminus U_0$,
 uniformly in $\la\ge \la_0$, and argues as in the construction of $\xj, \Uj, j\le -1$, above.

Set $\Uj=\{(z, x_2)\ |\ x_2\in U_{0j}, |z|\le M/\Pst(\xj)\}$,
$\Uoj= \{(z, x_2)\ |\ x_2\in (U_{0j})_{M_1^{1/2}}, |z|\le
2M/\Pst(\xj)\}$. Further, take $\mu\in C^\infty(\bR)$, $\mu(z)=1,
z\le 1$, $\mu(z)=0, z\ge 1$, and for $j=1,3, \ldots, $ set
$\chi_j(z, z_2)=\chi_{0j}(x_2)\mu(M_1^{-1/2}\Pst(\xj)(|z|-M))$.
Finally, repeat these constructions with $\Sigma^-$ and
$j=2,4,\ldots$. Clearly, for $j\ge 1$, \eq{e:partun}--
\eq{e:partun1}, \eq{e:partun4} and \eq{e:partun5} hold, and an
analog of \eq{e:psi4} is
\begin{equation}\label{e:psi4p}
|\dd^\al \chi_j|\le C_\al (M_1^{-1/2}\Pst(\xj))^{|\al|}, \quad
\forall\ \al\in\bZ_+^n.
\end{equation}
 Now it is clear that properties  \eq{e:partunm}--
\eq{e:partun1}, \eq{e:partun4} and \eq{e:partun5}  hold for all
$j\in \bZ$, and that uniformly in $j\in\bZ$ and $\la$,
\eq{e:partun2} holds with $\eps=C_1M_1^{-1/2}$.  Applying
\eq{e:defpsij},
 we finish the  construction of open sets $\Uj\subset\Uoj$ and functions $\psi_j$
with properties \eq{e:partunm}--\eq{e:partun3}; hence, we have
bounds \eq{e:bbelow} and \eq{e:babove}. In it easy to show that
there exist $c, C_1>0$ such that
\begin{equation}\label{e:HD0}
\langle (H(a)+V-C\eps(\Pst)u, u\rangle\ge c\langle(-\Delta-
C_1M_0^p)u, u\rangle,\quad u \in C^\infty_0(U_0),
\end{equation}
therefore from the classical estimate
 for the counting function of the
  Dirichlet Laplacian, we obtain a bound
  \begin{equation}\label{e:ND0}
  N_D(\la; H(a)+V-C\eps(\Pst)^2; \Uoj)\le  C_2 M^q \la^2,
  \end{equation}
  where $C_2, q$ are independent of $\la$.  Since $M=\log \la$,  \eq{e:ND0}
  and \eq{e:finhom3}
  imply that the
  $j=0$ terms in \eq{e:bbelow} and
\eq{e:babove} do not contribute to the leading term of the
asymptotics.

\subsection{Estimates for individual terms
in \eq{e:bbelow} and \eq{e:babove}, for $j\le -1$: the case of
general $V$ and $B$}\label{ss:jlem1} The proof is valid for the
general case provided \lemm{l:estpsi0} and \lemm{l:m-bounds} have
been proved already. First, we make a gauge transformation used in
\cite{HM}, \cite{MN}.  Set
\begin{equation}\label{e:axj}
a_{\xj}(x)=\sum_{\al}\frac{(x - \xj)^{\al}}{\al!(|\al|+2)}
(\dd^{\al}B)(\xj)\cdot(x - \xj),
\end{equation}
and write $V\ $in the form
\begin{equation}\label{e:vxj}
V(x)=\sum_{\al}\frac{(x - \xj)^{\al}}{\al!}(\dd^{\al}V)(\xj).
\end{equation}
One easily checks that $d a_{\xj} = B$, therefore there exists a
gauge transformation $u\mapsto \exp (-i\phi_{j})u\ $such that
\begin{equation}\label{e:gauge}
\lan H(a)u, u\ran=\lan H(a_{\xj})\exp (-i\phi_{j})u, \exp
(-i\phi_{j})u\ran,
\end{equation}
and so we may study the $j\le -1$-terms in \eq{e:bbelow} and
\eq{e:babove} with $a_{\xj}$ instead of $a$. On $\Uoj(\supset
\Uj)$, the absolute value of each term in \eq{e:axj} with
$|\al|\ge 1$ admits an upper bound via
\[ C(M_1^{1/2}\Pso(\xj)^{-1})^{|\al|+1}\Psi(\xj)^{|\al|+2};\]
using \eq{e:psidif}, we find an upper bound via $
M_1^{-1}\Pst(\xj)$. Similarly, on the same set, for each of the
terms in \eq{e:vxj} with $|\al|\ge 1$, we obtain an upper bound
via $ M_1^{-2}\Pst(\xj)^2$.

Introduce $a_{\xj}^{0}(x) = \frac{1}{2}B(\xj)\cdot (x - \xj),\
V_{\xj}^{0}(x)= V(\xj),\ H_{\xj}= H(a_{\xj}^{0}) + V_{\xj}^{0}$,
and label by $\eps, \eps_1,\ldots,$ any function of $\la$, which
tends to 0 as $\la\to+\infty$. Using the estimates for $|\al|\ge
1$ terms, which we just obtained, and \eq{e:lbound} and
\eq{e:slsgmm}, we conclude that on $C^\infty_0(\Uoj)$,
\[
(1-\eps_1)H_{\xj} \le H(a)+V-C\eps\Pst(x)^2\le H(a)+V\le
H_{\xj}(1+\eps_1).
\]
It follows that
\begin{equation}\label{e:lowjm}
N_D(\la(1-\eps_2); H_{\xj}, \Uj)\le N_D(\la; H(a)+V, \Uj),
\end{equation} and
\begin{equation}\label{e:upjm}
N_D(\la; H(a)+V-C\eps\Pst(x)^2, \Uoj)\le N_D(\la(1+\eps_2);
H_{\xj}, \Uoj).
\end{equation}
For $j\le -1$, $H_{\xj}\ $is a Schr\"odinger operator with uniform
magnetic potential and constant electric potential. By using
coverings of $\Uj\ $and $\Uoj$ by cubes of size
$M^{1/2}\Pso(\xj)^{-1}$, repeating the proof of Theorem 3.1 in
\cite{CdV}, and using \eq{e:u3} and \eq{e:u5}, we obtain the
estimates
\begin{equation}\label{e:lowjm2}
 (1 - \eps_3)\meas \Uj\cdot
\nu_{B(\xj)}(\la(1 - \eps_3) - V(\xj))\leq N_{D}(\la(1 - \eps_2),
H_{\xj}, \Uj),
\end{equation}
and
\begin{equation}\label{e:upjm2}
N_D(\la(1+\eps_2); H_{\xj}, \Uoj)\le (1 + \eps_3)\meas \Uj~
\nu_{B(\xj)}(\la(1 + \eps_3) - V(\xj)).
\end{equation}
It follows from \eq{e:lbound} and \eq{e:slsgmm} that there exists
$\eps_4$ such that for all $x\in\Uoj$,
\[
\nu_{B(x)}(\la(1-\eps_4)-V(x))\le \nu_{B(\xj)}(\la(1 \pm \eps_2) -
V(\xj))\le \nu_{B(x)}(\la(1+\eps_4)-V(x)),
\]
therefore by gathering \eq{e:lowjm}-\eq{e:upjm2}, we obtain the
estimates
\begin{equation}\label{e:lowm}
(1-\eps_3)\int_{\Sigma_{-M}}\nu_{B(x)}(\la(1-\eps_4)-V(x))dx\le
\sum_{j\le -1}N_D(\la, H(a)+V, \Uj),
\end{equation}
and
\begin{equation}\label{e:upm}
\sum_{j\le -1}N_D(\la, H(a)+V, \Uoj)\le
(1+\eps_3)\int_{\Sigma_{-M}}\nu_{B(x)}(\la(1+\eps_4)-V(x))dx.
\end{equation}

\subsection{Estimate for the contribution of $j\le -1$-terms: the
case $V=0$, $B(x)=x_1^2-x_2$}  Since the degeneration is strong
(this can be guessed or deduced from the general conjectural
formula as in \subsect{ss:inhomog}), the $j\le -1$ terms should
not contribute to the leading term of the asymptotics, and
therefore it suffices to obtain a crude upper bound. In view of
\eq{e:finhom3}, it suffices to show that the RHS in \eq{e:upm}
 is $o(\la^{7/2}).$ We have
\[
\nu_{B(x)}(\la)=(2\pi)^{-1}|x_1^2-x_2|\cdot{\rm card}\{k\ |\
(2k+1)|x_1^2-x_2|<\la \},
\]
therefore
\begin{equation}\label{e:upm2}
\sum_{j\le -1}N_D(\la, H(a), \Uoj)\le C \la \meas U(\la),
\end{equation}
where $U(\la)= \{x\in\Sigma_{-M},\ |x_1^2-x_2|\le \la\}.$ Thus, we
need to show that $\meas U(\la)=o(\la^{5/2})$. Clearly, the
measure of the part of $U(\la)$ below the line $x_2=1$ admits a
bound via $\meas\{x\ |\ x_1^2+|x_2|<\la\}\le C\la^{3/2}$,
therefore it suffices to obtain an upper bound for the measure of
the set $V(\la)$  defined by the following inequalities: $ x_1>0,
x_2>1, |x_1-\sqrt{x_2}|\ge M x_2^{-1/6}, |x_2-x_1^2|\le \la.$ The
last two inequalities taken together imply that on $V(\la)$,
$x_1+\sqrt{x_2}\le \la M^{-1}x_2^{1/6}$, and therefore, $x_2\le
(\la/M)^3$.

$V(\la)$ is the union of $V_1(\la)$ defined by $x_2>1, 0<
x_1<\sqrt{x_2}- M x_2^{-1/6}, x_2-x_1^2\le \la,$ and $V_2(\la)$
defined by $ x_2>1, x_1>\sqrt{x_2}+ M x_2^{-1/6}, x_1^2-x_2\le
\la.$ On $V_2(\la)$, $\sqrt{x_2}+ M x_2^{-1/6}\le x_1\le
\sqrt{\la+x_2}$, therefore
\[
\meas V_2(\la)\le \int_1^{(\la/M)^3} dx_2\cdot
(\sqrt{\la+x_2}-\sqrt{x_2})\le \la^{3/2}+\la \int_\la^{(\la/M)^3}
dx_2\cdot x_2^{-1/2}\le CM^{-3/2}\la^{5/2},
\]
which is $o(\la^{5/2})$. On $V_1(\la)$,
$\sqrt{x_2-\la}<x_1<\sqrt{x_2}- M x_2^{-1/6}$, and essentially the
same calculations give $\meas V_1(\la)=o(\la^{5/2})$ as well.

\subsection{Estimates for individual terms in \eq{e:bbelow} and \eq{e:babove},
for $j\ge 1$} Now \lemm{e:psisg} is no longer available, and we
cannot freeze $V$ and $B$ at $x=\xj$. However, we can freeze the
coordinate along $\Sigma$. In the case under consideration, $V=0$,
and we may assume that $a_1(x)=0, a_2(x)=x_1^3/3-x_1x_2$. By
applying the Taylor formula at $\xj=(y_1, y_2)\in\Sigma_M^+$, we
obtain
\[
x_1^3/3-x_1x_2=a_{21}(x)+a_{22}(x)+a_{\xj}(x),
\]
where $a_{21}(x):=y_1^3/3-y_1y_2-y_1(x_2-y_2)$ can be gauged away,
  $a_{22}(x):=-(x_1-y_1)(x_2-y_2)/2$ is
small on $\Uoj$:
\[
|a_{22}(x)|\le
CM(1+y_2^{1/6})^{-1}M_1(1+y_2^{1/6})^{-1}=O(M_1^{-\infty}),
\]
 and $a_{\xj}(x):=y_1(x_1-y_1)^2=y_2^{1/2}(x_1-y_2^{1/2})^2$.
 Denote $D_j=-i\dd/\dd x_j$, and introduce
\[H_{\xj}=H(a_{\xj})=D_1^2+(D_2+y_2^{1/2}(x_1-y_2^{1/2})^2)^2
=D_z^2+(D_{x_2}+y_2^{1/2}z^2)^2.
\]
 We have \eq{e:lowjm} and \eq{e:upjm} with $V=0$. In $(z, x_2)$
 coordinates,
 $\Uj=(-M\Pst(\xj)^{-1}, M\Pst(\xj)^{-1})\times U_{0j}$, and
 $\Uoj=(-2M\Pst(\xj)^{-1}, 2M\Pst(\xj)^{-1})\times
 (U_{0j})_{M_1^{1/2}}$, therefore,
 \begin{equation}\label{e:upjp}
N_D(\la(1+\eps_2); H_{\xj}, \Uoj)\le N_D(\la(1+\eps_2); H_{\xj},
\bR\times (U_{0j})_{M_1^{1/2}}).
\end{equation}
The operator on the RHS can be realized as the operator on
$(U_{0j})_{M_1^{1/2}}$, with the constant operator-valued symbol,
and the same arguments as in the case of Schr\"odinger operators
without magnetic potential used in \cite{L1}, \cite{L3}, and for
2D-Schr\"odinger operators with homogeneous potentials in
\cite{L4}, show that the RHS admits an upper bound
\begin{equation}\label{e:upjp2}
 N_D(\la(1+\eps_2); H_{\xj},
\bR\times (U_{0j})_{M_1^{1/2}})\le
(2\pi)^{-1}(1+\eps)\int_{(U_{0j})_{M_1^{1/2}}}d x_2 \int_{\bR}
d\xi_2\, N(\la(1+\eps), D_z^2+(\eta_2+y_2^{1/2}z^2)^2).
\end{equation}
Since $\meas ((U_{0j})_{M_1^{1/2}}\setminus U_{0j})=o(\meas
U_{0j})$, we may replace the integration over
$(U_{0j})_{M_1^{1/2}}$ with the one over $U_{0j}$ (and change the
$\eps$ in front of the integral). Further, on $U_{0j}$,
$|x_2-y_2|\le M\Pst(\xj)^{-1}$ is small, therefore we may replace
$y_2$ with $x_2$ (and change the $\eps$ in the integrand). After
that we sum up all the terms $j\ge 1$, and  taking into account
the symmetry of $B(x)$ and $\Sigma$, derive from \eq{e:upjm},
\eq{e:upjp}, and \eq{e:upjp2}
\begin{equation}\label{e:upp}
\sum_{j\ge 1}N_D(\la; H(a)+V-C\eps\Pst(x)^2, \Uoj)\le
\frac{1+\eps}{\pi}\int_{M_0}^\infty dx_2\int_{\bR} d\xi_2
N(\la(1+\eps), D_z^2+(\eta_2+x_2^{1/2}z^2)^2).
\end{equation}
We replace $M_0$ with 0 (the inequality remains valid, of course),
and change the variables $z\mapsto (\la(1+\eps))^{-1/2}z,\\
\xi_2\mapsto (\la(1+\eps))^{1/2}\xi_2$, $x_2\mapsto
(\la(1+\eps))^3 x_2$; the RHS becomes
$\ka(H)\la^{7/2}(1+\eps))^{9/2}$, where $\ka(H)$ is given by
\eq{e:defka}. Since the sum of the terms with the labels $j\le 0$
grows slower than $\la^{7/2}$, we conclude that
\begin{equation}\label{e:nlaup}
N(\la, H)\le (1+\eps)\ka(H)\la^{7/2}
\end{equation}
(with a new $\eps$.) To finish the proof of \eq{e:finhom3}, we
need to obtain a lower bound of the form
\begin{equation}\label{e:nlaum}
N(\la, H)\ge (1-\eps)\ka(H)\la^{7/2}.
\end{equation}
Construct a partition of unity $\phi_1^2+\phi_2^2=1$ on $\bR$ with
properties $0\le \phi_j\le 1, \phi_1(z)=1$ for $|z|\le
M/(2\Pst(\xj))$, $\phi_2(z)=1$ for $|z|\ge M/(\Pst(\xj))$, and
$|\phi'_j(z)|\le CM^{-1}/\Pst(\xj)$. If $\eps_2$ decays slower
than $\eps_1$ and any negative power of $M_1$, then we can use
this partition of unity to obtain an upper bound
\begin{equation}\label{e:l1}
N_D(\la(1-\eps_3); H_{\xj}, \bR\times U_{0j})\le
N_D(\la(1-\eps_2); H_{\xj}, U_j)+ N_D(\la(1-\eps_2); H_{\xj},
\{|z|\ge M/(2\Pst(\xj))\}\times  U_{0j}).
\end{equation}
Similarly to \eq{e:upp}, we derive
\begin{equation}\label{e:l2}
\sum_{j\ge 1}N_D(\la(1-\eps_3); H_{\xj}, \bR\times U_{0j})\ge
\frac{1-\eps_4}{\pi}\int_{M_0}^\infty dx_2\int_{\bR} d\xi_2\cdot
N(\la(1-\eps_4), D^2_z+(\xi^2+x^{1/2}z^2)^2).
\end{equation}
At the end of the section, we show that
\begin{equation}\label{e:l3}
\sum_{j\ge 1}N_D(\la(1-\eps_3); H_{\xj},\{|z|\ge
M/(2\Pst(\xj))\}\times  U_{0j})=o(\la^{7/2}),
\end{equation}
and the argument used to derive \eq{e:finhom3} from \eq{e:finhom2}
and prove the convergence of the integral \eq{e:defka}, give
\begin{equation}\label{e:l4}
\int_0^{M_0}dx_2\int_{\bR}d\xi_2\cdot N(\la;
D^2_z+(\xi_2+x_2^{1/2}z^2)^2)=o(\la^{7/2}).
\end{equation}
Estimates \eq{e:l1}--\eq{e:l4} taken together give an estimate
from below
\begin{equation}\label{e:l5}
\sum_{j\ge 1}N_D(\la(1-\eps_2); H_{\xj}, U_j)\ge
\frac{1-\eps_4}{\pi}\int_0^\infty dx_2\int_{\bR} d\xi_2\cdot
N(\la(1-\eps_4), D^2_z+(\xi^2+x_2^{1/2}z^2)^2)+o(\la^{7/2}).
\end{equation}
Changing the variables $z\mapsto (\la(1-\eps_4))^{-1/2}z,
\xi_2\mapsto (\la(1-\eps_4))^{1/2}\xi_2$, $x_2\mapsto
(\la(1-\eps_4))^3 x_2$, we obtain \eq{e:nlaum}.

The proof of \eq{e:finhom3} is complete.

\subsection{Proof of \eq{e:l3}} $H_{\xj}$ is a \sch operator
with magnetic tensor $2y_2^{1/2}z$ (and $V=0$), where $y_2$ is a
constant (the second coordinate of the fixed point $\xj$). For
this operator, $\Pso(z, x_2)=|2y_2^{1/2}z|^{1/2}$, and $\Psi(z,
x_2)=|2y_2^{1/2}|^{1/3}+1$. On the set $\{|z|\ge
M/(2\Pst(\xj))\}\times  U_{0j}$, we have $|z|\ge cM y_2^{-1/6}$;
therefore
\[
\Pso(z, x_2)/\Psi(z, x_2)\ge c_2y^{1/12}|z|^{1/2}\ge c_3M^{1/2},
\]
whence we can repeat the argument used in \S\ref{ss:AVpartunity}
and \S\ref{ss:jlem1} for the localization of $H(a)+V$ on
$\Sigma_{-M}$ and derive an upper bound of the form
\[
N_D(\la(1-\eps_3); H_{\xj},\{|z|\ge M/(2\Pst(\xj))\}\times
U_{0j})\le C \meas\, \Uoj \int\int_{|z|\ge cMy_2^{-1/6}}
\nu_{2y_2^{1/2}z}(\la)dz dy_2.
\]
Clearly, $ \nu_{2y_2^{1/2}z}(\la)=2y_2^{1/2}|z|\cdot\# \{k\ |\
(2k+1)2y_2^{1/2}|z|<\la\} $ is bounded, and vanishes unless
$c_2My_2^{-1/6}/2\le |z|\le \la y_2^{-1/2}$. Summing w.r.t. $j\ge
1$, we obtain that the LHS in \eq{e:l3} admits an upper bound via
the measure of the set $\{(z, y_2)\ |\ y_2>1, \ c_2My_2^{-1/6}\le
|z|\le \la y_2^{-1/2}\}$. It is easy to see that this measure is
$O((\la/M)^3)=o(\la^{7/2}).$

\section{Classical formulas}\label{s:examples}

\subsection{The classical Weyl formula}\label{ss:Weyl} In this subsection we
consider a \sch operator $H=H(0)+V$ with zero magnetic potential
and quasi-homogeneous weakly degenerate electric potential $V(x)$.
We aim at the following result.
\begin{prop}\label{p:weyl}
Assume that $V(x)\geq 0$ is a quasi-homogeneous polynomial, and
that the integral on the RHS of the classical Weyl formula
\eqref{e:weyl1} converges. Then the formula for the leading term
of the asymptotics of $N(\la,H)$ provided by our conjecture
coincides with the classical Weyl formula.
\end{prop}

\sbr

We begin the proof of Proposition \ref{p:weyl} with some general
estimates that will also be useful for us later. Observe that
Fubini's theorem easily implies that
\[
\meas\cur{ (\xi,x)\in\bR^{2n} \st \norm{\xi}^2+V(x)\leq\la }
<\infty \quad\text{for all }\la>0
\]
if and only if
\[
\meas\cur{ x\in\bR^{n} \st V(x)\leq\la } <\infty \quad\text{for
all }\la>0.
\]
Now suppose $V(x)$ is quasi-homogeneous of weight
$\ga=(\ga_1,\dotsc,\ga_n)$, $\ga_j\in\bQ$, $\ga_j>0$. Introduce
the ``quasi-dilation'' $\de_\la:\bR^n\to\bR^n$,
\[
\de_\la : x=(x_1,\dotsc,x_n)\longmapsto \pare{\la^{\ga_1}x_1,
\dotsc, \la^{\ga_n}x_n}.
\]
Thus $V(\de_\la x)=\la\cdot V(x)$ for all $\la>0$, $x\in\bR^n$.
\begin{lem}\label{l:A}
If $\al=(\al_1,\dotsc,\al_n)$, $\al_j\in\bZ$, $\al_j\geq 0$, then
\[
(\dd^\al V)(\de_\la x) =
\begin{cases}
\la^{1-\langle\al,\ga\rangle}\cdot \dd^\al V(x) & \quad \text{if}\ \ \pair{\al}{\ga}\leq 1, \\
0 & \quad \text{if}\ \  \pair{\al}{\ga}>1,
\end{cases}
\]
where $\pair{\al}{\ga}=\al_1\ga_1+\dotsb\al_n\ga_n$.
\end{lem}
\begin{proof}
This follows by induction on $\abs{\al}=\al_1+\dotsc+\al_n$,
observing that if $P(x)$ is a polynomial such that $P(\de_\la
x)=\la^\be\cdot P(x)$ for some $\be>0$, then
\[
\frac{\dd P}{\dd x_j}(\de_\la x) = \la^{-\ga_j}\cdot
\frac{\dd}{\dd x_j} \pare{P(\de_\la x)} = \la^{\be-\ga_j}\cdot
\frac{\dd P}{\dd x_j}(x).
\]
\end{proof}
\begin{lem}\label{l:B}
Let $P(x),Q(x)\geq 0$ be measurable functions on $\bR^n$, finite
almost everywhere, such that $\meas\{P(x)\leq\la\}$ and
$\meas\{Q(x)\leq\la\}$ are finite for each $\la>0$. Fix
$\ga=(\ga_1,\dotsc,\ga_n)$, $\ga_j\in\bR$, $\ga_j>0$, and define
the quasi-dilation $\de_\la$ as before. Assume that there exist
constants $0<\eps<1$ and $C>0$ such that $P(\de_\la x)=\la\cdot
P(x)$ and $Q(\de_\la x)\leq C\la^\eps\cdot Q(x)$ for all
$x\in\bR^n$. Then
\[
\meas\{P(x)+Q(x)\leq\la\} \sim \meas\{P(x)\leq\la\} =
\la^{\abs{\ga}}\cdot \meas\{P(x)\leq 1\} \quad \text{as }\la\to
+\infty,
\]
where $\abs{\ga}=\ga_1+\dotsb+\ga_n$.
\end{lem}
The proof is given in Appendix B.

\sbr

\begin{cor}\label{c:Weyl}
Let $V\geq 0$ satisfy the assumptions of Proposition \ref{p:weyl}.
Introduce the functions $\Phi^*(x)$, $\Psi^*(x)$ as in Section
\ref{s:conj}:
\[
\Phi^*(x) = \sum_{\al\geq 0} \abs{\dd^\al V(x)}^{1/2}, \quad
\Psi^*(x) = \sum_{\al\geq 0} \abs{\dd^\al V(x)}^{1/(2+\abs{\al})}.
\]
Then
\[
\lim\limits_{\la\to +\infty}
\dfrac{\meas\{\Phi^*(x)\leq\la\}}{\meas\{V(x)^{1/2}\leq\la\}} = 1
= \lim\limits_{\la\to +\infty}
\dfrac{\meas\{\Psi^*(x)\leq\la\}}{\meas\{V(x)^{1/2}\leq\la\}}.
\]
\end{cor}
\begin{proof}
Define
\[
\Phi(x) = \sum_{\al>0} \abs{\dd^\al V(x)}^{1/2}, \quad \Psi(x) =
\sum_{\al>0} \abs{\dd^\al V(x)}^{1/(2+\abs{\al})}.
\]
If $\eps=\max_{1\leq j\leq n} \{1-\ga_j\}<1$, then by Lemma
\ref{l:A}, we have
\[
\Phi\pare{\de_{\la^2}x}\leq \la^\eps\cdot \Phi(x)
\quad\text{and}\quad \Psi\pare{\de_{\la^2}x}\leq \la^\eps\cdot
\Psi(x)
\]
for all $\la>1$, and since $V\pare{\de_{\la^2}x}^{1/2}=\la\cdot
V(x)^{1/2}$, Lemma \ref{l:B} applies.
\end{proof}

\sbr

We see that we are in case (b) of Conjecture 2. Moreover, the
normalizing factor $G_2(\la)/G_1(\la)$ is asymptotically equal to
$1$. It remains to compute the quotient measure $\mu_0$. The first
step here is very similar to the first step in the computation
presented in \S\ref{ss:inhomog}.

\sbr

Let $\fg$ be the Lie algebra associated to the operator
$H=-\De+V$. Choose a basis $P_1(x),\dotsc,P_K(x)$ of the vector
space spanned by all mixed partial derivatives (of all orders) of
$V(x)$, {\em not including} $V(x)$ itself. Then we obtain a basis
$\cur{L_1,\dotsc,L_n,v,p_1,\dotsc,p_K}$ of $\fg$ such that the
tautological representation of $\fg$ maps $L_j\mapsto\dd/\dd x_j$,
$v\mapsto\sqrt{-1}\cdot V(x)$, $p_k\mapsto\sqrt{-1}\cdot P_k(x)$.
As in \S\ref{ss:inhomog}, we use this basis to identify $\fg^*$
with $\bR^{n+K+1}$. Write
$\fh=\Span_\bR\cur{v,p_1,\dotsc,p_K}\subset\fg$, and define a
linear functional $f_0:\fg\to\bR$ by $f_0(L_j)=0$ ($1\leq j\leq
n$), $f_0(v)=V(0)=0$, $f_0(p_k)=P_k(0)$ ($1\leq k\leq K$). It
follows from Proposition \ref{p:realization} that the induced
representation $\rho_{f_0,\fh}$ coincides with the tautological
representation of $\fg$ (at least up to changing the action of
$L_1$ and $L_2$, so that $H$ is replaced by a gauge-equivalent
operator). In particular, $\rho_{f_0,\fh}$ is irreducible (by
Theorem \ref{t:reps}), so $\fh$ is a polarization of $\fg$ at
$f_0$. Let $\Om=G\cdot f_0$ be the coadjoint orbit of $f_0$ in
$\fg^*$, and let $\mu_\Om$ be the corresponding Kostant measure.
Again, to simplify notation, we will implicitly identify $\mu_\Om$
with its extension by zero to all of $\fg^*$.

\sbr

By Proposition \ref{p:measure}, the orbit $\Om$ is parameterized
by the map $\varphi:\bR^{2n}\to\bR^{n+K+1}\cong\fg^*$,
$(\xi,x)\mapsto\pare{\xi_1,\dotsc,\xi_n,V(x),P_1(x),\dotsc,P_K(x)}$,
and, moreover, we have $\mu_\Om=(2\pi)^{-n}\cdot\varphi_*(d\xi
dx)$. To find $\mu_0$, we fix $F\in C_c(\fg^*)$ and consider the
integral
\[
\int_{\bR^{2n}} F\pare{ \la^{-1}\xi_1,\dotsc,\la^{-1}\xi_n,
\la^{-1} V(x), \la^{-1} P_1(x),\dotsc, \la^{-1} P_K(x)} d\xi dx.
\]
The change of variables $\xi=\la\xi'$, $x=\de_\la x'$ transforms
the integral into
\[
\la^{n+\abs{\ga}}\cdot\int_{\bR^{2n}} F\pare{
\xi'_1,\dotsc,\xi'_n, V(x'), \la^{-1} P_1(\de_\la x'),\dotsc,
\la^{-1} P_K(\de_\la x')} d\xi' dx'.
\]
By assumption, $\meas\{x'\in\bR^n\big\vert V(x')\leq\la\}$ is
finite for every $\la>0$. Now Fubini's theorem implies that for
every compact subset $\cK\subset\bR^{n+1}$, we have
\[
\meas\cur{ (\xi',x')\in\bR^{2n} \st (\xi',V(x'))\in\cK }<\infty.
\]
On the other hand, for a fixed $x'\in\bR^n$, we have $\la^{-1}
P_k(\de_\la x')\to 0$ as $\la\to +\infty$ by Lemma \ref{l:A}. Thus
Proposition \ref{p:measlim} applies, and we see that there exists
a weak limit
\[
\mu_0=\lim\limits_{\la\to +\infty}
\la^{-(n+\abs{\ga})}\cdot\mu_\la,
\]
given by
\begin{equation}\label{e:mu0}
\int_{\fg^*} F\, d\mu_0 = (2\pi)^{-n}\cdot \int_{\bR^{2n}}
F\pare{\xi'_1,\dotsc,\xi'_n,V(x'),0,\dotsc,0} d\xi dx
\end{equation}
for every $F\in C_c(\fg^*)$. In particular, $\mu_0$ is supported
on the annihilator $\fa^\perp\subset\fg^*$ of the ideal
\[\fa=\Span_\bR\{p_1,\dotsc,p_K\}\subset\fg.\] But $\fg/\fa$ is
clearly abelian, so all \unirreps of the corresponding Lie group
$\exp(\fg/\fa)$ are $1$-dimensional and are obtained by
exponentiating linear functionals on $\fg/\fa$. Given
$f\in\fa^\perp\cong(\fg/\fa)^*$, the image of
$H^\circ=-(L_1^2+\dotsb+L_n^2)-\sqrt{-1}\cdot v\in\cU(\fg)_\bC$
under the representation corresponding to $f$ is the scalar
$f(L_1)^2+\dotsc+f(L_n)^2+f(v)$. Thus we see that our conjecture
produces the following formula:
\[
N(\la,H)\sim (2\pi)^{-n}\cdot \meas\cur{ (\xi',x')\in\bR^{2n} \st
\norm{\xi'}^2+V(x')\leq \la } \quad \text{as } \la\to +\infty,
\]
which coincides with the classical Weyl formula \eqref{e:weyl1}.

\begin{rem}\label{r:weyl}
We can now justify the remark of \S\ref{ss:pushforward} of the
Introduction. Indeed, note that the measure $\mu_0$ given by
\eqref{e:mu0} can also be obtained simply by taking the
pushforward of the canonical measure $\mu_\Om$ with respect to the
projection $\fg^*\to\fg^*$ induced by the map $\fg\to\fg$ taking
$L_j\mapsto L_j$, $v\mapsto v$ and $p_k\mapsto 0$ for $1\leq k\leq
K$. Furthermore, if we simply apply Conjecture 1 (i.e., formula
\eqref{e:conj1}) with the normalization factor $\kappa$ equal to
$1$, then we immediately recover Weyl's formula. Also, this
``derivation'' of Weyl's formula does not rely on the assumption
that $V$ is quasi-homogeneous. On the other hand, the computation
presented above (via the scaling construction) uses the
quasi-homogeneity assumption in an essential way.
\end{rem}

\subsection{Colin de Verdi\`ere's formula}\label{ss:CdV}
In this subsection we consider a \sch operator $H=H(a)$ with zero
electric potential and polynomial magnetic tensor
$B(x)=\brak{b_{jk}(x)}_{j,k=1}^n$ such that $\norm{B(x)}\to\infty$
as $\norm{x}\to\infty$. We suppose that $B$ is quasi-homogeneous
of weight $\ga=(\ga_1,\dotsc,\ga_n)$ and use the notation
$\de_\la$ in the same way as was done in \S\ref{ss:Weyl}. Our goal
is to show that the formula provided by our conjecture coincides
with the formula \eqref{e:CdV} of Colin de Verdi\`ere (see also
\cite{CdV}, Th\'eor\`eme 4.1).

\sbr

We begin by observing that since $\norm{B(x)}\to\infty$ as
$\norm{x}\to\infty$, we have, a fortiori, $\meas\cur{ \sum_{j,k}
\abs{b_{jk}(x)}^{1/2}\leq\la }<\infty$ for all $\la>0$, and
therefore the same type of argument as in \S\ref{ss:Weyl} above
applies. Thus we are in case (b) of Conjecture 2 (in fact, there
is no degeneration at all), and, moreover, the normalizing factor
$G_2(\la)/G_1(\la)$ is asymptotically equal to $1$.

\sbr

Next, let $\fg$ be the Lie algebra corresponding to $H$. We have
the elements $L_j=\dd/\dd x_j + \sqrt{-1}\cdot a_j(x)\in\fg$. We
will denote by $B_{jk}$ the elements of $\fg$ that map to
$\sqrt{-1}\cdot b_{jk}(x)$ under the tautological representation
of $\fg$. We also choose a basis $P_1(x),\dotsc,P_K(x)$ of the
space of polynomials spanned by all mixed partial derivatives of
all {\em nonzero} orders of all the $b_{jk}(x)$, and we let
$p_k\in\fg$ denote the element mapping to $\sqrt{-1}\cdot P_k(x)$
under the tautological representation of $\fg$. As before, up to
replacing $H$ with a gauge equivalent \sch operator, the
tautological representation of $\fg$ lifts to a \unirrep of
$G=\exp\fg$, which corresponds to the coadjoint orbit
$\Om\subset\fg^*$ parameterized by the map
$\varphi:\bR^{2n}\to\fg^*$ defined by
\[
\varphi(\xi,x)(L_j)=\xi_j,\quad
\varphi(\xi,x)(B_{jk})=b_{jk}(x),\quad \varphi(\xi,x)(p_k)=P_k(x).
\]
Note that the polynomials $b_{jk}(x)$ are not linearly
independent, so that $\{L_1,\dotsc,L_n,B_{jk},p_k\}$ spans $\fg$
but is not a basis of $\fg$. Still, the formula above is
meaningful and correct. Moreover, the canonical measure $\mu_\Om$
is given by
\[
\mu_\Om = (2\pi)^{-n}\cdot \varphi_*(d\xi dx).
\]

\sbr

Now, using quasi-homogeneity as in \S\ref{ss:Weyl}, we see that
there exists a weak limit
\[
\mu_0=\lim_{\la\to +\infty} \la^{-n-\abs{\ga}}\cdot\mu_\la =
(2\pi)^{-n}\cdot \psi_*(d\xi dx),
\]
where $\psi:\bR^{2n}\to\fg^*$ is defined by
\[
\psi(\xi,x)(L_j)=\xi_j,\quad \psi(\xi,x)(B_{jk})=b_{jk}(x),\quad
\psi(\xi,x)(p_k)=0.
\]
In particular, $\mu_0$ is supported on the annihilator $\fa^\perp$
of the ideal $\fa\subset\fg$ spanned by $p_1,\dotsc,p_K$. Let
$\fgb=\fg/\fa$, and let $\Lb_j$, $\Bb_{jk}$ denote the images of
the elements $L_j$, $B_{jk}$ in $\fgb$. Thus
$[\Lb_j,\Lb_k]=\Bb_{jk}$ in $\fgb$, and the elements $\Bb_{jk}$
are central in $\fgb$.
\begin{rem}\label{r:CdV}
As in Remark \ref{r:weyl}, we note that the measure $\mu_0$ could
have been obtained from the canonical measure $\mu_\Om$ as the
pushforward with respect to a suitable projection map
$\fg^*\to\fa^\perp$. The computation of the ``quotient measure''
$\nu$ presented below does not make use of the quasi-homogeneity
assumption. In fact, if we take $Q=\supp\mu_0/G$, then
\eqref{e:sd0} becomes precisely the formula of Colin de Verdi\`ere
\eqref{e:CdV}.
\end{rem}

\sbr

The rest of the computation is based on the following two lemmas.
The proof of Lemma \ref{l:measurable} is a simple exercise in
linear algebra, and the proof of Lemma \ref{l:spectral} is a
straightforward application of Fubini's theorem. Therefore both
proofs are omitted.
\begin{lem}\label{l:measurable}
Let $O(n,\bR)$ denote the Lie group of orthogonal $n\times n$
matrices over $\bR$, and let $\fo(n,\bR)$ be its Lie algebra,
i.e., the space of all skew-symmetric $n\times n$ matrices over
$\bR$. Then there exists a measurable map
\[
\fC:\fo(n,\bR)\rar{} O(n,\bR)
\]
such that for every $B\in\fo(n,\bR)$, the matrix $\fC(B)\cdot
B\cdot \fC(B)^T$ is block-diagonal, with $r$ blocks along the
diagonal of the form $\matr{0}{b_j}{-b_j}{0}$, where $b_1\geq
b_2\geq \dotsb\geq b_r>0$, and which has all the other elements
equal to $0$.
\end{lem}
\begin{lem}\label{l:spectral}
Let $r,d\geq 1$ be integers, and let $b_1,\dotsc,b_r>0$ be fixed
real numbers. Then we have
\[
\int_{\bR^d} \#\Bigl\{m_1,\dotsc,m_r\in\bZ_+\st \sum_{j=1}^r
(2m_j+1) b_j\leq\la-\norm{\xi}^2\Bigr\}\,d\xi =
\abs{v_d}\cdot\sum_{m_j\geq 0} \Bigl(\la-\sum_{j=1}^r
(2m_j+1)b_j\Bigr)_+^{d/2}.
\]
\end{lem}

\mbr

Let $2r$ be the maximal possible rank of $B(x)$, as $x$ runs
through all of $\bR^n$. The set of points $x\in\bR^n$ where $B(x)$
has rank less than $2r$ is defined by a collection of polynomial
equations, and hence has Lebesgue measure $0$. Thus it can be
ignored both in our conjectural asymptotic formula, and in the
formula of Colin de Verdiere. Put $U=\{x\in\bR^n\st \rank
B(x)=2r\}$. For each $x\in U$, we let $b_1(x)\geq\dotsb\geq
b_r(x)>0$ denote the positive eigenvalues of $\sqrt{-1}\cdot
B(x)$. Also, let $\fC$ be the measurable map provided by Lemma
\ref{l:measurable}, and set
$\cL=\Span_{\bR}\{L_1,\dotsc,L_n\}\subset\fg$. For every point
$x\in U$, we have a basis $\{L_k'(x)\}_{k=1}^n$ of $\cL$ defined
by $L_k'(x)=\sum_{j=1}^n p_{kj}(x) L_j$, where $p_{kj}(x)$ are the
entries of the matrix $\fC(B(x))$. Thus we can define a new map
$\psi':\bR^{2n}\to\fg^*$ by
\[
\psi'(\xi,x)(L_j'(x))=\xi_j,\quad
\psi'(\xi,x)(B_{jk})=b_{jk}(x),\quad \psi'(\xi,x)(p_k)=0.
\]
Note that $\psi'$ is no longer a polynomial map, but it is still
measurable. Since the matrix $\fC(B(x))$ is orthogonal, it is easy
to check that $\mu_0=(2\pi)^{-n}\cdot \psi'_*(d\xi dx)$, and also
$L_1'(x)^2+\dotsb+L_n'(x)^2=L_1^2+\dotsb+L_n^2$ in $\cU(\fg)$ for
all $x\in U$.

\sbr

By abuse of notation, we view $\psi'$ as a map
$\bR^{2n}\to\fgb^*$, where $\fgb=\fg/\fa$ is defined as above.
Given $\xi\in\bR^n$, let us write $\xi'=(\xi_1,\dotsc,\xi_{2r})$
and $\xi''=(\xi_{2r+1},\dotsc,\xi_n)$. It is easy to check, using
the definition of $\psi'$ and Propositions \ref{p:orbit} and
\ref{p:measure} that, for fixed $\xi''\in\bR^{n-2r}$ and $x\in U$,
the map $\psi_{\xi'',x}:\bR^{2r}\to\fgb^*$ defined by
$\psi_{\xi'',x}(\xi')=\psi'(\xi',\xi'',x)$ parameterizes a single
coadjoint orbit $\Theta_{\xi'',x}\subset\fgb^*$ whose associated
canonical measure $\mu_{\xi'',x}$ is given by
\[
\mu_{\xi'',x}=(2\pi)^{-r}\cdot b_1(x)^{-1}\dotsm b_r(x)^{-1}\cdot
\pare{\psi_{\xi'',x}}_*(d\xi_1\dotsm d\xi_{2r}).
\]
Now it is clear that as the quotient measure space $(Q,\nu)$ we
can take $Q=\bR^{n-2r}\times U$ and
\[
\nu(\xi'',x)=b_1(x)\dotsm b_r(x)\cdot d\xi_{2r+1}\dotsm d\xi_n
dx_1\dotsm dx_n.
\]
Moreover, the representation of $\fgb$ corresponding to the orbit
$\Theta_{\xi'',x}$ can be realized in $L^2(\bR^r)$ in such a way
that $L_j\mapsto \frac{\dd}{\dd y_j}$ for $1\leq j\leq r$,
$L_j\mapsto\sqrt{-1}\cdot b_j(x)\cdot y_j$ for $r+1\leq j\leq 2r$,
and $L_j\mapsto\sqrt{-1}\cdot\xi_j$ for $2r+1\leq j\leq n$. The
rest follows immediately from Lemma \ref{l:spectral} and the
formula
\[
\sg\left( -\De_y+\sum_{j=1}^r b_j^2 y_j^2 + \norm{\xi''}^2 \right)
= \left\{ \sum_{j=1}^r (2m_j+1) b_j+\norm{\xi''}^2 \,\Big\vert\,
m_1,\dotsc,m_r\in\bZ_+ \right\}
\]
which is straightforward from the well-known formula for the
spectrum of the 1D \sch operator corresponding to the harmonic
oscillator.

\subsection{Weak and Intermediate Degeneration Cases: A General Theorem}
For $U\subset \rn$, define $G(\la, \Pst, U)$ by \eq{e:defG}, and
set
\begin{equation}
F(\lambda, a, V, U) = \int_{U}\nu_{B(x)}(\lambda - V(x)) dx.
\end{equation}
For  $d>>1$, define $W_{-d}=\bigl\{x\ |\ \Psi(x)\le
d^{-1}\Pst(x)\bigr\}$, $W_d=\rn\setminus W_{-d}$. The next theorem
states that if for large $d$, the sets $W_d$ are relatively small
in a certain sense, then the classical Weyl formula (when there is
no magnetic field) or Colin de Verdi\`ere's formula hold but with
the integration over $W_{-d}$, where $d$ is a fixed positive
number.
\begin{thm}\label{thmwint}
Let \eq{e:necsuf} hold, and let
 there exist a function $M\to +\infty\
$as $\lambda\to +\infty$, such that for any $C>0$,
\begin{equation}\label{e:genwint1}
G( C\lambda, \Pst, W_M)= o(G(\lambda, \Pst, \rn))\quad {\rm as}\
\lambda\to +\infty.
\end{equation}
Then
\begin{enumerate}[(a)]
\item there exists a function $\eps=\eps(\lambda)\to 0\ $as
$\lambda\to +\infty\ $ such that for any $d>0$,
\begin{equation}\label{e:genwint2}
(1 - \eps)F(\lambda(1 - \eps), a, V, W_{-d})\leq N(\lambda, H(a) +
V) \leq (1 + \eps)F(\lambda(1 + \eps), a, V, W_{-d});
\end{equation}
\item if, in addition, for any $\eps=\eps(\lambda)\to 0\ $as
$\lambda\to +\infty$,
\begin{equation}\label{e:genwint3}
F(\lambda(1+\eps), a, V, W_{-d})\sim F(\lambda, a, V, W_{-d}),
\quad {\rm as}\  \lambda\to +\infty,
\end{equation}
then
\begin{equation}\label{e:genwint4}
 N(\lambda, H(a)+V)\sim F(\lambda, a, V, W_{-d})
 \quad {\rm as}\  \lambda\to +\infty.
 \end{equation}
\end{enumerate}
\end{thm}
\begin{proof}
As $M\to+\infty$, the set $W_{-M}$ satisfies the crucial condition
\eq{e:psidif} for $\Sigma_{-M}$, which was needed to construct
sets $\Uj, \Uoj$ and functions $\chi_j, j\le -1$ (see
\eq{e:defchij}) satisfying properties
\eq{e:partun}-\eq{e:partun2}, and \eq{e:partun4}-\eq{e:partun5}.
Set $\chi_0=1-\sum_{j\le -1}\chi_j$, and define $
\psi_j=\chi_j/\sum_{j\le 0}\chi_j^2,$ $U_0=\{x\ |\ \chi_0(x)>0\}$,
$U_{1,0}=W_{CM}$, where $C>0$ sufficiently large so that $\supp
\chi_0\subset W_{CM}$ (if $W_M=\emptyset$ for large $M$, this step
is not needed, and only $j\le -1$ are involved). Then all the
conditions \eq{e:partunm}-\eq{e:partun3} are satisfied, and we
have the estimates \eq{e:bbelow} and \eq{e:babove} for $N(\la;
H(a)+V)$. The $j\le -1$ terms are treated exactly as
\subsect{ss:jlem1}, and we derive estimates \eq{e:lowm} and
\eq{e:upm}. In view of \eq{e:genwint1}, to finish the proof, it
suffices to shaw that there exist constants $ C, C_1, C_2$ such
that
\begin{equation}\label{e:estpst}
N_D(\la; H(a)+V, W_M)\le C\cdot G(C_1\la, \Pst, W_{C_2M}).
\end{equation}
With $\rn$ in place of $W_M$ and $W_{C_2M}$, this is the estimate
obtained in \cite{MN}. Since the proof in \cite{MN} is obtained by
using an appropriate partition of unity, it can be used to derive
\eq{e:estpst}, and we are done.
\end{proof}
If the potentials are quasi-homogeneous, and $\Pso$ has no zero
outside the origin, then for large $d$, $U_d=\emptyset$ (hence,
\eq{e:genwint1} holds). Due to the quasi-homogeneity, $F(\lambda,
a, V, \rn)$ is of the form $const\cdot \la^p$, hence
\eq{e:genwint3} is satisfied, and \eq{e:genwint4} gives
\begin{equation}\label{e:genwint5}
 N(\lambda, H(a)+V)\sim F(\lambda, a, V, \rn)
 \quad {\rm as}\  \lambda\to +\infty.
 \end{equation}
 In particular, we recover the classical Weyl formula  or
Colin de Verdi\`ere's formula (for no degeneration case). If
$\Pso$ has zeroes outside the origin, we apply the following
result, whose proof is given in Appendix B.
\begin{lem}\label{l:weakdeg}
Let $V(x)$ and $B(x)$ be quasi-homogeneous of the same weight
$\ga=(\ga_1,\dotsc,\ga_n)$. Assume also that
\begin{equation}\label{e:sufWeyl}
\meas\cur{ (\xi,x)\in\bR^{2n} \st \Psi_0(x)^2+\norm{\xi}^2 \leq\la
} <\infty \quad\forall\,\la>0.
\end{equation}
Then
\begin{equation}\label{e:dagger}
\meas\cur{ \Psi_0(x)^2+\norm{\xi}^2\leq\la } \,\sim\, \meas\cur{
\Psi^*(x)^2+\norm{\xi}^2\leq\la } \quad\text{as } \la\to +\infty,
\end{equation}
and \eqref{e:genwint1} holds.
\end{lem}
\begin{thm}\label{t:weyl}
Let $V$ be quasi-homogeneous, $a=0$, and  \eq{e:sufWeyl} hold.

Then the classical Weyl formula holds. \end{thm}
\begin{proof} Under assumptions of the theorem, the conclusion of
\lemm{l:weakdeg} holds, and from \theor{thmwint}, we know that the
estimate \eq{e:genwint2} holds as well. If $V$ is non-degenerate,
then $W_{-d}=\rn$, and since there exist $\ka(V)$ and $p>0$ such
that $F(\la; 0, V, \rn)\sim \ka(V)\la^p$, the estimate
\eq{e:genwint2} gives the classical Weyl formula. If $V$ is
degenerate, we notice that if $\de\to 0$ sufficiently slowly, and
$C$ is sufficiently large, then outside a ball $B(0, C/\eps)$, the
$\de_\la$-invariant $\eps$-neighborhood of the degeneration set
$\Sigma$, call it $CW_\eps$, will contain $W_d$. Hence, $F(\la; 0,
V, W_d)\le F(\la; 0, V, CW_\eps)+O(\eps^{-1}\la^n)$ as $\eps\to
0$, uniformly in $\la>1$ . But  $F(\la; 0, V, CW_\eps)=o(\la^p)$.
Hence, on the LHS and RHS of \eq{e:genwint2}, we may replace
$W_{-d}$ with $\rn$, and add $o(\la^p)$. After that, the classical
Weyl formula is immediate.
\end{proof}

\section{Examples of Schr\"odinger operators with degenerate homogeneous
potentials}\label{s:example2}

\subsection{An example with strong degeneration}\label{ss:strong}
We consider a $2$-dimensional Schr\"odinger operator with zero
electric potential and magnetic tensor $b(x)=b_{12}(x)=x_1^k
x_2^l$, where $k>l\geq 1$ (the case $k=l$ is considered in
\S\ref{ss:intermed}, and leads to substantially different
computations). The corresponding Lie algebra $\fg$ has a natural
basis of the form
\[
\{L_1,L_2\}\cup\{L^{(ij)}\}_{0\leq i\leq k,\ 0\leq j\leq l},
\]
where, in the tautological representation of $\fg$, we have:
\[
L_1\mapsto \frac{\dd}{\dd x_1}+\sqrt{-1}\cdot a_1(x),\quad
L_2\mapsto \frac{\dd}{\dd x_2}+\sqrt{-1}\cdot a_2(x),
\]
\[
\frac{\dd a_2}{\dd x_1}-\frac{\dd a_1}{\dd x_2}=x_1^k x_2^l, \quad
L^{(ij)}\mapsto \sqrt{-1}\cdot x_1^i x_2^j.
\]
(Since this is a faithful representation, these formulas determine
the commutation relations between the basis elements $L_1$, $L_2$,
$L^{(ij)}$ of $\fg$.) Our \sch operator is the image of the
element $H^\circ=-L_1^2-L_2^2\in\cU(\fg)_\bC$ under this
representation. As usual, we assume that the tautological
representation of $\fg$ lifts to a \unirrep of $G=\exp\fg$, and we
write $\Om\subset\fg^*$ for the corresponding coadjoint orbit, and
$\mu_\Om$ for its canonical measure, extended by zero to a measure
on $\fg^*$. Our goal is the following
\begin{prop}\label{p:strong}
We have an asymptotic formula
\begin{equation}\label{e:2dstrong}
N(\la,H) \sim \frac{1}{\pi}\cdot \la^{(l+k+2)/2l}\cdot
\int_0^{+\infty} dx_2\, \int_{-\infty}^{+\infty} d\xi_2\,
N\left(1, -\frac{d^2}{dy^2}+ \Bigl( x_2^l\cdot
\frac{y^{k+1}}{k+1}+\xi_2 \Bigr)^2 \right) \quad\text{as }\la\to
+\infty.
\end{equation}
\end{prop}

\sbr

The direct variational proof of \eq{e:2dstrong} is an evident
modification of the proof in \sect{s:ex1}. The goal of this
subsection is to derive Proposition \ref{p:strong} from Conjecture
2. We leave it to the reader to check that $H$ exhibits strong
degeneration, so that we are in case (a) of Conjecture 2. The
computations are rather straightforward, and are easier, for
example, than the ones involved in Lemma \ref{l:estimates} below
whose detailed proof is given in Appendix B.

\sbr

In what follows, the notation will become quite messy. To simplify
it a little bit, we choose an ordering of our basis of $\fg$ as
follows:
\[
L_1,\ L_2,\ L^{(k,l)},\ L^{(k-1,l)},\dotsc,\ L^{(1,l)},\
L^{(0,l)},
\]
\[
L^{(k,l-1)},\dotsc,L^{(0,l-1)},\dotsc,L^{(k,0)},\dotsc,L^{(0,0)}.
\]
Using this ordering, we identify $\fg^*$ with
$\bR^{2+(k+1)(l+1)}$. This particular identification will be used
throughout our computation.

\sbr

It follows from Propositions \ref{p:orbit} and \ref{p:measure}
that the orbit $\Om$ is parameterized by the map
$\varphi:\bR^4\to\bR^{2+(k+1)(l+1)}\cong\fg^*$ defined by
\begin{eqnarray*}
\varphi(\xi_1,\xi_2,x_1,x_2)&=&(\xi_1,\xi_2,x_1^k
x_2^l,x_1^{k-1}x_2^l,\dotsc,x_1 x_2^l,x_2^l, \\
&& x_1^k x_2^{l-1},\dotsc,x_2^{l-1},\dotsc,x_1^k,\dotsc,x_1,1),
\end{eqnarray*}
and we have
\[
\mu_\Om = (2\pi)^{-2}\cdot\varphi_*(d\xi_1 d\xi_2 dx_1 dx_2).
\]

\sbr

Now we fix a continuous function $F$ with compact support on
$\fg^*\cong\bR^{2+(k+1)(l+1)}$, and for $\la>0$, we consider the
integral
\[
\la^{-2-1/l}\cdot\int_{\fg^*} F\,d\mu_\la =
(2\pi)^{-2}\cdot\la^{-2-1/l}\cdot \int_{\bR^4}
F(\la^{-1}\xi_1,\la^{-1}\xi_2,\la^{-1}x_1^k x_2^l,\dotsc,\la^{-1})
d\xi_1 d\xi_2 dx_1 dx_2.
\]
We make the following change of variables: $\xi_1'=\la\xi_1$,
$\xi_2'=\la^{-1}\xi_2$, $x_1'=x_1$, $x_2'=\la^{-1/l}\cdot x_2$.
This kills the factor $\la^{-2-1/l}$, and the integral becomes
\[
(2\pi)^{-2}\int_{\bR^4} F\bigl(\xi'_1,\xi'_2,(x'_1)^k (x'_2)^l,
\dotsc, (x'_2)^l, \la^{-1/l}\cdot (x'_1)^k
(x'_2)^{l-1},\dotsc,\la^{-1}\bigr) d\xi'_1 d\xi'_2 dx'_1 dx'_2.
\]
It is straightforward to prove that, for a fixed $R>0$, the set of
points $(\xi'_1,\xi'_2,x'_1,x'_2)\in\bR^4$ such that
$\norm{\xi'}\leq R$ and $\abs{{x'_1}^k {x'_2}^l}\leq R$, \ldots,
$\abs{{x'_2}^l}\leq R$ has finite Lebesgue measure (this is where
the assumption that $k>l$ is used in a crucial way). Now
Proposition \ref{p:measlim} implies that the limit as $\la\to
+\infty$ of the above integral is equal to
\[
(2\pi)^{-2}\int_{\bR^4} F\bigl(\xi'_1,\xi'_2,(x'_1)^k (x'_2)^l,
\dotsc, (x'_2)^l, 0, 0, \dotsc, 0\bigr) d\xi'_1 d\xi'_2 dx'_1
dx'_2,
\]
which implies that there exists a weak limit
\[
\mu_0 = \lim\limits_{\la\to +\infty} \la^{-2-1/l}\cdot\mu_\la =
(2\pi)^{-2}\cdot \psi_*(d\xi_1 d\xi_2 dx_1 dx_2),
\]
where $\psi:\bR^4\to\bR^{2+(k+1)(l+1)}\cong\fg^*$ is given by
\[
\psi(\xi_1,\xi_2,x_1,x_2)=(\xi_1, \xi_2, x_1^k x_2^l, x_1^{k-1}
x_2^l, \dotsc, x_1 x_2^l, x_2^l, 0, 0, \dotsc, 0).
\]

\sbr

Next we must compute the ``relevant'' coadjoint orbits, i.e.,
those that contribute to the RHS of our conjectural formula. We
begin by observing that $\mu_0$ is supported on the annihilator
$\fa^\perp\subset\fg^*$, where $\fa\subset\fg$ is the ideal
spanned by all the elements $L^{(ij)}$ with $0\leq i\leq k$,
$0\leq j\leq l-1$. Let us identify the quotient algebra
$\fgb=\fg/\fa$ explicitly. It has a basis which we will denote
(abusing notation) by $\bigl(\frac{\dd}{\dd
w},w^{k+1},w^k,\dotsc,w,1\bigr)$. The commutation relations are
obvious from the notation, and the quotient map $\fg\to\fg/\fa$ is
determined by $L_1\mapsto\frac{\dd}{\dd w}$,
$L_2\mapsto\frac{w^{k+1}}{k+1}$, and $L^{(jl)}\mapsto w^j$ for
$0\leq j\leq k$. Once again, we will use this chosen basis to
identify $\fgb^*$ with $\bR^{k+3}$. By abuse of notation, we will
identify $\mu_0$ with its restriction to
$\fgb^*\cong\fa^\perp\subset\fg^*$. Explicitly, we have
\[
\mu_0=(2\pi)^{-2}\cdot \psi_*(d\xi_1 d\xi_2 dx_1 dx_2),
\]
where $\psi:\bR^4\to\bR^{k+3}$ is given by
\[
\psi(\xi_1,\xi_2,x_1,x_2) = (\xi_1,(k+1)\xi_2, x_1^k x_2^l,\dotsc,
x_1 x_2^l, x_2^l).
\]
Note that the set
\[
\psi^{-1}\bigl(\{f\in\fgb^*\big\vert f(1)=0\}\bigr) = \bigl\{
(\xi_1,\xi_2,x_1,x_2)\in\bR^4 \big\vert x_2=0\bigr\}
\]
has zero Lebesgue measure. Hence, for the rest of our computation,
we can throw away the ``non-generic'' elements $f\in\fgb^*$ (i.e.,
such that $f(1)=0$), and concentrate our attention on the map
(which we still denote by $\psi$)
\[
\psi:\bR^3\times\bR^\times\longrightarrow \fgb^*_{gen},
\]
\[
(\xi_1,\xi_2,x_1,x_2)\longmapsto \bigl(\xi_1,(k+1)\xi_2,x_1^k
x_2^l, \dotsc, x_1 x_2^l, x_2^l),
\]
where $\fgb^*_{gen}=\bigl\{f\in\fgb^*\big\vert f(1)\neq 0\bigr\}$,
and $\bR^\times=\bR\setminus\{0\}$. Thus, the ``relevant'' orbits
for us are those which are contained in
$\psi(\bR^3\times\bR^\times)$ (we will see below that this set is
actually stable under the adjoint action on $\fgb^*$).

\sbr

Given any generic $f\in\fgb^*_{gen}$, it is immediate that
$\Span_\bR\{w^{k+1},w^k,\dotsc,1\}$ is a polarization of $\fgb$ at
$f$. Hence it follows from Proposition \ref{p:orbit} orbit that
for any fixed point $(\xi_2,x_2)\in\bR\times\bR^\times$, the
coadjoint orbit in $\fgb^*$ through $\psi(0,\xi_2,0,x_2)$ is
parameterized by the map
\[
(\eta,y)\mapsto \psi'(\eta,\xi_2,y,x_2):=\pare{ \eta,
(k+1)\xi_2+x_2^l y^{k+1}, x_2^l y^k,\dotsc, x_2^l y, x_2^l}.
\]
Moreover, by Proposition \ref{p:measure}, the canonical measure on
this orbit in terms of this parameterization is given by
$(2\pi)^{-1}\,d\eta dy$. Now consider the map
$\theta:\bR^3\times\bR^\times\to\bR^3\times\bR^\times$ given by
\[
\theta(\xi_1,\xi_2,x_1,x_2) = \left( \xi_1,
\xi_2-\frac{x_1^{k+1}x_2^l}{k+1}, x_1, x_2 \right).
\]
It is clear that $\theta$ is a diffeomorphism preserving the
Lebesgue measure; moreover, $\psi=\psi'\circ\theta$. Therefore
$\mu_0=(2\pi)^{-2}\cdot\psi'_*(d\eta d\xi_2 dy dx_2)$. Now we see
that, tautologically, the image of $\psi'$ (which is the same as
the image of $\psi$) is invariant under the coadjoint action on
$\fgb^*$, and the space of coadjoint orbits in the image of
$\psi'$ is naturally parameterized by the points
$(\xi_2,x_2)\in\bR\times\bR^\times$, with the quotient measure
$\nu$ given by $\nu=(2\pi)^{-1}\,d\xi_2 dx_2$. By Proposition
\ref{p:realization}, the representation of $\fgb$ corresponding to
the coadjoint orbit parameterized by $(\xi_2,x_2)$ can be realized
in the space $L^2(\bR^1)$ as follows:
\[\dd/\dd w\mapsto
d/dy,\quad w^{k+1}\mapsto \sqrt{-1}\cdot\pare{(k+1)\xi_2+x_2^l
y^{k+1}},\quad w^j\mapsto \sqrt{-1}\cdot x_2^l y^j\quad\text{for }
0\leq j\leq k.
\]
Finally, we see that our conjectural formula \eqref{e:conj1} gives
the following answer:
\[
N(\la,H) \sim (2\pi)^{-1} \int_{\bR\times\bR^\times} N\left( \la,
-\frac{d^2}{dy^2}+ \Bigl( x_2^l\cdot \frac{y^{k+1}}{k+1}+\xi_2
\Bigr)^2 \right) d\xi_2 dx_2 \quad \text{as }\la\to +\infty.
\]
To reduce this equation to the form \eqref{e:2dstrong}, we make
the change of variables $y\mapsto \la^{-1/2}y$,
$x_2\mapsto\la^{(k+2)/2l} x_2$, $\xi_2\mapsto\la^{1/2}\xi_2$, and
use symmetry to replace $\int_{\bR^\times\times\bR}dx_2 d\xi_2$
with $2\cdot \int_{0}^\infty dx_2 \int_{-\infty}^\infty d\xi_2$.

\subsection{An example with intermediate
degeneration}\label{ss:intermed} In this subsection we consider a
\sch operator $H$ on $L^2(\bR^2)$ with zero electric potential and
magnetic tensor $b(x)=b_{12}(x)=x_1^k x_2^k$ (since we are working
in 2D, the magnetic tensor has only one relevant component). We
assume that $k\geq 1$. We will use Conjecture 2 to derive
\begin{prop}\label{p:intermed}
We have the following asymptotic formula, as $\la\to +\infty$,
\begin{equation}\label{e:2interm}
N(\la,H) \sim \frac{2}{\pi k}\cdot \la^{1+1/k}\log\la\cdot
\sum_{j=0}^\infty (2j+1)^{-1-1/k}.
\end{equation}
\end{prop}
(Note that the equation \eq{e:2interm} can be derived from
\theor{thmwint}.)

 \sbr

The functions \eqref{e:phistar}, \eqref{e:psistar} associated to
the operator $H$ are given by
\[
\Phi^*(x)=\sum_{j,l=0}^k \left| \frac{k!}{(k-j)!}\cdot
\frac{k!}{(k-j)!}\cdot x_1^{k-j} x_2^{k-l}\right|^{1/2}
\]
and
\[
\Psi^*(x)=\sum_{j,l=0}^k \left| \frac{k!}{(k-j)!}\cdot
\frac{k!}{(k-j)!}\cdot x_1^{k-j} x_2^{k-l}\right|^{1/(2+j+l)}.
\]
Define $G_1(\la)$, $G_2(\la)$ as in \S\ref{ss:precise}. We begin
with the following
\begin{lem}\label{l:estimates}
We have
\[
G_1(\la)\sim \frac{8}{k}\cdot\la^{2/k}\log\la \quad\text{and}\quad
G_2(\la)\sim \frac{8(k+1)}{k}\cdot\la^{2/k}\log\la,
\]
and hence $G_2(\la)/G_1(\la)\sim k+1$ as $\la\to +\infty$.
\end{lem}
The proof is given in Appendix B.

\sbr

The lemma implies that we are in the intermediate degeneration
case, and the constant $\kappa$ in our formula \eqref{e:conj1} can
be replaced by $k+1$. Thus we have to compute $\mu_0$ and the
quotient measure $\nu$. We will use the same notation as in
\S\ref{ss:strong} (indeed, the only difference in the setup is
that $k=l$, which, however, does not affect the
representation-theoretic part of our discussion). In particular,
we have the basis $\cur{ L_1,L_2, L^{(k,k)},\dotsc,L^{(0,0)} }$ of
$\fg$ which yields an identification of $\fg^*$ with
$\bR^{2+(k+1)^2}$. Moreover, the orbit $\Om\subset\fg^*$
corresponding to the tautological representation of $\fg$ is
parameterized by the map
\[
\varphi:\bR^4\to\bR^{2+(k+1)^2}\cong\fg^*, \quad
(\xi_1,\xi_2,x_1,x_2)\mapsto (\xi_1,\xi_2,x_1^k x_2^k,\dotsc, 1),
\]
and the canonical measure on $\Om$ is given by
\[
\mu_\Om = (2\pi)^{-2}\cdot \varphi_*(d\xi_1 d\xi_2 dx_1 dx_2).
\]
In Appendix B, we use Proposition \ref{p:limborel} to prove
\begin{lem}\label{l:limit}
There exists a weak limit
\[
\mu_0 = \lim\limits_{\la\to +\infty} \la^{-2-1/k}
(\log\la)^{-1}\cdot\mu_\la.
\]
Moreover, $\mu_0$ is supported on the annihilator
$\fa^\perp\subset\fg^*$ of the ideal
\[
\fa=\Span_{\bR}\cur{ L^{(i,j)} \st i<k \text{ or } j<k }.
\]
In particular, we can view $\mu_0$ as a measure on $\fgb^*$, where
$\fgb=\fg/\fa$ is a $3$-dimensional Heisenberg algebra. If
$(\eta_1,\eta_2,y)$ are the coordinates on $\fgb^*$ defined by the
images of the elements $L_1,L_2,L^{(k,k)}\in\fg$ in $\fgb$, then
\[
\mu_0 = \frac{1}{2\pi^2 k^2} \cdot \abs{y}^{\frac{1-k}{k}}
\,d\eta_1 d\eta_2 dy.
\]
\end{lem}

\sbr

Now the image of $L^{(k,k)}$ in $\fgb$ is central, whence the set
$\{y=0\}\subset\fgb^*$ is invariant under the coadjoint action and
has $\mu_0$-measure zero. Thus it can be ignored in our
computation. The complement of this set is a union of
$2$-dimensional coadjoint orbits parameterized by points
$c\in\bR^\times$. The orbit corresponding to such a $c$ is defined
by $y=c$, and the corresponding canonical measure is
$\mu_c=(2\pi)^{-1}\cdot \abs{c}^{-1} d\eta_1 d\eta_2$ (we may view
$(\eta_1,\eta_2)$ as coordinates on the orbit). We see immediately
that the measure space $(Q,\nu)$ which appears in our conjectural
formula is given by $Q=\bR^\times$, $\nu=(\pi k^2)^{-1}\cdot
\abs{y}^{1/k} dy$. Combining this result with Lemmas
\ref{l:estimates} and \ref{l:limit}, we see that \eqref{e:conj1}
becomes
\begin{equation}\label{e:intermed}
N(\la,H)\sim \frac{k+1}{\pi k^2}\cdot\log\la\cdot
\int_{-\infty}^\infty \abs{y}^{1/k} N\left( \la, -\frac{\dd^2}{\dd
x^2} + y^2 x^2 \right) dy.
\end{equation}
To reduce this to the formula of Proposition \ref{p:intermed}, we
recall that the eigenvalues of the operator $-\dd/\dd x^2+y^2 x^2$
(for $y\neq 0$) are given by $\{(2j+1)y\big\vert j\in\bZ_+\}$, and
Fubini's theorem implies that for any $\la>0$,
\begin{eqnarray*}
 \int_0^\infty y^{1/k}\cdot N\pare{\la, -\dd/\dd x^2+y^2 x^2} dy
&=& \int_0^\infty y^{1/k} \cdot \#\cur{ j\in\bZ_+ \st
(2j+1)y\leq\la } dy \\
&=& \sum_{j=0}^\infty \, \int_{\cur{y\geq 0\st (2j+1)y\leq\la}}
y^{1/k} dy \\
&=& \frac{k}{k+1} \la^{1+1/k} \sum_{j=0}^\infty (2j+1)^{-1-1/k}.
\end{eqnarray*}
Substituting this into \eqref{e:intermed}, we obtain
\eq{e:2interm}.

\subsection{A three-dimensional example}\label{ss:3dim} As our last example we
briefly discuss the \sch operator $H=-\De+x_1^{2k} x_2^{2l}
x_3^{2p}$ on $L^2(\bR^3)$, where $1\leq p<k\leq l$. Then the set
of degeneration, $\Sigma$, is the union of coordinate planes
$\Sigma^{jk}$, $j\neq k\in\{1,2,3\}$. The scheme of \sect{s:ex1}
should be modified as follows.

After $U_0, U_{1,0}$ and $\chi_0$ are constructed, the sets $\Uj,
\Uoj$ and functions $\chi_j$ should be constructed separately for
$\Sigma_{-M}$, then for $M$-neighborhoods of coordinate planes (in
the metric $\Pst(x)^2 |\cdot|^2$) but outside
$M^{1/2}$-neighborhoods (in the same metric) of coordinate axis,
and finally, for $M^{1/2}$-neighborhoods of coordinate axis. Under
condition $p<l\le k$, the leading contribution comes from the sets
adjacent to the axis $\{l=0, k=0\}$, and we derive
\begin{equation}\label{e:3d}
N(\lambda, H(0)+ V)\sim (2\pi)^{-1}\int_{-\infty}^{+\infty} dy_{3}
 \int_{-\infty}^{+\infty}d\eta_{3}\cdot N(\lambda, \eta_{3}^{2} -
\Delta_{x_{1}, x_{2}}+ y_{3}^{2p}x_{1}^{2k}x_{2}^{2l}).
\end{equation}
 Making substitutions $x_{1}\mapsto
|y_{3}|^{-\al}x_{1},\ x_{2}\mapsto |y_{3}|^{-\al}x_{2},\
\eta_{3}\mapsto |y_{3}|^{\al}\eta_{3}$, where $\al=p/(k+l+1)$, we
obtain
\[ N(\lambda, H(0) + V)\sim \ka(V)\lambda^{(\al+1)/2\al},
\]
where
\begin{equation}\label{e:3Dka1}
\ka(V)=\frac{1}{\pi}\int_0^{+\infty} dy_{3}\cdot y_3^\al
 \int_{-\infty}^{+\infty}d\eta_{3}\cdot N(y_{3}^{-2\al}, \eta_{3}^{2} -
\Delta_{x_{1}, x_{2}}+ x_{1}^{2k}x_{2}^{2l}).
\end{equation}
 Let $0<\la_1\le \la_2\le\cdots$ be the eigenvalues of the
operator $H_1=-\Delta+x_{1}^{2k}x_{2}^{2l}$ in $L_2(\bR^2)$. Then
\begin{eqnarray*}
\ka(V)&=&\frac{1}{\pi}\int_0^{+\infty} dy_{3}\cdot y_3^\al
 \int_{-\infty}^{+\infty}d\eta_{3}\cdot\#\,\{j\ |\
 \eta_3^2<y_{3}^{-2\al}-\la_j\}\\
 &=&\sum_j\frac{2}{\pi}\int_0^{+\infty} dy_{3}\cdot y_3^\al
  ((y_{3}^{-2\al}-\la_j)_+)^{1/2}.
 \end{eqnarray*}
 Under the summation sign, change the variable $y_3\mapsto
 \la_j^{-1/2\al}y_3$:
\begin{equation}\label{e:3Dka2}
\ka(V)=\sum_{j\ge 1} \la_j^{-1/2\al} \ka_1(V),
\end{equation}
where
\[
\ka_1(V)=\frac{2}{\pi}\int_0^1 dy_{3}\cdot y_3^\al
  (y_{3}^{-2\al}-1)^{1/2}=\frac{2}{\pi}\int_0^1
  dy_{3}\cdot(1-y_3^{2\al})^{1/2}=\frac{1}{\pi\al}B(1/2\al, 3/2).
  \]
We need to prove the convergence of the series in \eq{e:3Dka2}. If
$k<l$, then similarly to \eq{e:2dstrong}, we have
\[
N(\la; H_1)\sim \ka(k,l)\la^{(l+k+1)/2l},
\]
where $\ka(k,l)>0$ is independent of $\la$, and hence,
 $\la_j\sim const\cdot j^{2l/(l+k+1)}$,
as $j\to +\infty$. Since $-(2l/(l+k+1))/(2\al)=-l/p<-1$, the
series converges. If $k=l$, then similarly to  \eq{e:2interm},
\[
N(\la; H_1)\sim \ka(l)\cdot \log\la\cdot \la^{1+1/2l},
\]
where $\ka(l)>0$ is independent of $\la$, and hence, $\la_j\sim
const\cdot j^{2l/(2l+1)}/\log j$, as $j\to +\infty$. Since
$-(2l/(2l+1))/(2\al)=-l/p<-1$, the series converges.

\mbr

We now wish to derive \eqref{e:3d} from our conjectural formula.
The computations are completely analogous to, and, in a way,
simpler (due to the absence of the magnetic potential) than those
in \S\ref{ss:strong}; therefore we restrict ourselves to outlining
the main steps. First one verifies without difficulty that we are
in the strong degeneration case of Conjecture 2. The Lie algebra
$\fg$ associated to $H$ has a basis
\[\cur{ L_1,L_2,L_3}\cup\cur{ L^{(jst)} \st 0\leq j\leq 2k,\,
0\leq s\leq 2l,\, 0\leq t\leq 2p},\] where, in the tautological
representation of $\fg$, we have $L_j\mapsto\dd/\dd x_j$
($j=1,2,3$) and $L^{(jst)}\mapsto \sqrt{-1}\cdot x_1^j x_2^s
x_3^t$. The canonical measure on the orbit corresponding to the
tautological representation of $\fg$ is given by
$\mu_\Om=(2\pi)^{-3}\cdot\varphi_*(d\xi dx)$, where
$\varphi:\bR^6\to\fg^*$ is defined by
\[
\langle \varphi(\xi,x), L_j\rangle = \xi_j \ \  (j=1,2,3), \quad
\langle \varphi(\xi,x), L^{(jst)}\rangle = x_1^j x_2^s x_3^t.
\]
As in \S\ref{ss:strong}, one finds that there exists a weak limit
\[
\mu_0 = \lim\limits_{\la\to +\infty} \la^{-3-1/(2p)} \mu_\la =
(2\pi)^{-3}\cdot \psi_*(d\xi dx),
\]
where $\psi:\bR^6\to\fg^*$ is defined by
\[
\langle \psi(\xi,x), L_j\rangle = \xi_j \ \  (j=1,2,3), \quad
\langle \psi(\xi,x), L^{(j,s,2p)}\rangle = x_1^j x_2^s x_3^{2p},
\quad \langle \psi(\xi,x), L^{(jst)}\rangle = 0 \ \ \text{if }
t<2p.
\]
In particular, $\mu_0$ is supported on the annihilator
$\fa^\perp\subset\fg^*$ of the ideal \[\fa=\Span_\bR\cur{
L^{(jst)} \st t<2p } \subset\fg.\]

\sbr

If we denote by ``bar'' the quotient map $\fg\to\fgb=\fg/\fa$,
then we have a direct sum decomposition of Lie algebras
$\fgb=\bR\cdot \Lb_3\oplus\fg'$, where
$\fg'=\Span_\bR\cur{\Lb_1,\Lb_2,\Lb^{(j,s,2p)}}$ is isomorphic to
the Lie algebra associated to the 2D \sch operator
$-\De_{x_1,x_2}+x_1^{2k} x_2^{2l}$. Using this and the explicit
formula for $\mu_0$, it is easy to see that Conjecture 2(a)
predicts  \eqref{e:3d}.

\appendix\section{Representation theory}\label{s:A}

\subsection{Review of induced representations for nilpotent Lie
groups}\label{ss:induced} The orbit method invented by Kirillov is
based on the notion of an induced representation. Let $G$ be a Lie
group and $H\subset G$ a closed subgroup. To every unitary
representation $U$ of $H$ one can associate a unitary
representation $T=\Ind_H^G(U)$ of $G$. In general, the
construction of $T$ is somewhat complicated, due to the fact that
the quotient space $G/H$ might not possess a $G$-invariant
measure, and the projection map $G\to G/H$ might have no
continuous sections. However, if $G$ is nilpotent, the situation
becomes much simpler. Let us assume that $G$ is a connected and
simply connected Lie group, and $H\subset G$ a closed connected
subgroup ($H$ is then also simply connected). Let $\fh\subset\fg$
be the corresponding Lie algebras. A {\em coexponential basis} for
$\fh$ in $\fg$ is a set of elements $X_1,\dotsc,X_n\in\fg$ such
that the map
\[
\fh\times\bR^n\to G, \quad (\xi,t_1,\dotsc,t_n)\mapsto \exp(t_1
X_1+\dotsb+ t_n X_n)\cdot\exp(\xi)
\]
is a diffeomorphism. Such a basis always exists. It can be
constructed as follows. Considering the adjoint action of $\fh$ on
$\fg/\fh$ and using the fact that $\fh$ is a nilpotent Lie
algebra, we see that there exists a subspace $\fa\supset\fh$ of
$\fg$ such that $[\fh,\fa]\subseteq\fh$, and $\fh$ has codimension
$1$ in $\fa$. It is then immediate that $\fa$ is a subalgebra of
$\fg$, and $\fh$ is an ideal of $\fa$. Next we apply this
construction to $\fa$ in place of $\fh$, and continue inductively.
We find that there is a chain
$\fh=\fa_0\subset\fa_1\subset\dotsb\subset\fa_n=\fg$ of
subalgebras of $\fg$ such that each $\fa_j$ is an ideal of
codimension $1$ in $\fa_{j+1}$. Now if we choose arbitrary
elements $X_j\in\fa_j$ such that $X_j\not\in\fa_{j-1}$ for all
$1\leq j\leq n$, then it is not hard to show, by induction on $n$,
that $\{X_j\}$ is a coexponential basis for $\fh$ in $\fg$.

\sbr

We fix one coexponential basis $\{X_j\}$ for $\fh$ in $\fg$. Let
$X=G/H$ denote the quotient space, let $\pi:G\to X$ be the natural
projection, and let $e:\bR^n\to G$ be defined by $e(t)=\exp(t_1
X_1 + \dotsb + t_n X_n)$. The composition $\bR^n\rar{e} G\rar{\pi}
X$ is a diffeomorphism. Composing the inverse of this
diffeomorphism with the map $e$ defines a smooth section $s:X\to
G$ of the projection $\pi$. We also fix a $G$-invariant positive
Borel measure $\nu$ on $X$; one can prove that since $G$ is
nilpotent, such a measure always exists, and is unique up to a
positive constant multiple. If $U$ is a unitary representation of
$H$ in a (complex) Hilbert space $\cH$, then the induced
representation $T=\Ind_H^G(U)$ of $G$ is realized in the space
$L^2(X,\cH;\nu)$ of square-integrable functions $F:X\to\cH$ with
respect to $\nu$, according to the following explicit formula:
\[
(g\cdot F)(x) = U\pare{ s(x)^{-1}\cdot g\cdot s(g^{-1} x)}
F(g^{-1} x) \qquad (g\in G,\,x\in X).
\]
Here, $g^{-1}x$ denotes the action of the element $g^{-1}\in G$ on
the point $x\in X$. A more detailed discussion of induced
representations can be found, e.g., in \cite{Ki}, \S1 or
\cite{reps}, Chapter V. Note that Kirillov in \cite{Ki} uses the
right homogeneous $G$-space $H\backslash G$ for the construction
of the induced representation, whereas we've used $G/H$, as in
\cite{reps}; of course, the two approaches are equivalent.

\subsection{Description of unitary
representations}\label{ss:proofsreps} In this subsection we prove
the results contained in \S\ref{ss:reps} and
\S\ref{ss:realizations}. We begin with the
\begin{proof}[Proof of Lemma \ref{l:polarization}]
By definition, $\fg(f)$ is the kernel of the form $B_f$, and $\bR
L_0+[\fg,\fg]$ is isotropic with respect to $B_f$ because it is
abelian. Hence $\fg(f)+\bR L_0+[\fg,\fg]$ is an isotropic subspace
of $\fg$ with respect to $B_f$. Let $\fh$ be any maximal isotropic
subspace of $\fg$ containing $\fg(f)+\bR L_0+[\fg,\fg]$. Then
$\fh$ is an ideal of $\fg$ (and a fortiori a subalgebra) because
$\fh$ contains $[\fg,\fg]$. Moreover, we claim that
$[\fh,\fh]\subseteq\fg(f)$. Let $x,y\in\fh$ and $z\in\fg$. Then
Jacobi's identity implies that
\[
\brak{z,[x,y]}=\brak{[z,x],y}+\brak{x,[z,y]}.
\]
But $[z,x]\in\fh$ and $[z,y]\in\fh$ because $\fh$ is an ideal,
whence $f$ annihilates both terms on the RHS of the last equation,
as $\fh$ is isotropic with respect to $B_f$. Since $z\in\fg$ is
arbitrary, we find that $[x,y]\in\fg(f)$, as desired.
\end{proof}

\mbr

Our strategy in what follows will be to give a proof of Theorem
\ref{t:reps} and obtain the explicit formula of Proposition
\ref{p:realization} as a by-product of our discussion. We use the
notation introduced in \S\ref{ss:realizations} (cf. especially the
paragraph preceding the statement of Proposition
\ref{p:realization}).

\sbr

Define $\fa_k=\fh\oplus\bR L_1\oplus\dotsb\oplus \bR L_k$; since
$\fh\supset[\fg,\fg]$, we see that
$\fh=\fa_0\subset\fa_1\subset\dotsb\subset\fa_n=\fg$ is a chain of
ideals of $\fg$, each of codimension $1$ in the next. It follows
easily that $L_1,\dotsc,L_n$ is a {\em coexponential basis} for
$\fh$ in $\fg$, i.e., the map
\[
\varphi : H\times\bR^n \to G, \quad (h,x)\mapsto \exp\pare{-(x_1
L_1+\dotsb+ x_k L_k)}\cdot h,
\]
is a diffeomorphism. To simplify our formulas, we introduce the
following notation: if $x\in\bR^n$, then $x.L=x_1 L_1 + \dotsb +
x_n L_n\in\cL_0$. Since $\fh\supset[\fg,\fg]$, we see that $G/H$
is an abelian Lie group. Hence the composition
\[
\pi:=\text{proj}_2\circ\phi^{-1} : G\xrightarrow{\phi^{-1}}
H\times\bR^n\xrightarrow{\text{proj}_2} \bR^n
\]
is a Lie group homomorphism with kernel $H$, so we obtain an
identification $G/H\cong\bR^n$ (as Lie groups). Moreover, if
$s:\bR^n\to G$ is the map $s(x)=\exp(-x.L)$, then the composition
$\bR^n\rar{s} G\rar{\pi}\bR^n$ is the identity on $\bR^n$. Also,
the $G$-action on $G/H\cong\bR^n$ factors through the action of
$\bR^n$ on itself by translations, so the Lebesgue measure on
$\bR^n$ is $G$-invariant. It follows from \S\ref{ss:induced} above
that $\Ind_H^G(\chi_f)$ can be realized as a unitary
representation of $G$ on $L^2(\bR^n)$ defined by the following
explicit formula:
\[
(g\cdot F)(x) = \chi_f\pare{ s(x)^{-1}\cdot g\cdot s(g^{-1}
x)}\cdot F(g^{-1} x),
\]
where $g\in G$, $F\in L^2(\bR^n)$, $x\in \bR^n$, and $g^{-1}x$
denotes the action of $g^{-1}$ on $x$. From the previous remarks,
we see that $g^{-1}x=-\pi(g)+x$, so we can rewrite our formula as
follows:
\begin{equation}\label{e:action}
(g\cdot F)(x) = \chi_f\pare{ s(x)^{-1}\cdot g\cdot
s(-\pi(g)+x)}\cdot F(-\pi(g)+x).
\end{equation}

\sbr

The proof of Theorem \ref{t:reps}(a) is now a straightforward
computation. Since $\fg=\fh\oplus\Span\{L_1,\dotsc,L_n\}$, we
analyze separately the action of $\fh$ and of the elements $L_j$
on $L^2(\bR^n)$. First, if $h\in H=\exp\fh$, then $\pi(h)=0$, so
we see from \eqref{e:action} that $h$ maps to the operator of
multiplication by the function $x\mapsto\chi_f(s(x)^{-1} h s(x))$.
Differentiating with respect to $h$, we find that an element
$\xi\in\fh$ maps to the operator of multiplication by the function
$i\cdot f\pare{ (\Ad s(x)^{-1})\xi}$. By definition,
$s(x)^{-1}=s(-x)=\exp(x.L)$. Thus, to find the image of $\xi$
explicitly, we must compute $\exp\pare{x_1\cdot\ad
L_1+\dotsb+x_n\cdot\ad L_n}(\xi)$. Now, even though the operators
$\ad L_j$ do not commute in general, each of their iterated
commutators will have the form $\ad Z$ for some
$Z\in[\fg,\fg]\subset\fh$. Hence, we will have $f([Z,\xi])=0$. So
we find that
\begin{equation}\label{e:h}
i\cdot f\pare{ (\Ad s(x)^{-1})\xi} = i\cdot
\sum_{\al_1,\dotsc,\al_n\geq 0} \frac{1}{\al_1!\dotsm\al_n!}\cdot
f\pare{(\ad L_1)^{\al_1}\dotsm(\ad L_n)^{\al_n}(\xi)}\cdot
x_1^{\al_1}\dotsm x_n^{\al_n}.
\end{equation}

\sbr

Next we compute the action of $L_j$ on $L^2(\bR^n)$, for $1\leq
j\leq n$. If $g=s(-y)\in G$, then \eqref{e:action} becomes
\[
(s(-y)\cdot F)(x) = \chi_f\pare{ s(x)^{-1}\cdot s(-y)\cdot
s(y+x)}\cdot F(y+x).
\]
Differentiating with respect to $y_j$ at $y=0$, we obtain, using
the product rule,
\[
(L_j\cdot F)(x) = \frac{\dd F}{\dd x_j}(x) + F(x)\cdot
\frac{\dd}{\dd y_j}\bigg\vert_{y=0} \Bigl(\chi_f\pare{
s(x)^{-1}\cdot s(-y)\cdot s(y+x)}\Bigr).
\]
The second term on the RHS has the form $\sqrt{-1}\cdot
a_j(x)\cdot F(x)$, where $a_j(x)$ is a certain real polynomial
which can easily be found by using the known formula for the
differential of the exponential map $\exp:\fg\to G$ at an
arbitrary point of $\fg$. However, we are not interested in the
exact formula for $a_j(x)$. We have shown that the elements
$\xi\in\fh$ (and in particular $L_0,L_{n+1},\dotsc,L_N$) map to
multiplication operators of the form $\sqrt{-1}\cdot p_\xi(x)$ on
$L^2(\bR^n)$, where $p_\xi(x)$ is a real polynomial, and that
$L_j$ maps to $\frac{\dd}{\dd x_j}+\sqrt{-1}\cdot a_j(x)$ for
$j=1,\dotsc,n$. This proves all the assertions of Theorem
\ref{t:reps}(a) except for the discreteness of spectrum. To
complete the proof, we apply the following
\begin{prop}\label{p:discspec}
A \sch operator \eqref{e:schr} with polynomial potentials has
discrete spectrum if and only if there is no rotation of the
coordinate axes in $\bR^n$ making $V(x)$ and all components of the
magnetic tensor $b_{jk}(x)$ independent of one of the coordinates.
\end{prop}
The ``only if'' direction is trivial; for the ``if'' direction,
see \cite{MN}. To apply the criterion of the proposition, note
that, by definition, $b_{jk}(x)$ is the image of $-\sqrt{-1}\cdot
[L_j,L_k]$, and $V(x)$ is the image of
$-L_{n+1}^2-\dotsb-L_N^2-\sqrt{-1}\cdot L_0$ under the
representation constructed above. We assume that $-\sqrt{-1}\cdot
L_0$ maps to a nonnegative polynomial; then $V(x)$ is the sum of
the image of $-\sqrt{-1}\cdot L_0$ and the images of $(\sqrt{-1}
\cdot L_s)^2$ ($n+1\leq s\leq N$). Each $\sqrt{-1}\cdot L_s$ maps
to a real polynomial. We see that if there exists a rotation of
the coordinate axes with the property defined in the lemma, then
it makes the images of $[L_j,L_k]$ and $L_0$, $L_s$ ($n+1\leq
s\leq N$) separately independent of one of the coordinates. But
the operator $L_j$ maps to $\dd/\dd x_j+\sqrt{-1}\cdot a_j(x)$ for
all $1\leq j\leq n$, so we see that there exists a rotation of the
coordinate axes as above if and only if there exists an element
$D=y_1 L_1+\dotsb+y_n L_n\in\cL_0$ of unit length such that the
image of $D$ commutes with the images of $[L_j,L_k]$, $L_0$, $L_s$
and all their derivatives. Since the operators
$L_0,L_1,\dotsc,L_N$ generate $\fg$, we see that the image of
$[D,\xi]$ in our representation is zero for each $\xi\in\fh$, and
a fortiori, $f([D,\xi])=0$ for all $\xi\in\fh$. But this
contradicts the assumption that $\fh$ is a maximal totally
isotropic subspace with respect to the form $B_f$.

\mbr

Comparing \eqref{e:action} with \eqref{e:action-of-h}, we see that
we have obtained a proof of Proposition \ref{p:realization}.

\mbr

We now prove Theorem \ref{t:reps}(b). Suppose we are given a \sch
operator \eqref{e:schr}, let $\fg$ be the corresponding Lie
algebra and $S\in\cU(\fg)_\bC$ the corresponding sublaplacian. Let
$\fh\subset\fg$ be the subspace spanned by $\sqrt{-1}\cdot
b_{jk}(x)$, $\sqrt{-1}\cdot V(x)$ and all their derivatives. Thus
$\fh$ is an abelian ideal of $\fg$ that contains $[\fg,\fg]$; it
can also be described as the part of $\fg$ consisting of operators
of order $0$. Let $f:\fh\to\bR$ be the linear functional given by
$f(\sqrt{-1}\cdot P(x))=P(0)$, and extend $f$ in an arbitrary way
to all of $\fg$. The construction of induced representation given
above works for the character $\chi_f$ of $H=\exp\fh$ defined by
$\chi_f(\exp\xi)=i\cdot f(\xi)$ without the assumption that $\fh$
is a real polarization at $f$. Thus, we obtain a unitary
representation of $G=\exp\fg$ on the space $L^2(\bR^n)$. Note that
if we set $L_j=\dd/\dd x_j+\sqrt{-1}\cdot a_j(x)$, then $\ad L_j$
acts on $\fh$ as $\dd/\dd x_j$, so comparing \eqref{e:h} to the
standard Taylor's formula, we find that the representation of
$\fg$ we have just defined coincides with the ``tautological
representation'' at least on the subalgebra $\fh\subset\fg$. In
particular, since the components of the magnetic tensor lie in
$\fh$, we see that the image of the sublaplacian $S$ in this
representation is at least gauge equivalent to the original \sch
operator \eqref{e:schr}. Finally, Kirillov's theory \cite{Ki}
implies that the representation of $G$ we have constructed is
irreducible if and only if $\fh$ is a maximal totally isotropic
subspace with respect to $B_f$. The fact that this condition is
equivalent to the discreteness of spectrum of $H$ follows from
Proposition \ref{p:discspec} using the same argument as in the
previous subsection.

\subsection{Parameterization of coadjoint orbits and Kostant measures}
Here we prove the results stated in \S\ref{ss:orbits}. We keep the
notation used in \S\ref{ss:realizations}, \S\ref{ss:orbits} and
\S\ref{ss:proofsreps}.
\begin{proof}[Proof of Proposition \ref{p:orbit}]
Let $H=\exp\fh$; this is an abelian normal subgroup of $G$. The
adjoint action of $G$ on $\fg$ leaves $\fh$ stable, whence $G$
acts on $\fh^*$. Restriction of linear functionals defines a
surjective $G$-equivariant linear map $\pi:\fg^*\to\fh^*$. Since
$H$ is abelian, it acts trivially on $\fh^*$, so the action of $G$
on $\fh^*$ factors through its abelian quotient $G/H$. In
particular, if $\Om=G\cdot f_0\subset\fg^*$ is the coadjoint orbit
of $f_0$, then $\pi(\Om)\subset\fh^*$ is the $G/H$-orbit of
$f_0\big\vert_{\fh}$.

\sbr

On the other hand, since $\fh$ is an ideal of $\fg$, we see that
for any $g\in G$, $\fh$ is a real polarization of $\fg$ at $(\Ad^*
g)(f_0)$. Now Pukanszky's criterion (see \cite{Pu} or \cite{reps},
Chapter VI) implies that if $f\in\fg^*$ is any point such that
$\fh$ is a real polarization at $f$, then the $G$-orbit of $f$
contains $f+\fh^\perp=\{f'\in\fg^*\st
f'\vert_{\fh}=f\vert_{\fh}\}$. With the notation of the previous
paragraph, we see that $\Om=\pi^{-1}(\pi(\Om))$. On the other
hand, by construction, $L_1,\dotsc,L_n$ is a complementary basis
to $\fh$ in $\fg$. So we see that it suffices to prove that the
map $\phi:\bR^n\to\fh^*$ given by
\begin{equation}\label{e:phi}
\pair{\phi(x)}{Y}=\frac{1}{\al_1!\dotsm\al_n!}\cdot f_0\pare{(\ad
L_1)^{\al_1}\dotsm(\ad L_n)^{\al_n}}\cdot x_1^{\al_1}\dotsm
x_n^{\al_n}\quad \text{for all } Y\in\fh
\end{equation}
is a diffeomorphism of $\bR^n$ onto $\pi(\Om)$.

\sbr

Now $\pi(\Om)$ is the $G/H$-orbit of $f_0\big\vert_{\fh}$ in
$\fh^*$. Recall that in \S\ref{ss:proofsreps} we have obtained an
identification of $G/H$ with the abelian Lie group $\bR^n$. More
precisely, if $s:\bR^n\to G$ is the map defined by
$s(x)=\exp(-x.L)$, then the composition $\bR^n\rar{s}G\rar{}G/H$
is a Lie group isomorphism. Hence $\pi(\Om)$ can be parameterized
by
\begin{equation}\label{e:param}
\bR^n\ni x \longmapsto \pare{\Ad^* \exp(-x_1 L_1-\dotsb-x_n
L_n)}(f_0\big\vert_{\fh}).
\end{equation}
But $\Ad^*\exp(-x.L)=\exp\pare{x_1\cdot(\ad
L_1)+\dotsb+x_n\cdot(\ad L_n)}$. Since $\fh$ is an abelian ideal
of $\fg$, the operators $\ad L_1,\dotsc,\ad L_n$ commute on $\fh$.
Thus we immediately see that the parameterization \eqref{e:param}
of $\pi(\Om)$ agrees with the one given by \eqref{e:phi}. Finally,
it remains to check that $\phi$ is a diffeomorphism, i.e., that
the stabilizer of $f_0\big\vert_{\fh}$ in $G/H$ is trivial. Let
$H'$ be the stabilizer of $f_0\big\vert_{\fh}$ in $G$. Since $\fh$
is a maximally isotropic subspace with respect to the form
$B_{f_0}$, it is easy to see that the Lie algebra of $H'$
coincides with $\fh$. Lastly, $H'$ is connected; this follows from
the general theory of unipotent representations of nilpotent Lie
groups (\cite{reps}, Ch. I).
\end{proof}

\begin{proof}[Proof of Proposition \ref{p:measure}]
Fix a point $f\in\Om$. Since $\fh$ is a maximal isotropic subspace
of $\fg$ with respect to $B_f$ and $L_1,\dotsc,L_n$ is a
complementary basis to $\fh$ in $\fg$, there exist elements
$Y_1,\dotsc,Y_n\in\fh$ such that $B_f(L_j,Y_k)=\de_{jk}$ (the
Kronecker delta). Define
\[
Z_j=-L_j+\sum_{k\neq j} f\pare{[L_k,L_j]}\cdot Y_k.
\]
For each $Y\in\fg$, let us write $\nu_Y$ for the tangent vector to
$\Om$ at $f$ generated by $Y$ (via the $G$-action on $\Om$), and
let us write $\eps_Y$ for the function $\Om\to\bR$, $f'\mapsto
f'(Y)$. We have the following
\begin{lem}
If $Y,Z\in\fg$, then $\nu_Y(\eps_Z)=\pair{f}{[Z,Y]}=B_f(Z,Y)$.
\end{lem}
The proof is a completely straightforward computation (see, e.g.,
\cite{reps}, Chapter III). The lemma implies that
$\nu_{Y_j}(\eps_{L_k})=B_f(L_k,Y_j)=\de_{jk}$, and also
\[
\nu_{Z_j}(\eps_{L_k}) = B_f\pare{L_j-\sum_{s\neq j}
f([L_s,L_j])\cdot Y_s,L_k}.
\]
If $j=k$, the last expression is trivially zero. If $j\neq k$, the
expression equals
\[
B_f(L_j,L_k)-f([L_k,L_j])\cdot B_f(Y_k,L_k)=0.
\]

\sbr

Now we wish to compute the vector field $\dd/\dd\xi_j$, $\dd/\dd
x_k$ on $\Om$ corresponding to the coordinates $(\xi,x)$. Note
that, with the notation above, we have $\xi_j=\eps_{L_j}$. Thus,
we see that
\[
\nu_{Y_j}(\xi_k)=\de_{jk} \quad\text{and}\quad \nu_{Z_j}(\xi_k)=0
\quad\text{for all } j,k.
\]
On the other hand, if $\xi\in\bR^n$ is fixed, then by
construction, the map $\bR^n\to\pi(\Om)\subset\fh^*$,
$x\mapsto\pi(\varphi(\xi,x))$, coincides with the map $x\mapsto
\pare{\Ad^*\exp(-x.L)}(\pi(\varphi(\xi,0)))$ (cf. the proof of
Proposition \ref{p:orbit}). Since $Y_k\in\fh$ for all $k$ and $H$
acts trivially on $\fh^*$, we find that $\exp(t\cdot Z_j)$ acts on
$\pi(\Om)\cong\bR^n$ as translation by $t$ in the $x_j$-direction.
It is then immediate that
\[
\nu_{Z_j}(x_k)=\de_{jk} \quad\text{and}\quad \nu_{Y_j}(x_k)=0
\quad\text{for all } j,k.
\]
We conclude that the values of the vector fields $\dd/\dd\xi_j$
and $\dd/\dd x_k$ at the point $f\in\Om$ are given by $\nu_{Y_j}$
and $\nu_{Z_k}$, respectively. Therefore
\[
\om_\Om\pare{\dd/\dd\xi_j,\dd/\dd\xi_k} = B_f(Y_j,Y_k)=0;
\]
\[
\om_\Om\pare{\dd/\dd\xi_j,\dd/\dd x_k} =
B_f(Y_j,Z_k)=B_f(L_k,Y_j)=\de_{jk};
\]
\[
\om_\Om\pare{\dd/\dd x_j,\dd/\dd x_k} = B_f(Z_j,Z_k)=f([L_k,L_j]).
\]
This proves \eqref{e:symplform}, and the formula for the Kostant
measure $\mu_\Om$ follows trivially from this.
\end{proof}

\subsection{Limits of polynomial measures and quotient measures}
In this subsection we prove the propositions stated in
\S\ref{ss:polymeas}.

\begin{proof}[Proof of Proposition \ref{p:measlim}]
For every compact set $K\subset\bR^N$, we have
$\tilde{\phi}_j^{-1}(K),\tilde{\phi}^{-1}(K)\subseteq\phi^{-1}(K')$,
where $K'$ is the projection of $K$ onto the first $N_1$
coordinates, so $\mu_j$ and $\mu$ are regular. Now let $F$ be a
continuous functions on $\bR^N$ that is supported on $K$; it
follows that
\[
\int_{\bR^N} F\,d\mu_j = \int_{\phi^{-1}(K')}
F(\tilde{\phi}_j(x))\,dm(x) \quad\text{and}\quad \int_{\bR^N}
F\,d\mu = \int_{\phi^{-1}(K')} F(\tilde{\phi}(x))\,dm(x),
\]
where $K'$ is defined as before. Since $\phi^{-1}(K')$ has finite
measure, and the integrands above are bounded by a constant
independent of $j$ (namely, the sup norm of $F$), the Dominated
Convergence Theorem applies. Thus, we may pass to the pointwise
limit as $j\to\infty$ inside the first integral, which proves the
proposition.
\end{proof}

\begin{proof}[Proof of Proposition \ref{p:limborel}]
For the first statement, it suffices to show that the limit on the
RHS of \eqref{e:star} exists and is finite for all $F\in
C_c(\bR^N)$. Indeed, it is then obvious that this limit defines a
positive linear functional on $C_c(\bR^N)$, so the Riesz
representation theorem yields the existence and uniqueness of
$\mu_0$.

\sbr

Fix $f\in C_c(\bR^N)$, $\eps>0$, and $R>0$ such that
\[
\Supp(f)\subseteq\cur{x\in\bR^N \st \norm{x}_\infty\leq R-1}.
\]
Now $f$ is uniformly continuous, so there exists $0<\de<1$ such
that if $x,y\in\bR^N$ and $\norm{x-y}_\infty\leq\de$, then
$\abs{f(x)-f(y)}\leq\eps$. Write
$B_{\infty,R}=\{x\in\bR^N\big\vert \norm{x}_\infty\leq R\}$. Then
there exist finitely many points $m_1,\dotsc,m_K\in\bZ^N$ such
that the sets $[\de\cdot m_k,\de\cdot(m_k+\vone))$ cover
$B_{\infty,R}$. Clearly, for each $y_0\in B_{\infty,1}$, the sets
$[\de\cdot m_k+y_0,\de\cdot(m_k+\vone)+y_0)$ will cover
$B_{\infty,R-1}$. Since $B_{\infty,1}$ is uncountable, we can find
$y_0\in B_{\infty,1}$ such that each of the coordinates of the
finitely many points $\cur{\de\cdot
m_k+y_0}_{k=1}^K\cup\cur{\de\cdot(m_k+\vone)+y_0}_{k=1}^K$ is not
contained in the countable set $E$. On the other hand, we still
have $[\de\cdot m_k+y_0,\de\cdot(m_k+\vone)+y_0)\subseteq
B_{\infty,R+1}$. By assumption,
\[
M:=\sup\limits_{\la>0} \mu_\la(B_{\infty,R+1}) \leq
\sup\limits_{\la>0}
\mu_\la\pare{[-(R+1)\cdot\vone,(R+2)\cdot\vone)} < +\infty.
\]
For each $1\leq k\leq K$, put $P_k=[\de\cdot
m_k+y_0,\de\cdot(m_k+\vone)+y_0)$, and
\[
u_k=\max\limits_{x\in\overline{P}_k} F(x), \quad
v_k=\min\limits_{x\in\overline{P}_k} F(x).
\]
By construction, $v_k\leq u_k\leq v_k+\eps$, and
\begin{equation}\label{e:dstar}
\sum_{k=1}^K v_k\cdot\mu_\la(P_k) \leq \int_{\bR^N}
F\,d\mu_\la\leq \sum_{k=1}^K u_k\cdot\mu_\la(P_k).
\end{equation}
This implies that
\[
\left( \limsup\limits_{\la\to +\infty} \int F\,d\mu_\la \right) -
\left( \liminf\limits_{\la\to +\infty} \int F\,d\mu_\la \right)
\leq \eps\cdot\sum_k \pare{\lim\limits_{\la\to +\infty}
\mu_\la(P_k)}\leq \eps\cdot M.
\]
Since $M$ does not depend on $\eps$, this proves the first part of
the proposition.

\sbr

For the second part, we proceed as above, and combine
\eqref{e:dstar} with the inequalities
\[
\sum_k v_k\cdot \nu_0(P_k)\leq \int F\,d\nu_0\leq \sum_k
u_k\cdot\nu_0(P_K),
\]
which show that both $\limsup\limits_{\la\to +\infty} \int
F\,d\mu_\la$ and $\liminf\limits_{\la\to +\infty} \int
F\,d\mu_\la$ are within $\eps\cdot M$ of $\int F\,d\nu_0$, for any
$\eps>0$.
\end{proof}

It remains to discuss Proposition \ref{p:quotient}. This result is
needed for the statements of Conjectures 1 and 2 to make sense.
However, it has little practical significance. Indeed, in all the
examples we have computed, it is obvious that the ``quotient
measure'' $\nu$ exists, since it is easy to write down an explicit
formula for it. Thus we only provide a sketch of the proof of
Proposition \ref{p:quotient}.

\sbr

The main result of \cite{Bon} provides a stratification of $\fg^*$
by locally closed (even in the Zariski topology) subsets, such
that each stratum is equipped with a triple $(r,q,p)$ of
vector-valued functions having the following properties:
\begin{itemize}
\item every stratum is a union of coadjoint orbits, which are
given by the level sets $r=const$;
\item for a fixed $r=r_0$, the pair $(q,p)$ is a global chart for the
corresponding orbit $\Om_{r_0}$, such that the canonical
symplectic form on the orbit is given by $\sum_{j=1}^\ell
dq_j\wedge dp_j$.
\end{itemize}
It is understood that the number $\ell$ depends on the stratum,
and the functions $r,q,p$ take values in $\bR^{n-2\ell}$,
$\bR^\ell$, $\bR^\ell$, respectively. It follows that the Kostant
measure on $\Om_{r_0}$ also has the simple form
\begin{equation}\label{e:kost}
(2\pi)^{-\ell}\cdot dq_1\dotsm dq_\ell dp_1\dotsm dp_\ell
\end{equation}
in the coordinates $(q,p)$. Now, it clearly suffices to prove the
existence and uniqueness of the ``quotient measure'' $\nu$ locally
(here, the word ``locally'' has to be interpreted as ``on every
stratum for a stratification of $\fg^*$ by $G$-invariant Borel
subsets''). But the facts we have just stated imply that locally
the projection $\rho:\fg^*\to\fg^*/G$ is isomorphic to the natural
projection $\pr_1:U\times\bR^{2\ell}\to U$, where $U$ is an open
set in $\bR^{n-2\ell}$, and moreover, the Kostant measures on the
fibers of $\pr_1$ all have the same form \eqref{e:kost}. Moreover,
the $G$-invariance of $\mu$ reduces to the statement that on
$U\times\bR^{2\ell}$, $\mu$ has the form $\mu_U\times (dq_1\dotsm
dq_\ell dp_1\dotsm dp_\ell)$, where $\mu_U$ is a (uniquely
determined) positive Borel measure on $U$. Thus we can define
$\nu$ to be $(2\pi)^\ell\cdot\mu_U$ on $U$, and the rest of the
proof easily follows.

\section{Auxiliary computations}\label{s:B}

\subsection{Proof of convergence of the integral \eq{e:defka}}
Change the variables $z\mapsto a^{-1/6}z, b\mapsto a^{1/6}b$, and
then $a^{7/6}=u$; obtain
\[
\ka(H)=\frac{6}{7\pi}\int_0^\infty du \int_{-\infty}^{+\infty} db
\cdot N\left(u^{-2/7},-\frac{d^2}{d z^2}+\pare{z^2+b}^2\right).
\]
Introduce $A_0=-d^2/dz^2+\pare{z^2+b}^2$ and $A=-d^2/dz^2+z^4$, as
self-adjoint operators in $L_2(\bR)$. Since $A_0$ is positive
definite, there exist $c, c_1>0$ such that for all $b\ge -c_1$, $
A_0\ge 2^{-2/7}(A+c+b^2),$ hence
\[
\int_0^\infty du \int_{-c_1}^{+\infty} db \cdot
N\left(u^{-2/7},-\frac{d^2}{d
z^2}+\pare{z^2+b}^2\right)\le\int_0^\infty du
\int_{-c_1}^{+\infty} db \cdot N\left((u/2)^{-2/7},
A+b^2+c\right).\] The classical Weyl formula gives $N(\la, A)\sim
{\rm const}\, \la^{1/2+1/4}$, as $\la\to+\infty$, therefore we
have an upper bound via
\[
C_1\int_{-c_1}^{+\infty} db\int_0^{2(c+b^2)^{-7/2}} du\,
\left((u/2)^{-2/7}-(c+b^2)\right)^{3/4}.
\]
Change the variable $u\mapsto 2(c+b^2)^{-7/2}u$ to obtain the
product of a constant and two integrals:
\[
2C_1\int_{-c_1}^{+\infty} db\, (c+b^2)^{-7/2+3/4}\int_0^1du\,
(u^{-2/7}-1)^{3/4},
\]
which converge.

It remains to prove the convergence of the integral
$\int_{-\infty}^{-c_1}db\int_0^\infty du$. For each $b<0$, the
potential $(b+z^2)^2$ has two wells, at $z=\pm \sqrt{-b}$, and
grows as $|z|^4$ at infinity. The leading term of the Taylor
expansion at $\pm\sqrt{-b}$ is $8(-b)(z\mp \sqrt{-b})^2$,
therefore $A_0\ge c_2(-b)^{1/2},$ and
 \[
N(\la, A_0)\le C N(\la, -d^2/dz^2+(-b)z^2)\le C_1\la(-b)^{-1/2}.
\]
We conclude that the integrand below vanishes for $b<-C_3
u^{-4/7}$, and therefore
\[
\int_{-\infty}^{-c_1}db\int_0^\infty
du\,N\left(u^{-2/7},-\frac{d^2}{d z^2}+\pare{z^2+b}^2\right)\le
C_2 \int_0^{C_3} du\int_{-C_3u^{-4/7}}^{-c_1} db\,
u^{-2/7}(-b)^{-1/2}.
\]
The convergence of the integral on the RHS is straightforward.

\subsection{Proof of Lemma \ref{l:B}}

Fix $\eps<\eps'<1$. We have
\[
\meas\{P(x)+Q(x)\leq\la\}\leq\meas\{P(x)\leq\la\}=\meas\pare{\de_\la(\{P(x)\leq
1\})}=\la^{\abs{\ga}}\cdot \meas\{P(x)\leq 1\}.
\]
On the other hand,
\begin{eqnarray*}
\meas\{P(x)+Q(x)\leq\la\} &\geq & \meas\cur{
P(x)\leq\la-\la^{\eps'},\, Q(x)\leq\la^{\eps'} } \\
&=& \la^{\abs{\ga}}\cdot \meas\pare{
\de_\la^{-1}(\{P(x)\leq\la-\la^{\eps'},\, Q(x)\leq\la^{\eps'}\})}
\\
&\geq & \la^{\abs{\ga}}\cdot \meas\cur{ P(x)\leq
1-\la^{\eps'-1},\, Q(x)\leq C^{-1}\la^{\eps'-\eps} }.
\end{eqnarray*}
Since $Q$ is finite a.e., we have
\[
\cur{ P(x)\leq 1 } = \bigcup_{\la\to +\infty} \cur{ P(x)\leq
1-\la^{\eps'-1},\, Q(x)\leq C^{-1}\la^{\eps'-\eps} }
\]
(increasing union), whence
\[
\lim\limits_{\la\to +\infty} \meas\cur{ P(x)\leq
1-\la^{\eps'-1},\, Q(x)\leq C^{-1}\la^{\eps'-\eps} } =
\meas\{P(x)\leq 1\},
\]
completing the proof.

\subsection{Proof of Lemma \ref{l:weakdeg}}
Introduce the dilation $\dti_\la$ in the space $\bR^{2n}$ as
follows: \[\dti_\la(\xi,x)=\pare{ \la^{1/2}\xi, \la^{\ga_1}
x_1,\dotsc,\la^{\ga_n} x_n},\] and let $\de_\la$ on $\bR^n$ be
defined as in \S\ref{ss:Weyl}. Set
\[
\Psti^*(\xi,x)=\Psi^*(x)^2+\norm{\xi}^2, \quad
\Psti_0(\xi,x)=\Psi_0(x)^2+\norm{\xi}^2 \quad\text{and}\quad
\Psti(\xi,x)=\Psti^*(\xi,x)-\Psti_0(\xi,x).
\]
Then it is easy to see that
$\Psti^*\pare{\dti_\la(\xi,x)}=\la\cdot\Psti^*(\xi,x)$, and also,
using Lemma \ref{l:A}, we find that there exists $0<q<1$ with
$\Psti\pare{\dti_\la(\xi,x)}\leq\la^q\Psti(\xi,x)$ for all
$\la>1$. Now the proof of \eqref{e:dagger} is the same as the
proof of Corollary \ref{c:Weyl}.

\sbr

To prove that \eqref{e:genwint1} holds, we have to show that
\begin{equation}\label{e:ddag}
\meas\cur{\Psi^*(x)<M\cdot\Psi(x),\ \Psi^*(x)^2+\norm{\xi}^2\leq
C\la } = o\left( \meas\cur{\Psi^*(x)^2+\norm{\xi}^2\leq\la }
\right)
\end{equation}
as $\la\to +\infty$ for every fixed $C>0$, where $M=M(\la)$ is a
suitable function such that $M(\la)\to +\infty$ as $\la\to
+\infty$.

\sbr

We assume that $M$ grows slower than any power of $\la$.
Certainly, if $\Psi^*(x)^2+\norm{\xi}^2\leq C\la$, then
$\Psti_0(\xi,x)\leq C\la$; combining this with \eqref{e:dagger},
we see that to prove \eqref{e:ddag}, it suffices to check that
\[
\meas\cur{\Psi_0(x)<M\cdot\Psi(x),\ \Psti_0(\xi,x)\leq C\la } =
o\left( \meas\cur{\Psti_0(\xi,x)\leq\la } \right) \quad\text{as
}\la\to +\infty.
\]
Now, the dilation $\dti_\la$ scales the Lebesgue measure by a
factor of $\la^{\abs{\ga}+n/2}$. Dividing both sides of the last
equality by $\la^{\abs{\ga}+n/2}$, we see that we can rewrite it
as
\begin{equation}\label{e:dddag}
\meas\cur{\Psi_0(\de_\la x)<M\cdot\Psi(\de_\la x),\
\Psti_0(\xi,x)\leq C } = o\left( \meas\cur{\Psti_0(\xi,x)\leq 1 }
\right) \quad\text{as }\la\to +\infty.
\end{equation}
But $\meas\cur{\Psti_0(\xi,x)\leq 1}$ is a finite nonzero
constant; on the other hand, the set on the left is contained in
the set of points $(\xi,x)\in\bR^{2n}$ for which
$\Psti_0(\xi,x)\leq C$ and $\Psi_0(x)<M \la^{q-1} \Psi(x)$. Since
the set where $\Psti_0(\xi,x)\leq C$ has finite measure by
assumption, we can take the limit as $\la\to +\infty$; since
$M\cdot\la^{q-1}\to 0$ as $\la\to +\infty$ by construction, we see
that the limit of the LHS of \eqref{e:dddag} as $\la\to +\infty$
is at most the measure of the set where $\Psi_0(x)=0$, which is
zero because the potentials are polynomial.

\subsection{Proof of Lemma \ref{l:estimates}} Throughout the proof,
we use the notation $M$ for a function of $\la$ that grows slower
than any power of $\log\la$ as $\la\to +\infty$, such as
$M=\log\log\la$ (the precise formula is not important). Also, $C$
or $C'$ appearing in an inequality will always stand for a
constant which can be chosen to make the inequality work, and if
the inequality in question depends on an additional variable such
as $x$, it is to be understood that $C$ or $C'$ can be chosen
uniformly for all $x$. Finally, $h(\la)$ will stand for any
function of $\la$ such that $h(\la)\to 1$ as $\la\to +\infty$. We
define
\[
G_0(\la)=\meas\cur{x\in\bR^2 \st \abs{x_1},\abs{x_2}\geq M,\
\abs{x_1^k x_2^k}\leq\la^2\cdot h(\la) }.
\]
One immediately computes (using Fubini's theorem) that
$G_0(\la)\sim (8/k)\cdot\la^{2/k}\log\la$ as $\la\to +\infty$. On
the other hand, we have $\Phi^*(x)=\abs{x_1^k x_2^k}^{1/2}+R(x)$,
where
\[
R(x)=\sum_{j>0 \text{ or } l>0} \left| \frac{k!}{(k-j)!}\cdot
\frac{k!}{(k-j)!}\cdot x_1^{k-j} x_2^{k-l}\right|^{1/2}.
\]
If $\abs{x_1}\geq M$ and $\abs{x_2}\geq M$, then $R(x)\leq
(C/M)\cdot\abs{x_1^k x_2^k}$. Therefore, if $\abs{x_1}\geq M$,
$\abs{x_2}\geq M$ and $\abs{x_1^k
x_2^k}\leq\la^2\cdot\pare{1-\frac{C}{M}}$, then
$\Phi^*(x)\leq\la$. Taking $h(\la)=1-C/M$ and applying the
statement about the asymptotics of $G_0(\la)$, we see that
$G_1(\la)=\meas\{\Phi^*\leq\la\}$ is bounded below by a function
which has the required asymptotics as $\la\to +\infty$. On the
other hand, we have
\[
\cur{\Phi^*(x)\leq\la}\subseteq\cur{\abs{x_1},\abs{x_2}\geq
M,\abs{x_1^k x_2^k}\leq\la^2} \,\cup\, \cur{\abs{x_1}\leq\la,\
\Phi^*(x)\leq\la} \,\cup\, \cur{\abs{x_2}\leq\la,\
\Phi^*(x)\leq\la}.
\]
As above, the measure of the first set on the RHS has the required
asymptotics, while the measure of the second set is bounded by
\[
\meas\cur{\abs{x_1}\leq M,\ \abs{x_2}\leq\la^{2/k} } =
2M\la^{2/k},
\]
which grows slower than $\la^{2/k}\log\la$. By symmetry, the same
statement holds for the third set. This provides the desired upper
bound and concludes the proof of the formula $G_1(\la)\sim
(8/k)\cdot\la^{2/k}\log\la$.

\sbr

Next we study the asymptotics of the function
$G_2(\la)\sim\meas\{\Psi^*\leq\la\}$. We will first prove that the
asymptotics of $G_2(\la)$ is controlled by the measure of the set
where $\Psi^*(x)\leq\la$ and $\Psi^*(x)$ is dominated by the first
term, in the sense that $\abs{x_1^k x_2^k}^{1/2}\geq M\cdot
\abs{x_1^{k-j} x_2^{k-l}}^{1/(2+j+l)}$ whenever $j>0$ or $l>0$. To
this end, we fix $(j,l)\neq (0,0)$ and consider the measure of the
set of points $(x_1,x_2)\in\bR^2$ such that $\Psi^*(x)\leq\la$ and
\begin{equation}\label{e:ineq}
\abs{x_1^k x_2^k}^{1/2}\leq M\cdot \abs{x_1^{k-j}
x_2^{k-l}}^{1/(2+j+l)}.
\end{equation}
It suffices to show that this measure grows slower than
$\la^{2/k}\log\la$, since $G_2(\la)\geq G_1(\la)$. By symmetry, we
may assume that $x_1,x_2>0$. We may rewrite \eqref{e:ineq} as
\begin{equation}\label{e:ineq2}
x_1\leq M^{\frac{2(2+j+l)}{kj+kl+2j}}\cdot
x_2^{-\frac{kj+kl+2l}{kj+kl+2j}}.
\end{equation}
There are two cases to consider.

\sbr

\noindent
1) We have $j=l>0$. Set $M'=M^{\frac{2(2+j+l)}{kj+kl+2j}}$. Thus
\eqref{e:ineq2} becomes $x_1\leq M'/x_2$. On the other hand, the
condition $\Psi^*(x)\leq\la$ forces, in particular,
$x_1^k\leq\la^{k+2}$ and $x_2^k\leq\la^{k+2}$. One easily checks
that the measure of the set where these two inequalities are
satisfied and where $x_1 x_2\leq M'$ grows as $C\cdot
M'\cdot\log\la\ll \la^{2/k}\log\la$.

\sbr

\noindent
2) We have $j\neq l$. By symmetry, we may assume that $j>l$. So
\eqref{e:ineq2} becomes $x_1\leq M'\leq x_2^{-\kappa}$, where
\[
0<\kappa=\frac{kj+kl+2l}{kj+kl+2j}=1-2\cdot\frac{j-l}{kj+kl+2j}<1.
\]
One easily checks that $\kappa\geq\frac{k}{k+2}$ (this minimum
value is achieved when $j=k$, $l=0$). Thus
$0<1-\kappa\leq\frac{2}{k+2}$. Now $\Psi^*(x)\leq\la$ implies, in
particular, $x_2\leq\la^{(k+2)/k}$. One easily bounds the measure
of the set where $0\leq x_2\leq 1$ and $x_1\leq M'\cdot
x_2^{-\kappa}$ by $C M'\ll\la^{2/k}\log\la$, and the measure of
the set where $1\leq x_2\leq\la^{(k+2)/k}$ and $x_1\leq M'\cdot
x_2^{-\kappa}$ by $C M' \la^{2/k}\ll\la^{2/k}\log\la$.

\mbr

To complete the proof of the lemma, we have to estimate the
measure of the set where $\Psi^*(x)\leq\la$ and $\abs{x_1^k
x_2^k}^{1/2}\geq M\cdot \abs{x_1^{k-j} x_2^{k-l}}^{1/(2+j+l)}$
whenever $j>0$ or $l>0$. Again, by symmetry, we may assume that
$0\leq x_1\leq x_2$. Clearly, the set where also $x_2\leq 1$ does
not contribute to the asymptotics. Therefore we only need to
consider the set where $x_2\geq 1$, $x_2\geq x_1\geq M'\cdot
x_2^{-\kappa}$ and $\Psi^*(x)\leq\la$. Here, $M'$ and $\kappa$ are
defined as above, and $(j,l)$ run over all pairs with $0\leq
j,l\leq k$. But since $\kappa\geq \frac{k}{k+2}$, all the
inequalities $x_1\geq M'\cdot x_2^{-\kappa}$ reduce to $x_1\geq
M'\cdot x_2^{-k/(k+2)}$ if $x_2\geq 1$. Under the assumptions
\begin{equation}\label{e:a1}
x_2\geq 1,\quad x_2\geq x_1\geq M'\cdot x_2^{-\frac{k}{k+2}},
\end{equation}
it is clear that the following implications hold:
\[
x_1^k x_2^k\leq \pare{1-C/M}\cdot\la^2 \quad\Longrightarrow\quad
\Psi^*(x)\leq\la \quad\Longrightarrow\quad x_1^k x_2^k\leq\la^2.
\]
Finally, one computes that for any function $h(\la)$ such that
$h(\la)\to 1$ as $\la\to +\infty$, the measure of the set where
\eqref{e:a1} holds and where $x_1^k x_2^k\leq h(\la)\cdot\la^2$ is
asymptotically equal to $\frac{k+1}{k}\cdot\la^{2/k}\cdot\log\la$.
Multiplying this by $8$ (due to the fact that we only considered
the set where $x_2\geq x_1\geq 0$), we see that the proof of Lemma
\ref{l:estimates} is complete.

\subsection{Proof of Lemma \ref{l:limit}}
Note that the measure $\mu=\mu_\Om$ has a ``trivial
component''---namely, $\mu$ is naturally a product of a measure
$\mu'$ on $\bR^2$ and a measure $\mu''$ on $\bR^{(k+1)^2}$, given
by $\mu'=(2\pi)^{-2}\cdot d\xi_1 d\xi_2$ and $\mu''=\phi_*(dx_1
dx_2)$, where $\phi:\bR^2\to\bR^{(k+1)^2}$ is defined by
$\phi(x_1,x_2)=(x_1^k x_2^k,\dotsc,1)$. It is clear that we can
apply our scaling constructions to the measures $\mu'$ and $\mu''$
separately, obtaining measures $\mu'_\la$ and $\mu''_\la$ for all
$\la>0$. Trivially, there is a weak limit
\[
\mu'_0=\lim\limits_{\la\to +\infty} \la^{-2}\mu'_\la = \mu' =
(2\pi)^{-2}\cdot d\xi_1 d\xi_2.
\]
Thus, we only need to concentrate on the measures $\mu''_\la$ on
$\bR^{(k+1)^2}$. The proof of Lemma \ref{l:limit} will be complete
if we show that there is a weak limit
\begin{equation}\label{e:mupp}
\mu''_0=\lim\limits_{\la\to +\infty} \la^{-1/k}
(\log\la)^{-1}\cdot \mu''_\la = \psi_* \pare{ \frac{2}{k^2}\cdot
\abs{y}^{\frac{1-k}{k}} dy },
\end{equation}
where $\psi:\bR\to\bR^{(k+1)^2}$ is given by
$\psi(y)=(y,0,\dotsc,0)$. Let us define $\mu''_0$ by the RHS of
\eqref{e:mupp}. The idea is to apply Proposition \ref{p:limborel}.
Since the function $\abs{y}^{(1-k)/k}$ is integrable near $0$, we
see that $\mu''_0$ is a c-finite Borel measure on $\bR^{(k+1)^2}$.
Let us first fix real numbers $R_{ij}>0$ ($0\leq i,j\leq k$) and
consider for every $\la>0$ the set
\[
A_\la=\cur{ (x_1,x_2) \st x_1,x_2\geq 0,\, \la^{-1} x_1^i
x_2^j\leq R_{ij} \ \forall\, 0\leq i,j\leq k }.
\]
It is easy to see that when $u>0$ is large enough (depending only
on the $R_{ij}$'s, but not on $\la$), we have
\[
A_\la\cap \cur{x_1\geq u,\, x_2\geq u} = \cur{x_1\geq u,\, x_2\geq
u,\, x_1^k x_2^k \leq \la R_{kk} },
\]
and then one computes (for fixed $R_{ij}$ and $u$)
\[
\meas (A_\la) \,\sim\, \meas\left(A_\la\cap \cur{x_1\geq u,\,
x_2\geq u} \right) \,\sim\, \frac{R_{kk}^{1/k}}{k} \cdot \la^{1/k}
\log\la \quad\text{as } \la\to +\infty.
\]
By symmetry, if
\[
B_\la = \cur{ (x_1,x_2)\in\bR^2 \st \la^{-1} \abs{x_1^i x_2^j}\leq
R_{ij} \ \forall\, 0\leq i,j\leq k },
\]
then
\begin{equation}\label{e:computation}
\meas (B_\la) \,\sim\, 4\cdot \frac{R_{kk}^{1/k}}{k} \cdot
\la^{1/k} \log\la \quad\text{as } \la\to +\infty.
\end{equation}
In particular, the answer does not depend on $R_{ij}$ if $i<k$ or
$j<k$. We now apply Proposition \ref{p:limborel} with $E=\{0\}$,
which is certainly countable. So let
$a=\pare{a_{ij}}\in\bR^{(k+1)^2}$ and
$b=\pare{b_{ij}}\in\bR^{(k+1)^2}$, with $a_{ij}<b_{ij}$ and
$a_{ij}\neq 0$, $b_{ij}\neq 0$ for all $0\leq i,j\leq k$. Two
situations are possible.

\sbr

\noindent
1) We have $0\not\in [a_{ij},b_{ij})$ for some $(i,j)$ such that
either $i<k$ or $j<k$. Then there exists $\de>0$ such that
$[a_{ij},b_{ij})$ does not intersect $[-\de,\de]$ and is contained
in $[\de^{-1},\de^{-1}]$. We apply the observation of the previous
paragraph with $R_{ij}=\de$ and $R_{ij}=\de^{-1}$, and conclude
that
\[
\lim\limits_{\la\to +\infty} \la^{-1/k}\log\la\cdot
\mu_\la\pare{[a,b)} = 0 = \mu''_0\pare{[a,b)}.
\]

\sbr

\noindent
2) We have $0\in [a_{ij},b_{ij})$ for all $(i,j)$ such that either
$i<k$ or $j<k$. Then, as before, there exists $\de>0$ such that
$[\de,\de]\subset[a_{ij},b_{ij})\subset [\de^{-1},\de^{-1}]$ for
all such $(i,j)$. Now if $a_{kk}=-b_{kk}$, then the computation
\eqref{e:computation} immediately implies that
\[
\lim\limits_{\la\to +\infty} \la^{-1/k}\log\la\cdot
\mu_\la\pare{[a,b)} = 4\cdot\frac{a_{kk}^{1/k}}{k} =
\mu''_0\pare{[a,b)}.
\]
The general case reduces to this one by additivity and symmetry,
noting that, for example, if $b_{kk}>a_{kk}>0$, then
\[
[a_{kk},b_{kk}) \cup [-b_{kk},-a_{kk}) = [-b_{kk},b_{kk})
\setminus [-a_{kk}, a_{kk}),
\]
and if $b_{kk}>-a_{kk}>0$, then
\[
[a_{kk},b_{kk} = [a_{kk},-a_{kk})\cup [-a_{kk},b_{kk}).
\]


\begin{thebibliography}{99}

\bibitem[BCD]{reps} P.~Bernat, C.~Conze, M.~Duflo,
N.~L\'evy-Nahas, M.~Rais, P.~Renouard and M.~Vergne,
``Repr\'esentations des Groupes de Lie R\'esolubles.''
Monographies de la Soc. Math. de France \textbf{4} (1972).

\bibitem[Bo]{Bon} P.~Bonnet, {\em Param\'etrisation du dual d'une
alg\`ebre de Lie nilpotente}, Ann. Inst. Fourier (Grenoble)
\textbf{38} (1988), no.~3, 169--197.

\bibitem[CdV]{CdV} Y.~Colin de Verdi\`ere, {\em L'asymptotique de Weyl
pour les bouteilles magnetiques}, Comm. Math. Phys. \textbf{105}
(1986), 327--335.

\bibitem[CFKS]{CFKS}
H.L.~Cycon, R.G.~Froese, W.~Kirsch and B.~Simon, ``Schr\"odinger
operators with applications to quantum mechanics and global
geometry.''  Berlin, New York, Heidelberg, London, Paris, Tokyo:
Springer-Verlag (1985).

\bibitem[Fe]{Fe} C.L.~Fefferman, {\em The uncertainty principle},
Bull. Amer. Math. Soc. \textbf{9} (1983), 129--206.

\bibitem[Gu]{G}
D.~Gurarie, {\em Non-classical eigenvalue asymptotics for
operators of Schr\"odinger type}, Bull. Am. Math. Soc. \textbf{15}
(2) (1986), 233--237.

\bibitem[HM]{HM} B.~Helffer and A.~Mohamed, {\em Caract\'erization
du spectre essentiel de l'op\'erateur de Schr\"odinger avec une
champ magn\'etique}, Ann. Inst. Fourier (Grenoble) \textbf{38}
(1988), no.~2, 95--112.

\bibitem[HN]{HN} B.~Helffer and J.~Nourrigat, ``Hypoellipticit\'e
Maximale pour des Op\'erateurs Polynomes de Champs de Vecteurs.''
Progress in Math., vol. 58, Birkh\"auser, Boston, 1985.

\bibitem[H]{H}
L.~H\"{o}rmander, ``The analysis of differential operators. 3.''
Berlin, New York, Heidelberg: Springer-Verlag (1985).

\bibitem[Ki]{Ki} A.A.~Kirillov, {\em Unitary representations of
nilpotent Lie groups}, Uspehi Mat. Nauk \textbf{17} (1962), no. 4
(106), 57--110.

\bibitem[Iv]{I} V.~Ivri\v{i}, {\em Estimate for the number of negative eigenvalues
of the Schr\"odinger operator with intense field}, Journ\'{e}es
\'{E}quations aux D\'{e}riv\'{e}es partielles de Saint-
Jean-de-Monts, Soc. Math. France (1987).

\bibitem[L1]{L1} S.Z.~Levendorski\v{i}, {\em Non-classical spectral asymptotics},
Russian Math. Surveys \textbf{43} (1988), no.~1, 123--157.

\bibitem[L2]{L2} S.Z.~Levendorski\v{i}, ``Asymptotic distribution of eigenvalues of
differential operators.'' Dordrecht: Kluwer Academic Publishers
(1990).

\bibitem[L3]{L3}  S.Z.~Levendorski\v{i},
``Degenerate elliptic equations.'' Dordrecht: Kluwer Academic
Publishers (1993).

\bibitem[L4]{L4} S.Z.~Levendorski\v{i}, {\em Spectral properties of Schr\"odinger
operators with irregular magnetic potentials, for a spin
$\frac{1}{2}$ particle}, J. Math. Anal. Appl. \textbf{216} (1997),
no.~1, 48--68.

\bibitem[LMN]{LMN} P.~Levy-Bruhl, A.~Mohamed and J.~Nourrigat,
{\em Spectral theory and representations of nilpotent groups},
Bull. Amer. Math. Soc. \textbf{26} (1992), no.~2, 299--303.

\bibitem[Ma1]{Ma1} D.~Manchon, {\em
Formule de Weyl pour les groupes de Lie nilpotents}, J. Reine
Angew. Math. \textbf{418} (1991), 77--129.


\bibitem[Ma2]{Ma2} D.~Manchon, {\em Weyl symbolic calculus on any Lie group},
Acta Appl. Math. \textbf{30} (1993), no. 2, 159--186.

\bibitem[Ma3]{Ma3} D.~Manchon, {\em Op\'erateurs pseudodiff\'erentiels et repr\'esentations unitaires des
groupes de Lie}, Bull. Soc. Math. France \textbf{123} (1995), no.
1, 117--138.

\bibitem[MN]{MN} A.~Mohamed and J.~Nourrigat, {\em Encadrement du
$N(\la)$ pour des op\'erateurs de Schr\"odinger avec champ
magn\'etique}, J. Math. Pures Appl. (9) \textbf{70} (1991),
87--99.

\bibitem[Ni1]{Ni1} N.~Nilsson, {\em Asymptotic estimates for spectral functions
connected with hypoelliptic differential operators}, Ark. Mat.
\textbf{5} (1965), 527--540.

\bibitem[Ni2]{Ni2} N.~Nilsson, {\em Some growth and ramification
properties of certain integrals on algebraic manifolds}, Ark. Mat.
\textbf{5} (1965), 463--476.

\bibitem[Pu1]{long} L.~Pukanszky, {\em Unitary representations of solvable Lie
groups}, Ann. Sci. \'Ecole Norm. Sup. (4) \textbf{4} (1971),
457--608.

\bibitem[Pu2]{Pu} L.~Pukanszky,
{\em On the theory of exponential groups}, Trans. Amer. Math. Soc.
\textbf{126} (1967), 487--507.

\bibitem[Ro]{R} D.~Robert, {\em Comportement asymptotique des valeurs propres
d'op\'erateurs de type Schr\"odinger \`a potentiel
``d\'eg\'en\'er\'e''}, J. Math. Pures Appl. \textbf{61} (1982),
275--300.

\bibitem[Roz]{Roz} G.V.~Rozenbljum, {\em Asymptotic behavior of the eigenvalues of the
Schrdinger operator}, Mat. Sb. (N.S.) \textbf{93} (135) (1974),
347--367, 487.

\bibitem[RSS]{RSS} G.V.~Rozenbljum, M.Z.~Solomyak and M.A.~Shubin, ``Spectral theory of
differential operators.'' Contemporary problems of mathematics
(Itogi Nauki i Tekhniki VINITI), v.64. Moscow: VINITI (1989).

\bibitem[Sim]{Sim} B.~Simon, {\em Nonclassical eigenvalue asymptotics}, J. Funct.
Anal. \textbf{53} (1983), no.~1, 84--98.

\bibitem[Sol]{S} M.Z.~Solomyak, {\em Asymptotics of the spectrum of the Schr\"odinger
operator with non-regular homogeneous potential}, Math. USSR
Sbornik, \textbf{55} (1986), no.~1, 19--37.

\bibitem[Tam]{T} H.~Tamura, {\em Asymptotic distribution of eigenvalues for Schr\"odinger
operators with magnetic fields}, Nagoya Math. J. \textbf{105}
(1987), 49--69.

\bibitem[TS]{TS} V.N.~Tulovski\v{i} and M.A.~Shubin, {\em
The asymptotic distribution of the eigenvalues of
pseudodifferential operators in $\bR^n$}, Mat. Sb. (N.S.)
\textbf{92}(134) (1973), 571--588, 648.

\end{thebibliography}
\end{document}